\documentclass[10pt,a4paper]{article}
\usepackage{amsthm}
\usepackage{cite}
\usepackage{mathtools}
\usepackage{enumitem}
\usepackage{hyperref}
\usepackage{graphicx}
\usepackage{amssymb,array}
\usepackage[table,dvipsnames]{xcolor}
\usepackage{longtable}
\usepackage{pdflscape}
\usepackage[all]{xy}
\usepackage{tikz}
\usetikzlibrary{arrows.meta}
\usepackage{wrapfig}
\usepackage{bbm}
\usepackage{wrapfig}
\usepackage[symbol]{footmisc}

\usepackage{caption,subcaption}
\usepackage{adjustbox}

\usetikzlibrary{positioning}
\usepackage{blkarray}
\usepackage{multirow}
\usetikzlibrary{shapes}

\usepackage{geometry}
\usepackage{float}

\theoremstyle{plain}
\newtheorem{thm}{Theorem}[section]

\newtheorem{crt}[thm]{Criterion}

\theoremstyle{definition}
\newtheorem{dfn}[thm]{Definition}

\newtheorem{rem}[thm]{Remark}

\newtheorem{algorithm}[thm]{Algorithm}
\newtheorem{notation}[thm]{Notation}

\newcommand{\Bmu}{\mbox{$\raisebox{-0.59ex}
  {$l$}\hspace{-0.18em}\mu\hspace{-0.88em}\raisebox{-0.98ex}{\scalebox{2}
  {$\color{white}.$}}\hspace{-0.416em}\raisebox{+0.88ex}
  {$\color{white}.$}\hspace{0.46em}$}{}}

\newcommand{\MM}{{\mathbb{M}}}
\newcommand{\cO}{\mathcal{O}}
\newcommand{\QQ}{{\mathbb{Q}}}
\newcommand{\RR}{{\mathbb{R}}}
\newcommand{\PP}{{\mathbb{P}}}

\newcommand{\ZZ}{{\mathbb{Z}}}

\newcommand{\FF}{\mathbb{F}}
\newcommand{\CC}{\mathbb{C}}
\newcommand{\GG}{\mathbb{G}}
\newcommand{\LL}{\mathbb{L}}
\newcommand{\NN}{{\mathbb{N}}}
\newcommand{\TT}{\mathbb{T}}

\renewcommand{\emptyset}{\varnothing}
\renewcommand{\tilde}{\widetilde}

\newcommand{\RelInt}{\operatorname{RelInt}}
\newcommand{\codim}{\operatorname{codim}}
\newcommand{\Pf}{\operatorname{Pf}}
\newcommand{\Cl}{\operatorname{Cl}}

\newcommand{\Spec}{\operatorname{Spec}}

\makeatletter
\newcommand{\subjclass}[2][2022]{%
  \let\@oldtitle\@title%
  \gdef\@title{\@oldtitle\footnotetext{#1 \emph{Mathematics subject classification.} #2}}%
}

\usepackage{authblk}

\begin{document}


\title{$\QQ$-Fano threefolds and Laurent Inversion}
\author{Liana Heuberger\thanks{Universit\'e d'Angers, l.heuberger@univ-angers.fr.}}
\date{}
\maketitle
\begin{abstract}

We construct families of non-toric $\QQ$-factorial terminal Fano ($\QQ$-Fano) threefolds of codimension $\geq 20$ corresponding to 54 mutation classes of rigid maximally mutable Laurent polynomials. From the point of view of mirror symmetry, they are the highest codimension (non-toric) $\QQ$-Fano varieties for which we can currently establish the Fano/Landau--Ginzburg correspondence. We construct 46 additional $\QQ$-Fano threefolds with codimensions of new examples ranging between 19 and 5. Some of these varieties will be presented as toric complete intersections, and others as Pfaffian varieties. 
\end{abstract}
\section{Introduction}

In this paper we use mirror symmetry methods to construct a number of new deformation families of three-dimensional Fano varieties with terminal orbifold singularities. Our main technique, Laurent inversion \cite{LaurentInversion}, is an inverse of the celebrated Givental/Hori-Vafa construction that associates a Laurent polynomial mirror to a Fano toric complete intersection. We construct our examples by applying Laurent inversion to specific class of Laurent polynomials. The resulting deformation families of $\QQ$-Fano threefolds are complete intersections or Pfaffians, of low codimension, in singular toric varieties. However, the embedding codimension (see below) of these examples is high, so they are well beyond the reach of classical constructions. Since the ambient toric varieties are highly singular, these examples also fall outside the range where Laurent inversion was previously known to work. Thus this paper has several implications:
\begin{itemize}
\setlength\itemsep{0.25em}
\item it systematically explores part of the ``landscape" of possible $\QQ$-Fano threefolds that so far has been out of reach;
\item it is one of the first occasions where mirror symmetry has been used to construct previously unknown algebraic varieties;
\item it gives further evidence for a surprising phenomenon: that terminal singularities, which arise so naturally in birational geometry, are also singled out by mirror symmetry.
\end{itemize} 

\subsection*{Outline}
$\QQ$-factorial terminal threefolds appear naturally in the minimal model program. Among these, we can hope to classify Fano varieties, as they come in a finite number of families. For simplicity, we only consider orbifold terminal quotient singularities, i.e. cyclic quotient singularities of the form $\frac{1}{r}(1,a,-a)$ with $(a,r)=1$. Possible Hilbert series of these objects have been inventoried in the Graded Ring Database (GRDB for short), yielding $54610$ cases in the ``Fano 3-folds" list \cite{GRDB-fano3} (see \cite{BKFano3,ABR02,BS07a,BS07b,suzuki2004Fano} for the theorems on which it is based and \cite{rep} for the data repository). Any polarised variety can be embedded inside a weighted projective space, and the GRDB also predicts the number of generators of such an embedding. The \textit{Fano index} of a variety $X$ is the largest positive index $r>0$ such that $-K_X=rA$ in the class group $\Cl X$. Throughout this paper we refer to the codimension of $X$ with respect to the predicted $A$-embedding as its \textit{GRDB-codimension}. We do not know how many entries in the ``Fano 3-folds" list correspond to existing varieties, especially since constructions using graded ring methods have had limited success in high codimension (i.e. $\geq 5$). This paper shows that a previously inaccessible part of the database is in fact constructible in a practical way. 

The examples are situated within the more general framework of providing evidence for the Fano/Landau--Ginzburg correspondence as formulated in \cite{Coates_2015}. This can be represented as a conjectural bijective correspondence between the two sets:

\[ \left\lbrace \begin{array}{c}
\textup{qG-deformation equivalence classes of}\\
\QQ\textup{-Fano threefolds of class TG}\\ 
\end{array} \right\rbrace \stackrel{\rm F/LG}{\longleftrightarrow} \left\lbrace
\begin{array}{c}
\textup{mutation classes of}\\
\textup{rigid MMLPs} 
\end{array}
\right\rbrace, \] 
where MMLP stands for \textit{maximally mutable Laurent polynomial} (see Section \ref{LPintro}). A $qG$ or \textit{$\QQ$-Gorestein} deformation of a $\QQ$-Gorenstein variety $X$ is a flat family $\pi \colon \mathcal{X} \to B$ over a scheme $B$ such that there exist $r\in \NN$ and a Cartier divisor $L$ on $\mathcal{X}$ with $L|_F=-r\cdot K_F$ for any fibre $F$ of  $\pi$ and $b\in B$ such that $X\simeq\pi^{-1}(b)$. If $X$ is a Fano variety, the Hilbert series is locally constant on the fibres of such deformations (see \cite{OrbidP} for a more detailed discussion). 

\begin{dfn}
Two $\QQ$-Fano varieties $X$ and $X'$ are \textit{qG-deformation equivalent} if there exists a qG deformation $\pi\colon\mathcal{X}\to B$ with terminal fibres such that $B$ is connected and there exist $b,b'\in B$ with $X\simeq \pi^{-1}(b)$ and $X'\simeq \pi^{-1}(b')$.
\end{dfn}

The F/LG correspondence can be understood more geometrically in terms of \textit{toric degenerations}:  the Newton polytope $P$ of a rigid MMLP $f$ spans a fan from which we build a toric Fano variety $X_P$. This is in general far from being terminal or $\QQ$-factorial, yet we expect $X_P$ to qG-deform to a  $\QQ$-Fano variety $X$. Working with toric Fano varieties is a crucial advantage, as they can be described via lattice polytopes, of which $674688$ appear in a list of toric \textit{canonical} varieties in the GRDB \cite{GRDB-toric3, Kasprzyk2010, repC3}. We assume that our varieties admit a toric degeneration, and call this property TG (from \textit{toric generization}). We further assume that such a degeneration contains a canonical toric variety in its mutation class. 

We say a $\QQ$-Fano $X$ and a Laurent polynomial $f$ are mirror to each other (via the F/LG correspondence) if the regularised quantum period of $X$ and the classical period of $f$ coincide (see \cite{MSandF} for the general picture). We refer to the coefficients of the classical period when expanded as a series as the \textit{period sequence} of the Laurent polynomial. 

In upcoming work, Coates, Kasprzyk and Pitton classify period sequences of rigid MMLPs supported on canonical polytopes. For our first result, we extract from the GRDB the Hilbert series of the candidate $\QQ$-Fano varieties of codimension $\geq 20$. We find all the canonical polytopes for each Hilbert series and look for rigid MMLPs supported on them. At the time of writing, there are 54 unique period sequences that can be obtained this way \cite{repMMLP}. 

\begin{thm}
\label{systemThm}
We construct all 54 $\QQ$-Fano families with GRDB-codimension at least 20 that correspond to a period sequence as described above. They are listed in Table \ref{tab1} together with their embedding data inside a toric variety. We verify that the examples constructed as complete intersections correspond under mirror symmetry to the Laurent polynomials in the ``Rigid MMLP" column. 
\end{thm}

In Theorem \ref{systemThm} we construct a family for every period sequence arising from a set of rigid MMLPs, thus providing strong evidence for the Fano/LG correspondence. However, we are not defining a set-theoretic function, in the sense that in principle there could exist more than one such family. We have no knowledge of this phenomenon actually occurring, neither in our work or in the literature.

There are many benefits to working in dimension three, from visualising polytopes to the specific form of singularities, which make it similar to the surface case. Nevertheless, there is an increased difficulty from the point of view of deformation theory (see \cite{petracci2019deformations} for a survey on the topic). In particular, we are missing two key ingredients: a concrete result in the style of \cite[Proposition 2.7]{SingContent}, which is a recipe for determining which toric surface singularities are smoothable, and \cite[Lemma 6]{OrbidP}, a useful local-to-global glueing result that is central when obtaining classifications. There have been attempts at settling these issues (e.g. \cite[Conjectures A\&B]{corti2020mirror} provide a conjectural answer to the first question in the Gorenstein case), yet it is apparent that they are relevant even in mildly singular cases (see \cite[Theorem 1.2]{petracci2019some} for an example where expected local-to-global glueing fails). Fortunately, the strategy we employ in this paper bypasses most of these subtle discussions by constructing the deformations explicitly. 

We are able to obtain our constructions precisely because we shift perspective from the high codimension embedding inside a projective space to a codimension $\leq 2$ complete intersection or a codimension 3 Pfaffian inside a higher-rank toric variety, as we explain below.

The method we use to put this in practice is  \textit{Laurent inversion}: in \cite{LaurentInversion}, Coates, Kasprzyk and Prince reverse-engineer of the so-called ``Przjyalkowski trick" \cite{prz2007}. Their algorithm is recalled in Section \ref{LIalg}: starting from a Laurent polynomial, it produces an embedded Fano variety.  In \cite{doran_harder_2016}, Doran and Harder obtain similar constructions in the complete intersection case, however the present article uses the language and notation in \cite{LaurentInversion}---the correspondence between the two is discussed in \cite[Section 12]{LaurentInversion}. 

Instead of considering its entire deformation family, the authors of \cite{LaurentInversion} construct a toric embedding of a Fano toric variety $X_P$ into another toric variety $Y$, either as a complete intersection or as a Pfaffian variety. Using the Cox coordinates on $Y$, one deforms the equations that give $X_P$ to obtain a less singular variety $X$. In general, $Y$ is not an orbifold (neither is $X_P$), and moreover its non-orbifold locus need not be isolated. However, if chosen conveniently, $Y$ is large enough so that $X$ completely avoids this locus and only has $\QQ$-factorial terminal singularities. Very explicit such computations are found in Sections \ref{2Dexamples} and \ref{3DPf}. Crucially, the codimension of $X$ inside $Y$ in the examples of this paper ranges from 1 to 3, dramatically decreasing the GRDB-codimension and thus allowing for more control on our variety.

One of our aims is to see how far we are able to go with these constructions, and to determine this we dive deeper into the GRDB. The key object in the Laurent inversion algorithm is the \textit{shape}, a smooth toric variety which encodes a decomposition of a given Laurent polynomial into summands. Its geometry rules many of the properties of $Y$ and the embedding of $X_P$ inside it. The dimension of the shape can in principle range from 1 to $\dim X_P$. For all but two varieties in Table \ref{tab1}, this shape is two-dimensional. We further explored the database in order to determine if the predominance of $2$D shapes is characteristic to varieties of high codimension.

\begin{thm}
\label{exploreThm}
In Table \ref{tab2} we construct $46$ $\QQ$-Fano threefolds with GRDB-codimensions of new examples ranging from $19$ to $5$, together with their singularity baskets and embedding data. We verify that the examples presented as complete intersections correspond under mirror symmetry to the Laurent polynomials in the ``Rigid MMLP" column. 
\end{thm}

\begin{rem} Eleven examples in Table \ref{tab2} use $3$D shapes in their construction, and we discuss the most complicated case in Section \ref{3DPf}.\end{rem}

Theorem \ref{exploreThm} extracts more information on the structure of the deformation spaces involved than Theorem \ref{systemThm}, in the following sense.

Mutating a polynomial implies mutating its Newton polytope, a combinatorial operation described in \cite[\S 3]{Akhtar_2012}. Geometrically, this operation produces a pencil deforming a toric Fano variety into another. This pencil was introduced in \cite{Ilten_2012}, was formulated in the setting of toric pairs in \cite[Theorem 1.3]{petracci2020homogeneous}, and the deformation was shown to be $\QQ$-Gorenstein in the surface case in \cite[Lemma 7]{OrbidP}. 

The notion of mutation for polynomials is finer than its combinatorial counterpart, which can be interpreted in the following geometric way from the point of view of mirror symmetry. Mutation-equivalent \textit{polytopes} determine toric varieties with the same Hilbert series. If both toric varieties deform to $\QQ$-Fano varieties, their Hilbert series will also coincide. Note that the converse is not true: two polytopes with the same Hilbert series are not necessarily mutation-equivalent. An easy example are the polytopes of $\FF_1$ and $\PP^1\times \PP^1$, which are the two smooth toric del Pezzo surfaces of degree $8$. 

A refinement of the statement above is that mutation-equivalent \textit{polynomials} correspond to $\QQ$-Fano varieties with the same Gromov--Witten invariants. Since the latter are deformation-invariant, this is a rephrasing of the F/LG correspondence as stated above: when we fix the mutation class of a polynomial, we are singling out a component of a deformation space. 

The Laurent inversion process starts from a Laurent polynomial $f$ and creates an embedding of a toric Fano variety $X_f$. In the case of toric complete intersections, we are able to check that the partial smoothing of $X_f$ to a $\QQ$-Fano inside this embedding is mirror to $f$: it has the same regularised quantum period as the classical period of $f$. 

Furthermore, if two different rigid MMLPs $f$ and $g$ with different period sequences are supported on the same polytope and they deform to $\QQ$-Fano threefolds $X_f$ and $X_g$, then these two varieties are not deformation-equivalent. In other words, from the point of view of mirror symmetry we have found different Fano varieties. Some such examples appear in Table \ref{tab2}, and crucially do not occur in the surface case (see \cite[Lemma 6]{OrbidP}) nor in the high codimension examples in Table \ref{tab1}. Combinatorial mutations alone cannot detect this phenomenon.

\paragraph*{Input of the construction.} 
\begin{enumerate}
    \item Parse the list of canonical polytopes in GRDB which have the same Hilbert series as a (conjectural) $\QQ$-Fano threefold. 
    \item Find rigid MMLPs on the polytopes. In particular, the highest codimension in which these exist is 23. 
    \item Some of these polynomials have the same period sequence, i.e. conjecturally they deform in the same family of $\QQ$-Fano threefolds. We treat these redundancies as follows: 
    \begin{enumerate}
    \item Given a fixed period sequence, if there is a terminal, $\QQ$-factorial toric variety supporting one of its polynomial representatives, we skip this case as this variety exists in the toric canonical 3-fold database. See Section \ref{Whynot} for a precise description of what is discarded in the systematic codimension $\geq 20$ case.
     \item Otherwise, we choose a single polynomial representative of this class.
     \end{enumerate}    
 \end{enumerate}
\paragraph*{Sketch of the construction.} 
\begin{enumerate}
  \setcounter{enumi}{3}

    \item For such a Laurent polynomial, use Laurent inversion to construct the $\QQ$-Fano threefold. This is either a complete intersection in a toric variety, or a codimension $3$ Pfaffian variety.
    \item Verify that the object we found has precisely the singularities that GRDB predicts for the $\QQ$-Fano threefold. If it does not, mutate the polynomial (see Section \ref{LPintro} for the precise definition) and retry step 4 with its mutations.
\end{enumerate}

\begin{rem}
\begin{itemize}
\item Steps 1--3 are the starting data of this paper, and are recorded in the database \cite{repMMLP}. The difficulty of the rest of the procedure is found in step 5, as we often construct varieties that are not $\QQ$-Fano (i.e. they have worse singularities). There are results in the literature \cite{Prince_2019,prince2019cracked} which guarantee that the setup will produce an orbifold ambient space, however they are not general enough to produce the examples of this paper. While these restrictions make sense when aiming for smooth varieties, they are no longer needed here and allow for a broader search.

\item Another significant advantage is implicit within step 5: the access to the mutation graph of a given polynomial is entirely automatised \cite{fanosearchcore} in the computational algebra system (or \emph{CAS}) Magma \cite{Magmapaper}, and each corresponding polytope can be visualised. We record in Tables \ref{tab1} and \ref{tab2} the iteration in the mutation graph (in the column MG\#), denoting the amount of times we have mutated before being successful with our construction. It ranges between 1 and 31. 
\end{itemize}
\end{rem}

There is of course no guarantee that a given $\QQ$-Fano threefold can be realised as a toric complete intersection or as a toric Pfaffian variety, yet it is very encouraging that many of them are, and is certainly a good indication of which entries in GRDB correspond to actual geometric objects. It is likely that some varieties in Table \ref{tab2} are birationally linked to one another via projections/unprojections (e.g. Q$33018$ and Q$33019$), as suggested in the GRDB. We chose not to explore this route, as our construction embeds them in ambient spaces which seem unrelated. Nonetheless, it is entirely conceivable to correlate blow-ups/contractions on shapes with projections/unprojections of $X$ inside $Y$. 

Our constructions have been especially easy in the case of varieties whose baskets contain only $\frac{1}{2}(1,1,1)$ singularities, which could be explored in the future as a systematic classification. It is unclear how much such an endeavour would overlap with Takagi's lists in \cite{takagi_2002}, because computing the Picard rank of a complete intersection inside a non-orbifold toric variety can be subtle. 

A natural direction to pursue is the automation of this process: indeed, many of the steps involved rely heavily on a CAS. We are however far from a systematic implementation: choosing a shape which yields a good construction of the ambient space is still very much a matter of experience as opposed to brute computational force. Our attempts at writing a program that generates all possible shapes (and then scaffoldings) on various iterated mutations of a given rigid MMLP currently yield overwhelmingly many pathological cases. Though ultimately this is our desired direction, we cannot yet propose an algorithm that is applicable at the large scale needed to construct all varieties in the GRDB.

\paragraph{Acknowledgements.} 
I thank Tom Coates for suggesting this project, for a very gentle introduction to the topic during my time as a postdoc at Imperial College and for encouraging me to make it (also) my own. I thank Alessio Corti for his invaluable help in many discussions in which he habitually provided clear answers to my numerous and often convoluted questions on toric geometry. I am very grateful to Thomas Prince, who, among other things, pointed out the construction in Section \ref{noptstrut} which unlocked many examples of this paper. I also thank Al Kasprzyk for being an excellent troubleshooter when it comes to polynomials, combinatorics and the GRDB, in particular when pointing out varieties in Section \ref{Whynot} belonged to Table~\ref{tab2} instead of Table \ref{tab1}. I am grateful to Giuseppe Pitton who, together with his aforementioned co-authors, provided me with a list of Laurent polynomials to start from. I was funded by the Projet \'Etoiles Montantes GeBi de la R\'egion Pays de la Loire, whose PI Susanna Zimmermann I thank for her generous support.

\section{Preliminaries}

\subsection{Laurent polynomials and mutations}
\label{LPintro}
Our starting point is a class of Laurent polynomials intrinsically linked with the concept of mutation. For the purpose of self-containment, we reproduce the terminology below, most of which we extract from \cite{coates2021maximally}, a paper entirely dedicated to these objects. 

Let $P$ be a full-dimensional polytope in a lattice $N$. It is a \textit{Fano polytope} if $0\in P\setminus \partial P$ and its vertices are primitive elements in $N$. This implies that $X_P$, the toric variety whose fan is the \textit{spanning fan} of $P$, is a (usually very singular) Fano variety. We consider Laurent polynomials $f=\sum\limits_{v\in P\cap N} c_v x^v$ with coefficients $c_v\in \NN$ such that Newt$(f)=P$, and we say they are \textit{supported on $P$}. We always set $c_0=0$, as two polynomials which differ by a constant have the same mirror (see \cite[Remark 2.9]{LaurentInversion}). We remember this convention when using Laurent inversion, where we allow the origin to be part of any strut.

Lastly, we assume, as in Convention 2.3 in \cite{coates2021maximally}, that all Laurent polynomials (and mutation factors) have non-negative integer coefficients and that if $f\in \CC[N]$ then the exponents of monomials in $f$ generate $N$.

\begin{dfn}{\cite[Def.1.6.]{coates2021maximally}}
Let $N$ be a lattice and let $w \in M$ be a primitive vector in the dual lattice. Then $w$ induces a grading on $\CC[N]$. Let $a\in \CC[w^\perp \cap N]$ be a Laurent polynomial in the zeroth piece of $\CC[N]$, where $w^\perp \cap N= \{v \in N | w(v) = 0\}$. The pair $(w, a)$ defines an automorphism of $\CC(N)$ given by \[\mu_{w,a} \colon \CC(N) \to \CC(N), \ x^v \mapsto x^v a^{w(v)}.\] Let $f\in \CC[N]$. We say that \textit{$f$ is mutable with respect to $(w,a)$} if
\[g := \mu_{w,a}(f) \in \CC[N],\] in which case we call $g$ a mutation of $f$ and $a$ a factor.
\end{dfn}

\begin{dfn}{\cite[Def.2.4.]{coates2021maximally}} Given a Laurent polynomial $f$, consider the graph $G$ with vertex labels that are Laurent polynomials and edge labels that are pairs $L(w,a)$, defined as follows. Write $l(v)$ for the label of a vertex $v \in V(G)$, and $l(e)$ for the label of an edge $e \in E(G)$.
\begin{enumerate}[label=(\roman*)]
\item Begin with a vertex labelled by the Laurent polynomial $f$.
\item Given a vertex $v$, set $g:=l(v)$. For each $(w,a)$, $\deg a > 0$, such that $g$ is mutable with
respect to $(w,a)$ and either:
\begin{enumerate}[label=(\alph*)]
	\item there does not exist an edge with endpoint $v$ and label $L(w,a)$; or 
	\item for every edge $e=\overline{vv'}$ with $l(e) = L(w,a)$ we have that \[l(v')\notin \mu_{a,x^{w^\perp \cap N}}(g);\] 
pick a representative $g'\in \mu_{w,ax^{w^\perp}\cap N}(g)$ and add a new vertex $v'$ and edge $\overline{vv'}$ labelled by $g'$ and $L(w,a)$, respectively.
\end{enumerate}
\end{enumerate}
The mutation graph $\mathcal{G}_f$ of $f$ is defined by removing the labels from the edges of $G$ and changing
the labels of the vertices from $g$ to the $\textup{GL}(N)$-equivalence class of Newt$(g)$.
\end{dfn}


\begin{dfn}{\cite[Def.2.5.]{coates2021maximally}} We partially order the mutation graphs of Laurent polynomials by saying that $\mathcal{G}_f\prec \mathcal{G}_g$ whenever there is a label-preserving injection $\mathcal{G}_f \hookrightarrow \mathcal{G}_g$. A Laurent polynomial $f$ is \textit{maximally mutable} (or for short, $f$ is an MMLP) if Newt$(f)$ is a Fano polytope, the constant term of $f$ is zero, and $\mathcal{G}_f$ is maximal with respect to $\prec$.
\end{dfn}

\begin{dfn}{\cite[Def.2.6.]{coates2021maximally}} An MMLP $f$ is \textit{rigid} if the following holds: for all $g$ such that the constant term of $g$ is zero and Newt$(f) =$Newt$(g)$, if $\mathcal{G}_f=\mathcal{G}_g$ then $f=g$.
\end{dfn}

\subsection{Toric varieties, weight matrices and complete intersections}

The $\QQ$-Fanos that we are interested in live inside toric ambient spaces, either as complete intersections or as Pfaffian varieties.  As in the original Laurent inversion paper \cite[\S 2.]{LaurentInversion}, the most useful description of these ambient spaces is via their GIT data: having global coordinates and equations then enables us to compute singularities. This characterisation is equivalent to the description via fans and cones and it is the main output of our paper. The Laurent inversion construction typically involves three incarnations of toric varieties: the starting object $X_P$ is given by a fan, the shape $Z$ appears (almost exclusively) via the polytopes of its nef divisors, and the ambient space $Y$ appears as a GIT quotient, all of which are necessary when discussing examples. We present a minimal set of preliminaries, adapted from \cite[\S 3.1]{1/3} to include the non-orbifold case. 

Let $\TT\simeq \CC^{\times n}$ be an $n$-dimensional torus, $M=\textup{Hom}(\TT,\CC^\times)$ its character lattice and $N=\textup{Hom}(M,\ZZ)$ the lattice dual to $M$. 

Let $P$ be a Fano polytope with vertices in $N$, and $\Sigma$ its spanning fan. We construct the projective toric variety $F:=X_\Sigma$ by glueing together the affine open subsets $X_\sigma=\Spec \CC[M\cap \sigma^\vee]$, where $\sigma\in\Sigma$. Denote by $\rho_1, \ldots, \rho_m\in N$ the primitive generators of the rays of $\Sigma$, which we assume generate $N$ as a group. Note that $F$ is not necessarily an orbifold, as the fan need not be simplicial, and could have at worst klt singularities.

We have the following two dual exact sequences:

\begin{equation} \label{dual_exact_sequences}
  \begin{aligned}
    \xymatrix{
      0 \ar[r] & \LL \ar[r] & \ZZ^m \ar[r]^{\rho} & N \ar[r] & 0 \\
      0 & \ar[l] \LL^* & \ar[l]_D (\ZZ^m)^* & \ar[l]_-{\rho^T} M & \ar[l] 0 
    }
    \end{aligned}
\end{equation}

The lattice $\LL^*$ is in fact the $\TT$-class group of  $F$, we call $\rho$ the \textit{ray map} (which in this paper will always be surjective) and $D$ the \textit{divisor homomorphism} of $F$. If $\GG$ is the torus with character group $\LL^*$, then $D$ is dual to a morphism $\GG \to \CC^{\times m}$ and it is via this morphism that $F=\CC^m/\!\!/_\omega\GG$. In other words, a matrix of the morphism $D$ is a \textit{weight matrix} of the $\GG$-action on $\CC^m$.

More precisely, let $x_i$ be the standard basis of $(\ZZ^m)^*$. Note that since $F$ is projective, the cone $\mathcal{C}\subset \LL^*_\RR$ spanned by $D_i=D(x_i)$ is strictly convex. There is a wall-and-chamber decomposition of $\mathcal{C}$ called the \textit{secondary fan}, and we choose $\omega$ inside a (not necessarily full-dimensional) cone of this fan. Now set $V_\omega=V(\textup{Irr})$, where \[\textup{Irr}=\left( x_{i_1}\cdot \ldots \cdot x_{i_k} | \ \omega\in\RelInt \langle D_{i_1}, \ldots, D_{i_k}\rangle  \right)\] is the \textit{irrelevant ideal}. We then write $F$ as a quotient \[F=(\CC^m\setminus V_\omega)/\GG,\] by the $\GG$-action given by $D$.

\begin{rem}[Birational POV]
\begin{enumerate}\item The choice of $\omega$ coincides with choosing the ample cone inside the secondary fan - every divisor in the same cone-interior as $\omega$ is ample. When working with toric Fano varieties, we usually bypass this discussion by imposing that $\omega=\sum\limits_i D(x_i)=-K_{F}$, and the stability cone is the smallest-dimensional cone that contains it. When working with Fano complete intersections or codimension $3$ Pfaffian varieties, we impose similar conditions by using adjunction-type formulas (see Sections \ref{CIAdj} and \ref{PfK}, respectively).
\item The rays of $\Sigma$ are determined by the weight matrix and vice-versa, modulo appropriate basis changes. The quotient $\CC^m /\!\!/_{\omega}\GG$ also depends on the chamber in which $\omega$ lies, and two stability conditions that belong to the same chamber induce isomorphic varieties. Varying $\omega$ between two adjacent chambers inside the secondary fan produces varieties that are related either by small or divisorial contractions. In terms of the fan, choosing $\omega$ in a maximal-dimensional chamber implies all cones in $\Sigma$ are simplicial, i.e. $F$ is an orbifold. Choosing $\omega$ on a wall may produce a non-$\QQ$-Gorenstein variety.
\end{enumerate}
\end{rem}

The choice of a stability condition determines the $\GG$-invariant open sets $V_{i_1\ldots i_k}=\{ x_{i_1}\neq0, \ldots, x_{i_k}\neq 0\} \subset \CC^m$ such that \[F= \bigcup\limits_{\{(i_1,\ldots i_k)| \omega \in \RelInt(\langle D_{i_1}, \ldots, D_{i_k} \rangle)\}} [V_{i_1,\ldots, i_k}/\GG].\]

Since we don't restrict to the orbifold case, we do not impose that $k=r$, therefore the toric charts $V_{i_1\ldots i_k}$ are not necessarily quotients by finite groups.  Indeed, when working inside these charts we write $[V_{i_1\ldots i_k}/\GG]=[U_{i_1\ldots i_k}/\Bmu]$, where $U_{i_1\ldots i_k}=\{x_{i_1}=\ldots =x_{i_k}=1\}$ and $\Bmu$ is the subgroup of $\GG$ which fixes $U_{i_1\ldots i_k}$. More precisely, $\Bmu$ is the group whose character group is the cokernel of the homomorphism \[D_{i_1 \ldots i_k}\colon (\ZZ^k)^*\to \LL^*.\] 

\subsubsection{Complete intersections}
\label{CIAdj}
We consider complete intersections $X\subset F$, i.e. general elements of linear systems $|L_1|, \ldots, |L_c|$ on $F$, where $L_i\in \LL^*$ are  $\GG$-linearised line bundles on $F$. The space of sections $H^0(F,L_i)$ is a vector subspace of $\CC[x_1,\ldots, x_m]$ with basis consisting of monomials of homogeneity type $L_i$. Let $f_i\in H^0(F,L_i)$, then $V(f_1,\ldots, f_c)$ is stable under the action of $\GG$ and we consider the subvariety $X=(V(f_1, \ldots, f_c)\setminus V_\omega)/\GG\subset F$. 

\begin{dfn}\begin{enumerate}\item Assume $X\subset F$ avoids the non-$\QQ$-Gorenstein locus of $F$. Then $X$ is \textit{quasi-smooth} if \[V(f_i, \ldots, f_c)\setminus V_\omega \subset \CC^m\setminus V_\omega.\]
\item Suppose $X\subset F$ is quasi-smooth. We say that $X$ is \textit{well-formed} if for every toric stratum $S\subset F$ with nontrivial stabiliser, $S\subset X$ implies $\codim_XD\geq 2$.
\end{enumerate}
\end{dfn}
If $X\subset F$ is a complete intersection as above, its canonical class is \[-K_X=-K_F-\sum L_i.\] There is a similar formula for Pfaffian varieties, which we work out in detail in section \ref{PfK}.

\subsection{Laurent inversion}
\label{LIsec}
Finally, we recall the algorithm of Laurent inversion in \cite{LaurentInversion}, whose purpose is to provide an embedding $X_P \hookrightarrow Y$, where $X_P$ and $Y$ are toric varieties, constructed using a decomposition  into summands of a rigid MMLP supported on $P$. We then use the coordinates on this new ambient space to smooth $X_P$ into a $\QQ$-Fano variety $X$.

\begin{dfn}{\cite[Def.3.1.]{LaurentInversion}} \label{ScafDef} Fix the following data:
\begin{enumerate}[label=(\roman*)]
\item a lattice $N$ together with a splitting $N = \overline{N} \oplus N_U$;
\item the dual lattice $M := \textup{Hom}(N, \ZZ)$, with the dual splitting $M = \overline{M}\oplus M_U$;
\item a Fano polytope $P \subset N_\QQ$;
\item a projective toric variety $Z$ given by a fan in $\overline M$ whose rays span the lattice $\overline{M}$.
Given such data, a \textit{scaffolding} $S$ of $P$ is a set of pairs $(D,\chi)$ where $D$ is a nef divisor on $Z$ and $\chi$ is an element of $N_U$ , such that \[P=\textup{conv}( P_D+\chi \ | \ (D,\chi)\in S).\]
We refer to $Z$ as the \textit{shape} of the scaffolding, and the elements $(D, \chi) \in S$ as \textit{struts}.
\end{enumerate}
\end{dfn}

The algorithm in \cite{LaurentInversion} directly produces the weight matrix of $Y$ and assumes the following

\begin{rem}{\cite[Remark 3.3]{LaurentInversion}} \label{Uneliminated} Given the $N=\overline{N}\oplus N_U$ splitting in Definition \ref{ScafDef}, there are $u=|U|$ struts of type $(0,e_i)$ where $\{e_i\}_{i=1}^u$ is a basis of $N_{U}$. 
\end{rem}

What happens when this assumption is not fulfilled appears in \cite{Prince_2019} and a detailed example can be found in section \ref{noptstrut} of this paper. Many entries in the systematic Table \ref{tab1} use this modification, and it is likely that an automatisation of Laurent inversion should include it. 

We now state the algorithm, which appears in \cite{LaurentInversion} with a sign typo at step $(ii)$ (as confirmed by \cite{LaurentInversion}[Example 6.3, Example in \S 7], which are correct constructions). To follow along its steps in a simple example, go to Section \ref{basic}.

\begin{algorithm}{\cite[Alg.5.1]{LaurentInversion}} \label{LIalg} Let $S$ be a scaffolding of a Fano polytope $P$ with shape $Z$. Let $u = \dim N_U$, $r = |S|-u$ and assume Remark \ref{Uneliminated}. Let $R$ be the sum of $|S|$ and the number $z$ of rays of $Z$. We determine an $r \times R$ matrix $\MM$, which is the weight matrix of the toric variety $Y$, as follows. Let $m_{i,j}$ denote the $(i, j)$ entry of $\MM$. Fix an identification of the rows of $\MM$ with the $r$ elements $(D_i,\chi_i)$ of $S$ which do not correspond to the basis of $N_U$, and an ordering $\Delta_1, \ldots , \Delta_z$ of the toric divisors in $Z$. 
\begin{enumerate}[label=(\roman*)]
\item For $1\leq j \leq r$ and any $i$, let $m_{i,j} =\delta_{i,j}$.
\item For $1\leq j \leq u$ and any $i$, let $m_{i,r+j}$ be determined by the expansion
\[ \chi_i = -\sum\limits_{j=1}^{u}m_{i,r+j}e_j.\]
\item For $1 \leq j \leq z$, let $m_{i,|S|+j}$ be determined by the expansion 
 \[ D_i =\sum\limits_{j=1}^{z} m_{i,|S|+j} \Delta_j. \]
 \end{enumerate}
The weight matrix $\MM$ alone does not determine a unique toric variety -- we also need to choose a stability condition $\omega$. Let $Y$ denote the toric variety determined by this choice. Unless otherwise stated, we take $\omega$ to be the sum of the first $|S|$ columns in $\MM$. This translates to the adjunction-type formulas we mentioned earlier, to ensure that the complete intersection inside $Y$ is a Fano variety.
\end{algorithm}
\vspace{-0.5cm}
\paragraph{Warning.} This last remark does not mean all varieties in Tables \ref{tab1} and \ref{tab2} have stability conditions given by the sum of the first $|S|$ columns of the ambient weight matrix, as we may have solved for one of those columns to eliminate an equation. When reading data from the tables, it is safer to apply the adjunction formulas in \ref{CIAdj} or \ref{PfK} to determine $\omega$.
\medskip

The construction of the line bundles $L_i$ is described explicitly in \cite[Proposition 12.2]{LaurentInversion}, which connects the algorithm to the method in \cite{doran_harder_2016} in the case of complete intersections. Each ``level" of the tower determines a corresponding line bundle: note that by (iii) in \ref{LIalg} the last $z$ columns of $\MM$ correspond to rays of $Z$. Since the rays in the fan of $Z$ are partitioned into the $c$ factors of the tower, so are then the Cox coordinates in this block. The $\GG$-weight of a line bundle will then be the same as that of $\textup{div}(\Pi x_j)$, where $x_j$ are Cox coordinates on $Y$ corresponding to a factor of the partition. In \cite{doran_harder_2016}, the embedding construction of $X_P$ inside $Y$ depends on the order of the $L_i$, and implicitly on the order of the factors. As we are looking to construct the general section of $\oplus_{i=1}^c L_i$ (and are not particularly concerned with a specific toric section), it is then enough to add all monomials in the same $(\CC^\times)^r$-eigenspace as $L_i$ for each $i=1\ldots c$. In other words, as soon as we choose $Z$ and $S$, we directly determine the equations of the smoothest possible $X$ given by this scaffolding.

\subsubsection{Pfaffians}
\label{PfK}
By the Buchsbaum-Eisenbud structure theorem \cite[Theorem 2.1]{BuchsbaumEisenbud}, a codimension three ideal inside a Gorenstein ring is given by the $2n \times 2n$ diagonal Pfaffians of a skew $(2n+1) \times (2n+1)$ matrix. In all such examples in this paper, $n=2$. 

Some entries in Tables \ref{tab1} and \ref{tab2} are in this format. Additionally, the codimension $4$ object in Section \ref{3DPf} is a complete intersection between a Pfaffian variety and a hyperplane, which we eventually reduce to the Pfaffian case. We now recall the information the format encodes in the case of weight matrices obtained by Laurent inversion. 

To understand the singularities and degree of such a variety, we determine the eigen-type of the anticanonical divisor of a variety embedded as a Pfaffian inside a toric variety; this is the equivalent of an adjunction theorem in the case of complete intersections. We then impose that the stability condition of the ambient space be in the same cone as the anticanonical divisor, ensuring that the resulting variety is Fano. 

We follow the reasoning in \cite[Section 3.3]{1/3}, which treats a surface example. In our context, $X$ is typically embedded inside an ambient space $Y$ with weight matrix:
\[\begin{array}{c cccccc}
&U & t_1 & t_2 & t_3 & t_4 & t_5 \\
\hline
I_r &C_U & C_1 & C_2 & C_3 & C_4 & C_5
\end{array}\]
where $C_U, C_i\in \CC^r$ are vectors representing the weights of the $\GG$-action on the column corresponding to $N_U$ (there will be only one, for dimension reasons) and the Cox coordinates $t_i$, respectively.  Then $-K_Y=\begin{pmatrix} 1\\ \vdots \\ 1 \end{pmatrix} + C_U + \sum C_i$. The antisymmetric matrix we consider is of the form

\[ A=\left( \begin{array}{ccccc}
0 & * & t_1+* & t_2+* & * \\
& 0 & * & t_3+*& t_4+* \\
&-sym& 0 & * &t_5+*\\
&&&0& * \\
&&&& 0
\end{array}\right)\]
and its $4\times 4$ Pfaffians, i.e. the equations of $X$, are
\begin{equation}
  \left\{
    \begin{aligned}
      t_3t_5+t_4*+\ldots &=0\\
      t_2t_5+t_1*+\ldots &=0\\
      t_2t_4+t_3*+\ldots &=0\\
      t_1t_4+t_5*+\ldots &= 0\\
      t_1t_3+t_2*+ \ldots &= 0
    \end{aligned}
  \right.
\end{equation}
where the forms ``$*$" are such that the entries of $A$ are of fixed weights and the Pfaffian equations are homogeneous. We denote the line bundles whose sections are given by these equations by $L_i$, $i=1\ldots 5$ respectively. Note that this determines the weights of the $\GG$-action on the line bundles (e.g. $L_1\sim C_3+C_5$). 

Let $E=\oplus_{i=1}^5L_i$. Then there exists a linearised line bundle $L$ on $Y$ such that $X\subset Y$ is the degeneracy locus of a general antisymmetric homomorphism $s:E \otimes L\to E^\vee$ defined by the matrix $A$.\footnote[3]{the order of the $L_i$ factors in $E$ is important: it must be compatible with the columns of $A$.} The specific format of $A$ implies that $L=-\dfrac{\sum_{i=1}^5 L_i}{2}=-\sum_{i=1}^5 C_i$.
\paragraph{Claim:} (the adjunction formula for Pfaffians) The weights of the $\GG$-action on $-K_X$  are given by \[-K_X=(-K_Y+L)|_X = (-K_Y-\dfrac{\sum_{i=1}^5 L_i}{2})|_X.\]

This computation requires that $-K_Y$ be locally free. We denote by $\omega_Y$ its restriction to the $\QQ$-factorial locus of $Y$. The locus which we avoid is at most of dimension one in all examples in this paper. A posteriori, $X$ avoids it completely (if we follow through with the computation and $X$ intersects it, then it is not $\QQ$-factorial and we discard the example).

 We use the homomorphism $A$ to construct a resolution of $\cO_X$ 
  \begin{equation} \label{resolution}
  \begin{aligned}
    \xymatrix{
      0 \ar[r] & L \ar[r]^{\Pf^\vee} & E\otimes L \ar[r]^{A} & E^\vee \ar[r]^{\Pf} & \cO_Y \ar[r] & \cO_X \ar[r] & 0 
    }
    \end{aligned}
\end{equation}
 We want to compute $\omega_X=\textup{\underline{Ext}}^3_{\cO_Y}(\cO_X,\omega_Y)$, so we use this resolution to induce a complex:
 \begin{align*} 
   \label{Homresolution}
      0 \rightarrow \textup{\underline{Hom}}_{\cO_Y}(\cO_Y,\omega_Y) \rightarrow  \textup{\underline{Hom}}_{\cO_Y}(E^\vee,\omega_Y) \rightarrow  \textup{\underline{Hom}}_{\cO_Y}(E\otimes L,\omega_Y)\rightarrow \textup{\underline{Hom}}_{\cO_Y}(L,\omega_Y) \rightarrow  0 
    \end{align*}
whose 3rd cohomology is $\omega_X$. This is the same as the self-dual complex \begin{equation} \label{chopped}
  \begin{aligned}
    \xymatrix{
      0 \ar[r] & L \ar[r]^{\Pf^\vee} & E\otimes L \ar[r]^{A} & E^\vee \ar[r]^{\Pf} & \cO_Y \ar[r] & 0 
    }
    \end{aligned}
\end{equation}
tensored by $L^\vee\otimes \omega_Y$. The complex \eqref{chopped} is everywhere exact except at  $\cO_Y$, where the cohomology is $\cO_Y/\Pf(E^\vee)=\cO_X$ (since \eqref{resolution} is a resolution), thus we conclude that $\omega_X=L^\vee\otimes \omega_Y\otimes \cO_X$. This proves the claim.

\subsection{The tables}

Table \ref{tab1} is ordered by codimension and $\QQ$-Fano ID in the Graded Ring Database, that is, the IDs which appear in the ``Fano 3-folds or big table" database. Table \ref{tab2} is ordered only by the latter. We also include the polytope ID (from the ``Toric canonical Fano 3-folds" database) to which we initially associate a rigid MMLP, and the expected singularities of each $\QQ$-Fano.

\begin{notation}
Since we regularly refer to objects in two different databases, to avoid confusion we refer to GRDB IDs in the $\QQ$-Fano threefolds database using the prefix Q (e.g. Q$38989$) and the GRDB IDs in the toric canonical Fano threefolds using the prefix P (e.g. P$519468$). We sometimes directly refer to such a variety by its polytope, which uniquely determines it.
\end{notation}

For simplicity, we do not record the original rigid MMLP (i.e. that is supported on a \textit{canonical} toric variety), as this may not give the successful partial smoothing. Instead, what we record is its mutation for which Laurent Inversion yields a $\QQ$-Fano threefold, in the ``Rigid MMLP" column. Its corresponding toric Fano might therefore not be canonical, but klt. 

\begin{notation}
When discussing Laurent inversion, we always denote the toric variety which we embed as $X_P$, the ambient space which we embed it in by $Y$, the shape of the embedding by $Z$ and the (partial) $\QQ$-Fano smoothing by $X$. If $X_P$ is embedded as a toric complete intersection given by line bundles $L_1,\ldots, L_c$, then $X$ is a general section of these line bundles. We use the coefficient convention for complete intersections in toric varieties: all coefficients of monomials in any equation equal $1$. 
\end{notation}

If a table entry is a complete intersection given by $\GG$-linearised line bundles $L_1,\ldots, L_c$, the stability condition is necessarily $\omega=-K_Y-\sum L_i=\sum \textup{div} (x_j)-\sum L_i$, where $\{x_j\}$ are the Cox coordinates on $Y$.  We discuss how to obtain the weight matrices in this case, as well as how to analyse the singularities of a general section of $\oplus L_i$ in \ref{2Dexamples}.  If a table entry is a codimension 3 Pfaffian variety given by the line bundles $L_1,\ldots, L_5$, then cf. Section \ref{PfK} the stability condition is $\omega=-K_Y-\frac{1}{2}\sum_{i=1}^{5} L_i=\sum \textup{div} (x_j) -\frac{1}{2}\sum_{i=1}^{5} L_i $, with the same notation for the $x_j$. A Pfaffian example is discussed in \ref{3DPf}.

The shape $Z$ we use to produce these weight matrices using Laurent inversion is also included in Tables \ref{tab1} and \ref{tab2}. We mention if we use point-struts as a basis of $N_U$, and the difference of how we obtain the ambient space in these cases is discussed in examples \ref{basic} and \ref{noptstrut}.

\begin{rem}As a sanity check after producing a weight matrix, we verify that the period sequence of the general element in the complete intersection (i.e. that obtained from its regularised quantum period) is the same as the period sequence of the Laurent polynomial (obtained by computing its classical period). The theory for computing the former exists only in the context of orbifold ambient spaces, however since we ultimately check that $X\subset Y$ completely avoids the non-orbifold locus, we can assume $Y$ is an orbifold by taking a partial crepant resolution if necessary. For the theory, we require an equivalent of \cite[Corollary D.5]{Coates_2016} in the orbifold case. This result is a combination of Givental's mirror theorem \cite[Theorem C.1]{Coates_2016} and the Quantum Lefschetz Theorem \cite[Theorem D.1]{Coates_2016}, the appropriate generalisations of which can be found in \cite[Theorem 31]{Coates_2015} and \cite[Thm 1.1]{Wang2019} respectively. We verify that the first 10 terms of both period sequences coincide using Magma, and we include the results in Tables \ref{tab1} and \ref{tab2}. 
\end{rem}

\begin{rem}A formula for computing the regularised quantum period of a codimension $3$ Pfaffian variety is not yet known. We however compute the anticanonical degree and singularities in each case (an example of this is thoroughly discussed in Section \ref{degree}), which, along with \cite[Thm 5.5]{LaurentInversion}, we consider as sufficient evidence of having constructed the correct objects.
\end{rem}

\subsection{How Table \ref{tab1} is constructed}
\label{Whynot}

We now sketch how the number $54$ is obtained in Theorem \ref{systemThm}. Assume we have followed Steps 1-3 in the ``Input of the construction" in the introduction, but only starting from polytopes $P$ such that there exists a candidate $\QQ$-Fano threefold of GRDB~-~codimension $\geq 20$ with the same Hilbert series as $X_P$. In the ``Fano 3-folds" list, no Fano index $>1$ varieties have GRDB-codimension~$\geq 20$, and we look to translate this into a condition on the polynomials we obtain after Steps~1-3.


In order to extract the correct period sequences we use the criterion below. This is not proved in general, but we can verify it by hand in our situation - there is a handful of cases to examine, and they are either trivial\footnote{High Fano index varieties tend to have low GRDB-codimension. For these examples, the embedding into weighted projective space is enough to establish the F/LG correspondence.} or included in Table \ref{tab3}.

\begin{crt} \label{perzeroes}
Assume a rigid MMLP $f$ is mirror to a $\QQ$-Fano threefold $X$ and let $(p_n)_{n\geq 0}$ be their common period sequence. Then $X$ is of Fano index $r$ iff $p_k=0$ for all $k\not\equiv 0 (\textup{mod } r)$. \end{crt}

In other words, periodic zeroes in the period sequence of such an $f$ suggests that a corresponding $\QQ$-Fano built via Laurent inverting $f$ does not belong to Table \ref{tab1}. After discarding high Fano index cases, we finally obtain 54 period sequences which can be found in Table \ref{tab1}. Note that some of these appear to have a periodicity of zeroes, however this disappears if we compute more than 10 terms of its period sequence. 

\subsection{High Fano index entries in Table \ref{tab2}}

We discuss one of the non-trivial examples illustrating Criterion \ref{perzeroes}. The other cases are similar and can be found in Table \ref{tab3}.

The polynomial \[f = \frac{x^6}{y^4z^5} + \frac{2x^3}{y^2z^2} +\frac{2x^2}{y^3z^2} + x + y + z + \frac{1}{y} + \frac{2z}{xy} + \frac{z}{x^2y^2}\] deceptively appears as an output of Steps 1-3 after the restriction to codimension $\geq 20$. It is the unique rigid MMLP supported on a mutation of the canonical polytope P$543951$. Its period sequence suggests a higher Fano index: \[[ 1, 0, 2, 0, 6, 0, 20, 0, 3430, 0, 75852 , \ldots].\] When Laurent inverting $f$ we obtain a hypersurface $X$ in the toric variety with weight matrix  
\[\begin{array}{cccccc|c} 1&0&0&0&0&1&0 \\ 0&3&1&5&2&3&6  \end{array},\] and as $-K_X=(2,4)$, indeed the variety is of Fano index $2$.

We now need to find the correct GRDB ID of the $\QQ$-Fano threefold we have built, which we deduce will appear polarised by $A=-K_X/2$. Keeping in mind that $-K_X=-K_{X_f}=-K_{X_{\textup{P}543951}}=168/5$, we search by the Fano index and the degree of the polarisation $A^2=1/4(-K_X)^2=21/5$ and find Q$40971$ as the only candidate.\footnote{If more precision is necessary, we may also search by a subsequence -- in this case, every second term -- of the Hilbert series associated to $A$ : $[1, *,  19, *, 87, *, 239, *, 509, \ldots]$ which we deduce from the Hilbert series of $X_f$.} For this reason, the object we build from P$543951$ appears as variety \hyperlink{40971}{39} in Table \ref{tab2} and does not belong in Table \ref{tab1} of high codimension objects. 

There are a total of five new Laurent inversion constructions in this situation among the ones in our classification, and three varieties (marked with an asterisk) that had previously appeared in \cite{TomStephen}, albeit in a different format (in fact, Brown and Suzuki already suggest Q$41200$ can be constructed in \cite{BS07b}). All of these entries are now part of Table \ref{tab2}.

\begin{center}
\arrayrulecolor{Gray}
\setlength{\arrayrulewidth}{0.2mm}
\rowcolors{1}{Gray!10}{Periwinkle!20}
\begin{longtable}[ht]{|  p{0.3cm} | >{\centering\arraybackslash}p{1.5cm} | >{\centering\arraybackslash}p{5.6cm}|>{\centering\arraybackslash}p{1cm} |>{\centering\arraybackslash}p{1.8cm} | }
\hline
\label{tab3}

 &  Polytope \newline ID & Weight matrix \newline \& bundles & Fano \newline index & $\QQ$-Fano \newline ID \\
 
 \hline
  
1  & P$543951$ &  \begin{math}\begin{array}{cccccc|c} 1&0&0&0&0&1&0\\ 0&3&1&5&2&3&6 \end{array}\end{math} & 2 & \hyperlink{40971}{Q$40971$} \newline codim 9 \\

\hline

2 & P$534760$ & \begin{math} \begin{array}{cccccc|c} 1&-3&-5&-1&-2&2&-6 \\ 0&3&5&1&2&-1&6 \end{array} \end{math}  & 2 &  \hyperlink{40948}{Q$40948$} \newline codim 9 \\
\hline
 
3 & P$473887$ & \begin{math}\begin{array}{ccccccc|c} 1&0&0&1&-1&-1&0&0 \\ 0&1&0&0&0&0&1&0\\ 0&0&1&-1&3&2&1&2 \end{array}\end{math} & 2 &  \hyperlink{40993}{Q$40993$} \newline codim 6 \\
 
 \hline
 
4 & P$413267$ & \begin{math}\begin{array}{cccccc|c} 1&0&1&-1&1&2&2\\ 0&1&1&2&-1&1&2 \end{array} \end{math}  & 2 &   \hyperlink{40988}{Q$40988$} \newline codim 6 \\ 

\hline

5 & P$402202$ & \begin{math}\begin{array}{cccccc|c} 1&0&2&-1&-1&4&2 \\ 0&1&-1&1&1&-2&0  \end{array} \end{math}   & 3 & \hyperlink{41251}{Q$41251$}\newline codim 5 \\
\hline

\multicolumn{5}{| l |}{}\\
\hline

6 &  P$547328$ & \begin{math}\begin{array}{ccccccc|cc} 1&0&2& 2&1&3&1& 3&4 \\ 0&1 & -1 &1&-1&-1&1&0&0\end{array}\end{math} & 3 & \hyperlink{41200}{Q$41200$}* \newline codim 4\\

\hline

 7  & P544064   & \begin{math}\begin{array}{ccccccc|cc} 1&0&3&0&1&-1&2& 3&0 \\ 0&1&-1&1&1&1&-1&0&2 \end{array}\end{math} & 3  &  \hyperlink{41218}{Q$41218$}*  \newline codim 4  \\
 
 \hline
 
8 & P543852  & \begin{math}\begin{array}{ccccccc|cc} 1&0&-1&1&1&1&-1& 0&2 \\ 0&1&4&0&1&-1&3&4&0 \end{array} \end{math} & 4&  \hyperlink{41334}{Q$41334$}* \newline codim 4  \\
\hline

\end{longtable}
\title{\textbf{Table \ref{tab3}: High Fano index varieties}}
\end{center}

\section{A systematic approach in high codimension}
\label{systematic}
In Table \ref{tab1} we provide a complete list of non-toric Fano varieties of codimensions between 20 and 23, which are mirror to rigid maximally-mutable Laurent polynomials, modulo the considerations in Section \ref{Whynot}.

\subsection{Typical examples}
\label{2Dexamples}

Recall Algorithm \ref{LIalg} directly produces the weight matrix of the ambient space $Y$, provided that one assumes Remark \ref{Uneliminated}. Throughout \cite{LaurentInversion}, the point-struts considered in Remark \ref{Uneliminated} are referred to as ```uneliminated variables". In \cite[Sec.3]{Prince_2019} Prince proves the embedding theorem without requiring the condition in Remark \ref{Uneliminated}. We discuss how to adjust the construction, as well as a basic incarnation of algorithm \ref{LIalg}, each on an example from Table \ref{tab1}.

\subsubsection{Laurent inversion with point-struts at the basis of $N_U$}
\label{basic}
As a warm-up, we demonstrate the construction in detail for threefold \#4 in Table \ref{tab1}. In other words, we start with a polynomial supported on the toric canonical Fano polytope P$519468$ and construct the $\QQ$-Fano threefold Q$38989$. The target variety has four isolated singular points: $2\times \frac{1}{2}(1,1,1)$, $\frac{1}{3}(1,1,2)$ and $\frac{1}{4}(1,1,3)$.

\begin{figure}[ht]
\centering
\begin{tikzpicture}[scale=0.35,every node/.style={scale=0.8}]
\draw[gray!80, thick] (0,3) -- (2,-2);
\draw[gray!80, thick] (5,-3) -- (2,-2);
\draw[gray!80, thick] (-6,-7) -- (2,-2);
\draw[gray!80, thick] (-11,-7) -- (2,-2);
\draw[gray!80, thick] (0,-9) -- (-6,-7);
\draw[gray!80, thick] (-11,-7) -- (-6,-7);

\draw[black, very thick] (0,3) -- (5,-3);
\draw[black, very thick] (0,-9) -- (5,-3);
\draw[black, very thick] (-2,-10) -- (5,-3);
\draw[black, very thick] (0,3) -- (-11,-7);
\draw[black, very thick] (0,-9) -- (-2,-10);
\draw[black, very thick] (-2,-10) -- (-11,-7);

\draw[gray!80, thick] (0.15,0.15) -- (-0.15,-0.15);
\draw[gray!80, thick] (-0.15,0.15) -- (0.15,-0.15);

\filldraw[gray] (0,0) circle (0pt) node[anchor=south] {\textcolor{black}{$(0,0,0)$}};
\filldraw[black] (0,3) circle (4pt) node[anchor=south] {$(0,0,1)$};
\filldraw[black] (5,-3) circle (4pt) node[anchor=north west] {$(1,1,-1)$};
\filldraw[gray] (2,-2) circle (4pt) node[label={[xshift=-0.9cm, yshift=-0.3cm]\textcolor{black}{$(0,1,-1)$}}]{};
\filldraw[black] (-11,-7) circle (4pt) node[label={[xshift=-0.3cm, yshift=-0.8cm]$(-2,-3,-2)$}]{};
\filldraw[gray] (-6,-7) circle (4pt) node[label={[xshift=1.3cm, yshift=-0.3cm]\textcolor{black}{$(-1,-2,-2)$}}]{};
\filldraw[black] (0,-9) circle (4pt) node[anchor=west] {$(1,-2,-2)$};
\filldraw[black] (-2,-10) circle (4pt) node[anchor=north] {$(1,-3,-2)$};

\filldraw[gray] (0,-3) circle (3pt) node[anchor=south east]{};
\filldraw[gray] (-2,-4) circle (3pt) node[anchor=south east]{};
\filldraw[gray] (-2,-1) circle (3pt) node[anchor=south east]{};
\filldraw[gray] (-3,-8) circle (3pt) node[anchor=south east]{};
\filldraw[black] (-5,-9) circle (4pt) node[anchor=south east]{};
\filldraw[black] (-8,-8) circle (4pt) node[anchor=south east]{};
\filldraw[black] (-4,-5) circle (3pt) node[anchor=south east]{};
\filldraw[black] (-7,-4) circle (3pt) node[anchor=south east]{};
\filldraw[gray] (-2,-4) circle (3pt) node[anchor=south east]{};
\end{tikzpicture}

\caption{The Newton Polytope of $g$}
\label{P519468}

\end{figure}
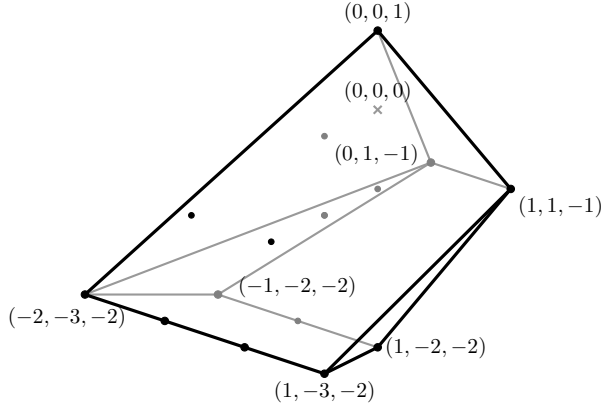

The only rigid MMLP supported on P$519468$ is:
\[f  =  x + \frac{xz^3}{y^2} + \frac{y^2}{z} + y + \frac{y}{xz^2} + \frac{y^3}{x^2z^3}.\] We use Magma to show that $f$ is mutation-equivalent (via $9$ mutation steps) to:
\[g=\frac{xy}{z} + \frac{x}{y^2z^2} + \frac{x}{y^3z^2} + \frac{y}{z} + z +\frac{ 2}{y^2z^2} + \frac{3}{y^3z^2} + 
\frac{1}{xy^2z^2} + \frac{3}{xy^3z^2} + \frac{1}{x^2y^3z^2},
\]
which is supported on the polytope in Figure \ref{P519468}.


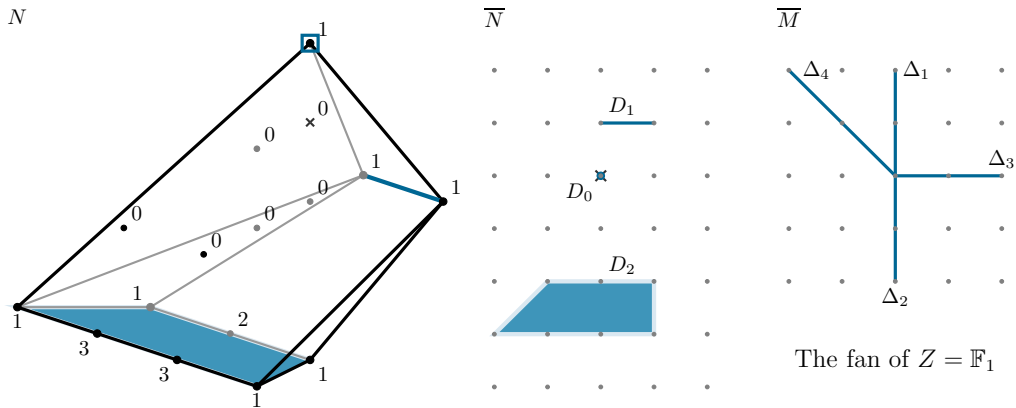
\begin{figure}[ht]
\begin{tikzpicture}[scale=0.35,every node/.style={scale=0.8}]
\centering
\node   at (-11,4) {$N$};
\filldraw[color=MidnightBlue!10, ultra thick, fill=MidnightBlue!60] (0,-9) -- (-6,-7) -- (-11,-7) -- (-2,-10); 
\draw[gray!80, thick] (0,3) -- (2,-2);
\draw[gray!80, thick] (5,-3) -- (2,-2);
\draw[gray!80, thick] (-6,-7) -- (2,-2);
\draw[gray!80, thick] (-11,-7) -- (2,-2);
\draw[gray!80, thick] (0,-9) -- (-6,-7);
\draw[gray!80, thick] (-11,-7) -- (-6,-7);

\draw[black, very thick] (0,3) -- (5,-3);
\draw[black, very thick] (0,-9) -- (5,-3);
\draw[black, very thick] (-2,-10) -- (5,-3);
\draw[black, very thick] (0,3) -- (-11,-7);
\draw[black, very thick] (0,-9) -- (-2,-10);
\draw[black, very thick] (-2,-10) -- (-11,-7);

\draw[black!80, thick] (0.15,0.15) -- (-0.15,-0.15);
\draw[black!80, thick] (-0.15,0.15) -- (0.15,-0.15);

\draw[MidnightBlue, ultra thick] (5,-3) -- (2,-2);


\filldraw[black] (0,3) circle (4pt) node[anchor=south west] {$1$};
\filldraw[black] (5,-3) circle (4pt) node[anchor=south west] {$1$};
\filldraw[gray] (2,-2) circle (4pt) node[anchor=south west] {\textcolor{black}{$1$}};
\filldraw[black] (-11,-7) circle (4pt) node[anchor=north] {$1$};
\filldraw[gray] (-6,-7) circle (4pt) node[anchor=south east] {\textcolor{black}{$1$}};
\filldraw[black] (0,-9) circle (4pt) node[anchor=north west] {$1$};
\filldraw[black] (-2,-10) circle (4pt) node[anchor=north] {$1$};

\filldraw[gray] (0,-3) circle (3pt) node[anchor=south west]{\textcolor{black}{$0$}};
\filldraw[gray] (-2,-1) circle (3pt) node[anchor=south west]{\textcolor{black}{$0$}};
\filldraw[gray] (-3,-8) circle (3pt) node[anchor=south west]{\textcolor{black}{$2$}};
\filldraw[black] (-5,-9) circle (4pt) node[anchor=north east]{$3$};
\filldraw[black] (-8,-8) circle (4pt) node[anchor=north east]{$3$};
\filldraw[black] (-4,-5) circle (3pt) node[anchor=south west]{\textcolor{black}{$0$}};
\filldraw[black] (-7,-4) circle (3pt) node[anchor=south west]{\textcolor{black}{$0$}};
\filldraw[gray] (-2,-4) circle (3pt) node[anchor=south west]{\textcolor{black}{$0$}};
\filldraw[gray] (0,0) circle (0pt) node[anchor=south west] {\textcolor{black}{$0$}};

\draw[MidnightBlue, very thick] (0.3,3+0.3) -- (-0.3,3+0.3) -- (-0.3,3-0.3) -- (0.3,3-0.3) -- cycle;


\end{tikzpicture}
\begin{tikzpicture}[scale=0.7,every node/.style={scale=0.8}, baseline=-110]  

\node   at (-2,2) {$\overline{N}$};
\draw[MidnightBlue,very thick] (1,0) -- (0,0);

\filldraw[color=MidnightBlue!10, fill=MidnightBlue!60, ultra thick] (1,-3) -- (-1,-3) -- (-2,-4) -- (1,-4) -- cycle; 
\draw[black!80, thick] (0.1,-0.9) -- (-0.1,-1.1);
\draw[black!80, thick] (-0.1,-0.9) -- (0.1,-1.1);

\filldraw[MidnightBlue] (0,-1) circle (2pt);
\node [anchor=south west]  at (0,0) {$D_1$};
\node [anchor=south west]  at (0,-3) {$D_2$};
\node [anchor= north east] at (0,-1) {$D_0$};

\filldraw[gray] (2,1) circle (1pt) node{};
\filldraw[gray] (0,1) circle (1pt) node{};
\filldraw[gray] (1,1) circle (1pt) node{};
\filldraw[gray] (-1,1) circle (1pt) node{};
\filldraw[gray] (-2,1) circle (1pt) node{};

\filldraw[gray] (2,-5) circle (1pt) node{};
\filldraw[gray] (0,-5) circle (1pt) node{};
\filldraw[gray] (1,-5) circle (1pt) node{};
\filldraw[gray] (-1,-5) circle (1pt) node{};
\filldraw[gray] (-2,-5) circle (1pt) node{};

\filldraw[gray] (2,-1) circle (1pt) node{};
\filldraw[gray] (0,-1) circle (1pt) node{};
\filldraw[gray] (1,-1) circle (1pt) node{};
\filldraw[gray] (-1,-1) circle (1pt) node{};
\filldraw[gray] (-2,-1) circle (1pt) node{};

\filldraw[gray] (2,-2) circle (1pt) node{};
\filldraw[gray] (0,-2) circle (1pt) node{};
\filldraw[gray] (1,-2) circle (1pt) node{};
\filldraw[gray] (-1,-2) circle (1pt) node{};
\filldraw[gray] (-2,-2) circle (1pt) node{};

\filldraw[gray] (2,0) circle (1pt) node{};
\filldraw[gray] (0,0) circle (1pt) node{};
\filldraw[gray] (1,0) circle (1pt) node{};
\filldraw[gray] (-1,0) circle (1pt) node{};
\filldraw[gray] (-2,0) circle (1pt) node{};

\filldraw[gray] (2,-3) circle (1pt) node{};
\filldraw[gray] (0,-3) circle (1pt) node{};
\filldraw[gray] (1,-3) circle (1pt) node{};
\filldraw[gray] (-1,-3) circle (1pt) node{};
\filldraw[gray] (-2,-3) circle (1pt) node{};

\filldraw[gray] (2,-4) circle (1pt) node{};
\filldraw[gray] (0,-4) circle (1pt) node{};
\filldraw[gray] (1,-4) circle (1pt) node{};
\filldraw[gray] (-1,-4) circle (1pt) node{};
\filldraw[gray] (-2,-4) circle (1pt) node{};

\end{tikzpicture}
\hspace{0.5cm}
\begin{tikzpicture}[scale=0.7,every node/.style={scale=0.8}, baseline=-110]  
\node   at (-2,2) {$\overline{M}$};

\draw [MidnightBlue,very thick] (0,-1) -- (2,-1);
\draw [MidnightBlue,very thick] (0,-1) -- (0,1);
\draw [MidnightBlue,very thick] (0,-1) -- (0,-3);
\draw [MidnightBlue,very thick] (0,-1) -- (-2,1);

\filldraw[gray] (-2,1) circle (1pt) node{};
\filldraw[gray] (-1,1) circle (1pt) node{};
\filldraw[gray] (0,1) circle (1pt) node{};
\filldraw[gray] (1,1) circle (1pt) node{};
\filldraw[gray] (2,1) circle (1pt) node{};

\filldraw[gray] (-2,0) circle (1pt) node{};
\filldraw[gray] (-1,0) circle (1pt) node{};
\filldraw[gray] (0,0) circle (1pt) node{};
\filldraw[gray] (1,0) circle (1pt) node{};
\filldraw[gray] (2,0) circle (1pt) node{};

\filldraw[gray] (-2,-1) circle (1pt) node{};
\filldraw[gray] (-1,-1) circle (1pt) node{};
\filldraw[gray] (0,-1) circle (1pt) node{};
\filldraw[gray] (1,-1) circle (1pt) node{};
\filldraw[gray] (2,-1) circle (1pt) node{};

\filldraw[gray] (-2,-2) circle (1pt) node{};
\filldraw[gray] (-1,-2) circle (1pt) node{};
\filldraw[gray] (0,-2) circle (1pt) node{};
\filldraw[gray] (1,-2) circle (1pt) node{};
\filldraw[gray] (2,-2) circle (1pt) node{};

\filldraw[gray] (-2,-3) circle (1pt) node{};
\filldraw[gray] (-1,-3) circle (1pt) node{};
\filldraw[gray] (0,-3) circle (1pt) node{};
\filldraw[gray] (1,-3) circle (1pt) node{};
\filldraw[gray] (2,-3) circle (1pt) node{};

\node at (0,-4.5) {\large{The fan of $Z=\FF_1$}};

\node [anchor=west, xshift=3pt]  at (-2,1) {$\Delta_4$};
\node [anchor=west]  at (0,1) {$\Delta_1$};
\node [anchor=south]  at (2,-1) {$\Delta_3$};
\node [anchor=north]  at (0,-3) {$\Delta_2$};

\end{tikzpicture}

\caption{A scaffolding of $f_9$}
\label{LP519468}

\end{figure}

Despite the existence of various lattice points inside this polytope, the only non-vanishing coefficients of $g$ are along its edges, suggesting the use of 2-dimensional struts. Moreover, the shape $\FF_1$ is a pertinent option given the quadrilateral face in the $(z=-2)$--plane. We scaffold the polytope as in Figure \ref{LP519468}, using $(0,0,1)$ as the 0-strut at a basis of $N_U$ cf. Remark \ref{Uneliminated}. We include the projection to the lattice $\overline{N}$ for convenience. The divisors in the shape $Z=\FF_1$ are ordered as in $\overline{M}$ in Figure \ref{LP519468}.


We now follow Algorithm \ref{LIalg} to produce a weight matrix. We have $z=4$, $u=1$, $r=s-u=2$, and $N_U=\langle (0,0,1) \rangle$. We have one 0-strut, i.e. $(D_0,\chi_0)=(0,(0,0,1))$, and two other struts corresponding to $(D_1,\chi_1)=(\Delta_1-\Delta_2+\Delta_3+\Delta_4,-(0,0,1))$ and $(D_2,\chi_2)=(-2\Delta_1+3\Delta_2+\Delta_3-\Delta_4,-2(0,0,1))$. The weight matrix will then be of size $r\times (s+z)=2\times 7$. If the first row corresponds to $D_1$ and the second to $D_2$, the column corresponding to $U$ is $(1,2)$ and we obtain:
\[
\begin{blockarray}{*{7}{c} | c c}
 \BAmulticolumn{2}{c}{} & U & \Delta_1 & \Delta_2 & \Delta_3 & \Delta_4 & L_1 & L_2\\
\BAhline
\begin{block}{*{7}{c} | cc}
  1 & 0 &  1 & 1 & -1 & 1 & 1 & 0 & 2 \\
  0 & 1 &  2 & -2 & 3 & 1 & -1 & 1 & 0\\
\end{block}
\end{blockarray}
\]
with stability condition $\omega=-K-L_1-L_2=(4,4)-(0,1)-(2,0)=(2,3)$. The line bundles are given by homogenising primitive relations in the fan of $Z$, in particular we have that $\Delta_1\Delta_2$ is a section of $L_1$ and $\Delta_3\Delta_4$ is a section of $L_2$. Notice the homogeneity type of the second column of the weight matrix is the same as $L_1$. This means a general section of $L_1$ is linear in this coordinate, thus we can solve for it and eliminate this equation altogether. We obtain an embedding of $X$ as a hypersurface in the rank two toric variety $Y$ given by:

\[
\begin{blockarray}{*{6}{c} | c}
 x_1 & x_2 & x_3 & x_4 & x_5 & x_6 & L\\
\BAhline
\begin{block}{*{6}{c} | c}
  1 & 1 & 1 & -1 & 1 & 1 & 2 \\
  0 & 2 & -2 & 3 & 1 & -1 & 0\\
\end{block}
\end{blockarray}
\]
with the same stability condition as before. The singular charts on $Y$ are $U_{12}$, $U_{14}$, $U_{23}$, $U_{26}$, $U_{45}$ and $U_{46}$, where $U_{ij}=\{x_i, x_j\neq 0\}$. The ambient $Y$ is also singular along two curves $C_1=\{x_1=x_2=x_3=0\}$ and $C_2=\{x_4=x_5=x_6=0\}$. 

The situation is quite tame, as $Y$ is an orbifold - this is the case for many of the ambient spaces in Table \ref{tab1}, yet for very few in Table \ref{tab2}. The main issue is that $C_1$ or $C_2$ could be in the base locus of $L$: once we establish that a general section of $L$ has isolated singularities, they typically coincide with the singularities predicted by the GRDB \cite{BKFano3}. 

A general section $s\in H^0(Y,L)$ is \[x_1x_3x_4x_6+x_3^2x_4x_5+x_3^4x_4^3x_6 +x_4x_6^3+x_2x_3+x_1^2+x_5x_6+ x_1x_3^3x_4^2+x_3^6x_4^4+x_3^2x_4^2x_6^2.\] This is a partial smoothing obtained by adding monomials in $\cO_Y(2,0)$ to the toric equation $x_{1}^{2}+x_{5}x_{6}=0$. Thus $X=V(s)$ is quasi-smooth: it is enough to check this using the Jacobian criterion, making sure its potential singularities all lie in the unstable locus.  

The curve $C_1$ belongs to the chart $U_{46}$, which is of type $\frac{1}{2}(1,1,1,0)_{1235}$. We see that the $\frac{1}{2}(1,1,1)$-action propagates along the $x_5$-axis, i.e. $C_1$ in this chart. If $X$ intersects it transversely, it acquires that precise type of singularity.

When restricting $s$ to $C_1$ we obtain $x_4x_6^3+x_5x_6=0$, thus $X$ contains the origin of $U_{45}$ with multiplicity one (this is a $\frac{1}{4}(1,1,3)$ point), and any other singularity belongs to $U_{46}$, the only other chart containing $C_1$. On $U_{46}$, the equation becomes $1+x_5=0$, whose solution is an isolated point, and (when comparing Jacobian matrices) it is clear that $X$ and $C_1$ intersect transversely here.

Similarly, the curve $C_2$ belongs to the charts $U_{12}$ and $U_{23}$. When restricting $s$ to $C_2$ we are left with $x_2x_3+x_1^2=0$, i.e. a point solving $x_3+1=0$ in the chart $U_{12}$. This chart is $U_{12}=\frac{1}{2}(0,1,1,1)_{x_3x_4x_5x_6}$ so we find another singularity of type $\frac{1}{2}(1,1,1)$.

According to the GRDB prediction, we should obtain another $\frac{1}{3}(1,1,2)$ point, which we find while examining the singular $0$-strata of $Y$: there are only two such points disjoint from $C_1\cup C_2$: $P_1=\{x_1=x_3=x_4=x_5=0\}$ and $P_2=\{x_2=x_3=x_5=x_6=0\}$. The equation of $s$ vanishes at $P_1$, which is the origin of the chart $U_{26}$ (a point of type $\frac{1}{3}(1,1,2)$), and does not vanish at $P_2$. 

\begin{rem}
When analysing the singularity on $C_{1}$, we make a point of examining it in the chart $U_{46}$, where we clearly see that the $\frac{1}{2}(1,1,1)$ singularity propagates. The origin of $U_{45}$ is a more singular point with a $\mu_4$-action, and in this chart the matter of determining which singularity $C_1$ acquires is more subtle. In general it may occur that both charts containing a curve such as $C_1$ have such points at the origin, though this is never true for the examples in this paper.
\end{rem}

\subsubsection{Laurent inversion without point-struts at the basis of $N_U$}
\label{noptstrut}

We explain how to build the weight matrix of $Y$ without the assumption in Remark \ref{Uneliminated}. The key fact is that the lattice of its fan is precisely $\textup{Div}_{T_{\overline{M}}}(Z) \oplus N_U$ (see discussion following Remark 5.4 in \cite{LaurentInversion}). In order to determine the weight matrix we need to specify the rays of this fan and write the sequences \eqref{dual_exact_sequences} for $Y$. 

Suppose $S=\{(D_i,\chi_i)\}_{i=1}^s$ is a scaffolding of $P$ with shape $Z$. Recall from Definition \ref{ScafDef} that this fixes a splitting of the ambient lattice of $P$ as $N=\overline{N}\oplus N_U$, where $\overline{N}$ is the character lattice of $Z$.  Denote by $u=\dim(N_U)$, $\overline{M}=\textup{Hom}(\overline{N},\ZZ)$ and recall that the toric divisors $\Delta_1,\ldots, \Delta_z$ on $Z$ are a $\ZZ$-basis of its divisor group $\textup{Div}_{T_{\overline{M}}}(Z)$.

 \begin{thm}[\cite{LaurentInversion}]
 The fan of $Y$ has ambient lattice $\textup{Div}_{T_{\overline{M}}}(Z)\oplus N_U\simeq \bigoplus\limits_{i=1}^z \ZZ \Delta_i\oplus \ZZ^u$. It has the following $r+s$ rays:
\begin{itemize}
\item $(\Delta_i,0)$, with $i=1,\dots, z$ and
\item $(-D_i,\chi_i)$, since $D_i\in \bigoplus\limits_{i=1}^z \ZZ \Delta_i$ by construction. 
\end{itemize}
 \end{thm}

It is inconvenient to think of high-dimensional toric varieties in terms of fans and cones as opposed to GIT quotients, which is why we usually avoid this description. The rows of the weight matrix output of Algorithm \ref{LIalg} span the kernel of the ray matrix $\rho$ appearing in the dual exact sequences \eqref{dual_exact_sequences} written for $Y$. This kernel can be computed directly because the ray matrix contains an $r \times r$ identity block \textit{if we have $u$ 0-struts at the basis of $N_U$}. This is ultimately the same as imposing that the weight matrix begin with the identity block in the algorithm. 

We illustrate this by computing the weight matrix in the example in Section \ref{basic} using its ray map. Recall $z=4$, $u=1$, $r=s-u=2$, and $N_U=\langle (0,0,1) \rangle$. There is one 0-strut at the origin of $N_U$, i.e. $(D_0,\chi_0)=(0,(0,0,1))$, and two other struts corresponding to $(D_1,\chi_1)=(\Delta_1-\Delta_2+\Delta_3+\Delta_4,-(0,0,1))$ and $(D_2,\chi_2)=(-2\Delta_1+3\Delta_2+\Delta_3-\Delta_4,-2(0,0,1))$, therefore in the basis $\{(\Delta_1,0), \ldots, (\Delta_z,0), (0,(0,0,1))\}$ the rays are:

\[\begin{array}{cccccrr}
\Delta_{1} & \Delta_{2} &\Delta_{3}& \Delta_{4} & -D_0 & -D_1 & -D_2 \\
\hline
1 & 0 & 0 & 0 & 0 & -1 & 2 \\ 
0 & 1 & 0 & 0 & 0 &  1 &-3 \\
0&0 &1 &0 &0 &-1& -1\\
0& 0& 0& 1& 0& -1&  1\\
0& 0& 0& 0& 1& -1& -2
\end{array}\]
Modulo a reordering of the rays, we obtain that its kernel is spanned precisely by the weight matrix in Example \ref{basic}. It is easy to see that the 0-strut condition is not necessary in order to have an identity block. We illustrate this in the example below.

\medskip

\begin{figure}[H]
\begin{center}
\begin{tikzpicture}[scale=0.8,every node/.style={scale=0.8}]
\draw[gray!80, thick] (-1,2.5) -- (0,-1);
\draw[gray!80, thick] (-1,-1) -- (0,-1);
\draw[gray!80, thick] (2,-4) -- (0,-1);
\draw[gray!80, thick] (2,-4) -- (0,-4);

\draw[gray!80, thick] (0.07,0.07) -- (-0.07,-0.07);
\draw[gray!80, thick] (-0.07,0.07) -- (0.07,-0.07);

\draw[black, thick] (-1,2.5) -- (-1,-1);
\draw[black, thick] (-1,2.5) -- (1,-6);
\draw[black, thick] (-1,2.5) -- (2,-4);
\draw[black, thick] (1,-6) -- (2,-4);
\draw[black, thick] (-1,-1) -- (0,-4);
\draw[black, thick] (-1,-1) -- (0,-4);
\draw[black, thick] (1,-6) -- (0,-4);
\draw[black, thick] (-1,-1) -- (0,-4);

\filldraw[black] (-1,2.5) circle (2pt) node[anchor=west] {$(1,-1,2)$};
\filldraw[black] (1,-6) circle (2pt) node[anchor=west] {$(-2,3,-3)$};
\filldraw[black] (2,-4) circle (2pt) node[anchor=west] {$(-2,1,-3)$};
\filldraw[black] (0,-4) circle (2pt) node[anchor=east] {$(0,1,-3)$};
\filldraw[black] (-1,-1) circle (2pt) node[anchor=east] {$(1,0,-1)$};
\filldraw[gray] (0,-1) circle (2pt) node[anchor= west]{\textcolor{black}{$(0,0,-1)$}};
\filldraw[black] (0,0) circle (0pt) node[anchor=west] {\textcolor{black}{$(0,0,0)$}};

\filldraw[black] (1.5,-5) circle (2pt) node[anchor=south east]{};
\filldraw[black] (0.5,-5) circle (2pt) node[anchor=south east]{};
\filldraw[gray] (1,-4) circle (1.5pt) node[anchor=south east]{};
\filldraw[gray] (1,-3) circle (1.5pt) node[anchor=south east]{};
\filldraw[gray] (0,-3) circle (1.5pt) node[anchor=south east]{};

\end{tikzpicture}
\end{center}
\caption{Newton polytope of $h$}
\label{P512391}
\end{figure}
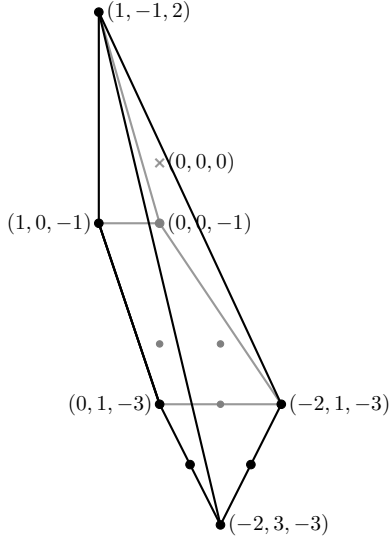

Consider variety \#21 in Table \ref{tab1}, which is a codimension $21$ $\mathbb{Q}-$Fano threefold with two singular points, of type $\frac{1}{3}(1,1,2)$ and $\frac{1}{4}(1,1,3)$ respectively. We construct this as a deformation of the toric canonical P$512391$ to the $\QQ$-Fano Q$38917$. The polytope of P$512391$ supports a single rigid MMLP: \[f = x + y + z + \frac{1}{xyz} + \frac{2}{xyz^2} + \frac{1}{xy^3z^4} + \frac{1}{x^2y^3z^3} + \frac{1}{x^2y^3z^4}\]  
whose 1-step mutation (modulo a change of basis) 
\[ h = \frac{xz^{2}}{y}+\frac{x}{z}+\frac{1}{z}+\frac{y}{z^{3}}+\frac{2y}{xz^{3}}+ \frac{y}{x^{2}z^{3}}+\frac{2y^{2}}{xz^{3}}+\frac{2y^{2}}{x^{2}z^{3}}+\frac{y^{3}}{x^{2}z^{3}}\] is supported on the Newton polytope $P$ in Figure \ref{P512391}.


Notice the bottom triangle suggests the shape $Z$ is a variety admitting a contraction to $\PP^2$ (the polytope of one of its nef divisors is this triangle, and $Z$ is a tower of projective bundles). As in the previous case, the coefficients of $h$ suggest a $2$-dimensional shape as $h$ is entirely supported on the edges of $P$.  We aim for a complete intersection, therefore the two options are $Z=\PP^{2}$ and $Z=\FF_{1}$. 

\begin{rem} The key issue here is that no single point proves to be a good candidate as the basis of $N_U$: the best option is $(1,-1,2)$, however $N$ does not decompose as $\overline{N}\oplus N_{U}$, where $\overline{N}$ is the 2-dimensional ambient lattice of the bottom triangle, and $N_U$ the lattice spanned by $(1,-1,2)$ (since $\overline{N}\oplus N_{U}$ is of index two inside $N$).  The cases of $N_U=\langle(0,0,-1)\rangle$ and $N_U=\langle(1,0,-1)\rangle$ are unsuccessful - here the lattice decomposes as needed but when analysing the weight matrix obtained as in Section \ref{basic} we obtain that $X_{P}$ does not deform to an orbifold. It is in general desirable to use a minimal amount of struts, and this is exactly why the example below works. \end{rem}
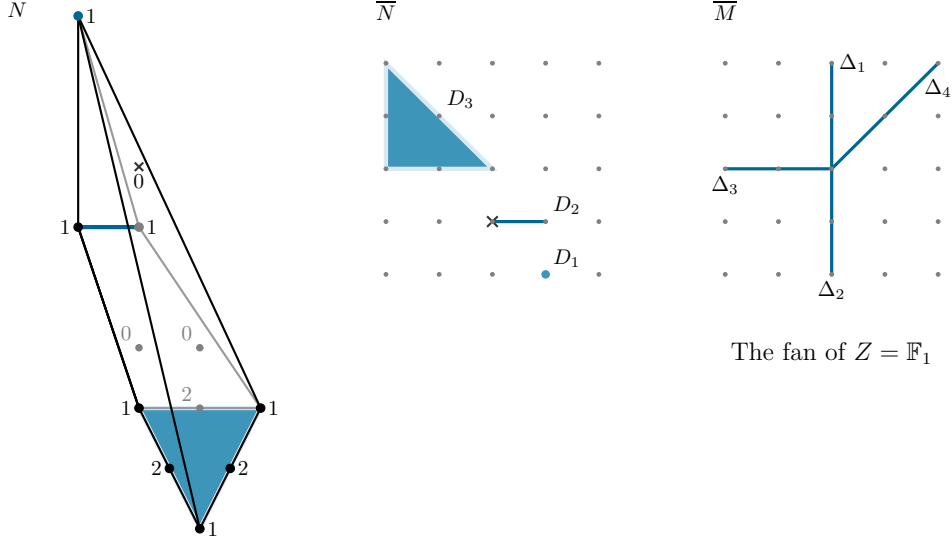
\begin{figure}[H]

\begin{minipage}{0.3\textwidth}
\vspace{-3.5cm}
\begin{center}
\begin{tikzpicture}[scale=0.8,every node/.style={scale=0.8}]
\node   at (-2,2.6) {$N$};
\filldraw[color=MidnightBlue!10, fill=MidnightBlue!60, ultra thick] (1,-6) -- (2,-4) -- (0,-4) -- cycle; 
\draw[gray!80, thick] (-1,2.5) -- (0,-1);
\draw[gray!80, thick] (-1,-1) -- (0,-1);
\draw[gray!80, thick] (2,-4) -- (0,-1);
\draw[gray!80, thick] (2,-4) -- (0,-4);

\draw[MidnightBlue, ultra thick] (0,-1) -- (-1,-1);

\draw[black!80, thick] (0.07,0.07) -- (-0.07,-0.07);
\draw[black!80, thick] (-0.07,0.07) -- (0.07,-0.07);

\draw[black, thick] (-1,2.5) -- (-1,-1);
\draw[black, thick] (-1,2.5) -- (1,-6);
\draw[black, thick] (-1,2.5) -- (2,-4);
\draw[black, thick] (1,-6) -- (2,-4);
\draw[black, thick] (-1,-1) -- (0,-4);
\draw[black, thick] (-1,-1) -- (0,-4);
\draw[black, thick] (1,-6) -- (0,-4);
\draw[black, thick] (-1,-1) -- (0,-4);

\filldraw[MidnightBlue] (-1,2.5) circle (2pt) node[anchor=west] {\textcolor{black}{$1$}};
\filldraw[black] (1,-6) circle (2pt) node[anchor=west] {$1$};
\filldraw[black] (2,-4) circle (2pt) node[anchor=west] {$1$};
\filldraw[black] (0,-4) circle (2pt) node[anchor=east] {$1$};
\filldraw[black] (-1,-1) circle (2pt) node[anchor=east] {$1$};
\filldraw[gray] (0,-1) circle (2pt) node[anchor= west]{\textcolor{black}{$1$}};
\filldraw[black] (0,0) circle (0pt) node[anchor=north] {\textcolor{black}{$0$}};

\filldraw[black] (1.5,-5) circle (2pt) node[anchor=west]{$2$};
\filldraw[black] (0.5,-5) circle (2pt) node[anchor=east]{$2$};
\filldraw[gray] (1,-4) circle (1.5pt) node[anchor=south east]{$2$};
\filldraw[gray] (1,-3) circle (1.5pt) node[anchor=south east]{$0$};
\filldraw[gray] (0,-3) circle (1.5pt) node[anchor=south east]{$0$};

\end{tikzpicture}
\end{center}
\end{minipage}
\begin{minipage}[t]{0.3\textwidth}
\begin{center}
\begin{tikzpicture}[scale=0.7,every node/.style={scale=0.8}, baseline=-110]  

\node   at (-2,2) {$\overline{N}$};
\draw[MidnightBlue,very thick] (1,-2) -- (0,-2);
\filldraw[color=MidnightBlue!10, fill=MidnightBlue!60, ultra thick] (-2,1) -- (-2,-1) -- (0,-1) -- cycle; 
\draw[black!80, thick] (0.1,-0.9-1) -- (-0.1,-1.1-1);
\draw[black!80, thick] (-0.1,-0.9-1) -- (0.1,-1.1-1);

\filldraw[gray] (2,1) circle (1pt) node{};
\filldraw[gray] (0,1) circle (1pt) node{};
\filldraw[gray] (1,1) circle (1pt) node{};
\filldraw[gray] (-1,1) circle (1pt) node{};
\filldraw[gray] (-2,1) circle (1pt) node{};


\filldraw[gray] (2,-1) circle (1pt) node{};
\filldraw[gray] (0,-1) circle (1pt) node{};
\filldraw[gray] (1,-1) circle (1pt) node{};
\filldraw[gray] (-1,-1) circle (1pt) node{};
\filldraw[gray] (-2,-1) circle (1pt) node{};

\filldraw[gray] (2,-2) circle (1pt) node{};
\filldraw[gray] (0,-2) circle (1pt) node{};
\filldraw[gray] (1,-2) circle (1pt) node{};
\filldraw[gray] (-1,-2) circle (1pt) node{};
\filldraw[gray] (-2,-2) circle (1pt) node{};

\filldraw[gray] (2,0) circle (1pt) node{};
\filldraw[gray] (0,0) circle (1pt) node{};
\filldraw[gray] (1,0) circle (1pt) node{};
\filldraw[gray] (-1,0) circle (1pt) node{};
\filldraw[gray] (-2,0) circle (1pt) node{};

\filldraw[gray] (2,-3) circle (1pt) node{};
\filldraw[gray] (0,-3) circle (1pt) node{};
\filldraw[MidnightBlue!60] (1,-3) circle (2pt) node{};
\filldraw[gray] (-1,-3) circle (1pt) node{};
\filldraw[gray] (-2,-3) circle (1pt) node{};


\node [anchor=south west]  at (-1,0) {$D_3$};
\node [anchor=south west]  at (1,-2){$D_2$};
\node [anchor=south west]  at (1,-3) {$D_1$};

\end{tikzpicture}
\end{center}

\end{minipage}
\begin{minipage}[t]{0.3\textwidth}
\begin{center}
\begin{tikzpicture}[scale=0.7,every node/.style={scale=0.8}, baseline=-110]  
\node   at (-2,2) {$\overline{M}$};

\draw [MidnightBlue,very thick] (0,-1) -- (2,1);
\draw [MidnightBlue,very thick] (0,-1) -- (0,1);
\draw [MidnightBlue,very thick] (0,-1) -- (0,-3);
\draw [MidnightBlue,very thick] (0,-1) -- (-2,-1);

\filldraw[gray] (-2,1) circle (1pt) node{};
\filldraw[gray] (-1,1) circle (1pt) node{};
\filldraw[gray] (0,1) circle (1pt) node{};
\filldraw[gray] (1,1) circle (1pt) node{};
\filldraw[gray] (2,1) circle (1pt) node{};

\filldraw[gray] (-2,0) circle (1pt) node{};
\filldraw[gray] (-1,0) circle (1pt) node{};
\filldraw[gray] (0,0) circle (1pt) node{};
\filldraw[gray] (1,0) circle (1pt) node{};
\filldraw[gray] (2,0) circle (1pt) node{};

\filldraw[gray] (-2,-1) circle (1pt) node{};
\filldraw[gray] (-1,-1) circle (1pt) node{};
\filldraw[gray] (0,-1) circle (1pt) node{};
\filldraw[gray] (1,-1) circle (1pt) node{};
\filldraw[gray] (2,-1) circle (1pt) node{};

\filldraw[gray] (-2,-2) circle (1pt) node{};
\filldraw[gray] (-1,-2) circle (1pt) node{};
\filldraw[gray] (0,-2) circle (1pt) node{};
\filldraw[gray] (1,-2) circle (1pt) node{};
\filldraw[gray] (2,-2) circle (1pt) node{};

\filldraw[gray] (-2,-3) circle (1pt) node{};
\filldraw[gray] (-1,-3) circle (1pt) node{};
\filldraw[gray] (0,-3) circle (1pt) node{};
\filldraw[gray] (1,-3) circle (1pt) node{};
\filldraw[gray] (2,-3) circle (1pt) node{};

\node at (0,-4.5) {\large{The fan of $Z=\FF_1$}};

\node [anchor=north]  at (-2,-1) {$\Delta_3$};
\node [anchor=west]  at (0,1) {$\Delta_1$};
\node [anchor=north, yshift=-4]  at (2,1) {$\Delta_4$};
\node [anchor=north]  at (0,-3) {$\Delta_2$};

\end{tikzpicture}

\end{center}

\end{minipage}
\captionof{figure}{A scaffolding of $f_{2}$}
\label{LP512391}

\end{figure}

Consider the scaffolding in Figure \ref{LP512391}, with shape $\FF_{1}$ (so $z=4$ just as before) and three struts $S_i$. Decompose $N=\overline{N}\oplus N_{U}$ with $N_{U}=\langle (0,0,1) \rangle$ and follow the above method for producing the ray map of $Y$: the ambient space is of dimension $z+u=5$ and its fan is given by $z+s=4+3=7$ rays in the lattice $\bigoplus_i \ZZ \Delta_{i}\oplus N_{U}$. Proceeding exactly as in the previous case, we obtain:
\[\begin{array}{ccccrrr}
\Delta_{1} & \Delta_{2} &\Delta_{3}& \Delta_{4} & -D_1 & -D_2 & -D_3 \\
\hline
1 & 0& 0& 0&  1&0&-3\\
0 & 1& 0& 0& -1&0& 1\\
0& 0& 1& 0&  1&0&-2\\
0& 0& 0& 1& 0&-1&-1\\
0& 0& 0& 0& -2& 1& 3
\end{array} \] and the weight matrix is simply the two-dimensional kernel of the ray map. 

The homogeneity type of the two line bundles which give $X_{P}$ are again given by the relations in the fan of $Z$, written in terms of the Cox coordinates. One obtains the following, modulo a reordering of the rays:
\[\begin{blockarray}{*{7}{c} | c c}
 x_1 & x_2 & x_3 & x_4 & x_5 & x_6 & x_7 & L_{1} & L_{2}\\
\BAhline
\begin{block}{*{7}{c} | c c}
  1 & 0 & -3 & 3 & -1 & 2 &  -2 & 2 & 0 \\
  0 & 1& 2 & -1 & 1 & -1 & 2 & 0 & 1\\
\end{block}
\end{blockarray}
\]
with the stability condition $\omega=(-2,3)$. We can eliminate $L_{2}$ by solving for $x_{2}$ and the rest of the verification is almost identical to the case in Section \ref{basic}, down to the fact that the ambient space is an orbifold which is singular along two curves (in this case, $C_{1}=\{x_{4}=x_{5}=x_{6}=0\}$ and $C_{2}=\{x_{1}=x_{3}=x_{6}=0\}$). We leave the details of this calculation to the reader.

\newpage

\newpage
\newgeometry{top=1cm, bottom=1cm,left=1cm, right=1cm}
\pagestyle{empty}
\begin{landscape}
\begin{center}
\title{\textbf{Table 2: Systematic constructions of high codimension examples}}
\maketitle
\arrayrulecolor{Gray}
\setlength{\arrayrulewidth}{0.2mm}
\rowcolors{1}{Gray!10}{Periwinkle!20}
\begin{longtable}[t]{| m{0.3cm} | c | c | >{\centering\arraybackslash}m{2.2cm}| >{\centering\arraybackslash}m{0.5cm} | >{\centering\arraybackslash}m{3.3cm} | c | m{3cm}|  >{\centering\arraybackslash}m{1.6cm} | }
\hline
 \label{tab1}
 & \textbf{$\QQ$-Fano} & \textbf{Polytope} & \textbf{Singularities } & \textbf{MG \#} & \textbf{Rigid MMLP} & \textbf{Weight Matrix and Bundles} & \textbf{Period Sequence} &  \textbf{Shape}\\
 \hline
 \hline
 
  \multicolumn{9}{|l|}{\textbf{Codim 23}}\\
  
  \hline 


1 & Q38950 & P519483 &  $\frac{1}{4}(1,1,3)$, $\frac{1}{5}(1,2,3)$ & 2 & $xy^8z^5 + 2xy^6z^4 + xy^4z^3 + x + 3y^4z^3 + 4y^2z^2 + y + z + 3z/x 
+ 2/(xy^2) + 1/(x^2y^4z)$ & \begin{math}\begin{array}{cccccc|c} 1&0&2&-2&1&1&2\\ 0&1&-1&3&1&-1&0 \end{array}\end{math} & [1, 0, 0, 0, 24, 0, 0, 0, 2520, 7560, 0] & $\FF_1$ + 0-st. \\

\hline

2 & Q38985 & P519424 &  $\frac{1}{2}(1,1,1)$, $\frac{1}{3}(1,1,2)$,  $\frac{1}{5}(1,2,3)$ & 7 & $x + y + z + 1/y + z/(xy) + 1/(xy^5z^2) + 3/(x^2y^5z) + 3/(x^3y^5) + 
z/(x^4y^5)$ & \begin{math}\begin{array}{cccccc|c} 0&1&0&1&1&0&1 \\ 1&2&2&1&-3&3&3 \end{array}\end{math} & 
  [1, 0, 2, 0, 6, 0, 20, 0, 70, 7560, 252] & $\PP^1\times\PP^1$ + 0-st. \\ 

\hline
3 & Q38989 & P519468 & $2\times\frac{1}{2}(1,1,1)$, $\frac{1}{3}(1,1,2)$, $\frac{1}{4}(1,1,3)$ & 2 &  $x + y + z + y/(x^4z^3) + 2/(x^4z^2) + 1/(x^4yz) + y/(x^5z^4) + 3/(x^5z^3) 
+ 3/(x^5yz^2) + 1/(x^5y^2z)$ & \begin{math}\begin{array}{cccccc|c} 1&1&1&-1&1&1&2 \\ 0&2&-2&3&1&-1&0 \end{array}\end{math} & [1, 0, 0, 0, 0, 0, 0, 420, 0, 7560, 0] & $\FF_1$ + 0-st. \\

\hline 
 4 & Q39260 & P516023 &  $\frac{1}{5}(1,2,3)$ & 9 & $x + y + z + y^2/(x^2z^3) + 2y^2/(x^2z^4) + y^2/(x^2z^5) + y/(x^2z^2) + 
4y/(x^2z^3) + 3y/(x^2z^4) + 2/(x^2z^2) + 3/(x^2z^3) + 1/(x^2yz^2)$ & \begin{math}\begin{array}{cccccc|c}  1&0&1&-1&1&1&2\\ 0&1&-1&3&2&-2&0 \end{array}\end{math} & [1, 0, 0, 0, 0, 60, 360, 0, 0, 0, 18900]
 & $\FF_1$ + 0-st. \\
 
\hline
 
5 & Q39266 & P541268 & $2 \times \frac{1}{3}(1,1,2)$ & 2 & $xyz^2 + x + y^2z + y + 2/(xz) + 1/(x^2y^2z^3)$ & \begin{math}\begin{array}{cccccc|c} 0&1&0&1&1&-1&0 \\ 1&-4&2&-2&-1&5&4 \end{array} \end{math}
 & [1, 0, 0, 0, 0, 120, 0, 420, 0, 0, 88200]
 & $\FF_2$ + 0-st. \\

6 & & P541702 & $2 \times \frac{1}{3}(1,1,2)$ & 2 & $xyz^3 + x + 2yz + y + 1/y + y/(xz)
$ & \begin{math}\begin{array}{cccccc|c} 1&0&2&-2&3&1&4 \\ 0&1&0&1&0&0&0 \end{array} \end{math}
 & [1, 0, 2, 0, 6, 120, 20, 1680, 70, 15120, 88452 ] & $\FF_2$ + 0-st. \\

\hline

7 & Q39329 & P429996 & $\frac{1}{2}(1,1,1)$, $\frac{1}{3}(1,1,2)$  & 3 & $xy^3z^2 + x + 2yz + y + z + 1/(xy) + 2/(xy^2) + 1/(x^2y^4z)$  &  \begin{math}\begin{array}{cccccc|c} 1&0&1&-1&1&-1&0 \\ 0&1&-1&2&-1&3&2 \end{array}\end{math} & [1, 0, 0, 6, 24, 0, 90, 1680, 2520, 1680, 75600] & $\PP^1\times \PP^1$ + 0-st. \\ 

\hypertarget{eight}{8} &  & P516841 &   & 1 & $x + y + z + 1/(y^2z^3) + 1/(xyz^3) + 1/(xy^2z^2) + 1/(xy^2z^3)$ & \begin{math}\begin{array}{cccccc|c} 1&0&0&0&0&1&1 \\ 0&1&1&2&3&-5&1 \end{array}\end{math} & [1, 0, 0, 0, 0, 0, 360, 420, 0, 0, 0] & $\PP^3$ \\

\hline

 \multicolumn{9}{|l|}{\textbf{Codim 22}}\\
 
\hline
 
9 & Q38746 & P517792 & $\frac{1}{4}(1,1,3)$, $\frac{1}{5}(1,1,4)$ & 11 & $xy^4z^5 + xy^3z^4 + x + 3y^3z^3 + 2y^2z^2 + y + 1/y + 3y^2z/x + y/x + 
y/(x^2z)$ & \begin{math}\begin{array}{ccccccc|c} 1&0&-1&-2&2&-1&-1&-2 \\ 0&0&0&-2&3&-1&1&0 \\ 0&1&1&2&-2&1&1&2 \end{array} \end{math} & [1, 0, 2, 6, 6, 60, 110, 420, 5110, 4200, 94752] & $\FF_1$ + 0-st.\\

10 & & P518340 & & 1& $xy^4z^5 + xy^3z^4 + x + 2yz^2 + y + z + 1/y + 1/(xy^2z)$ & \begin{math}\begin{array}{ccccccc|c} 1&0&0&3&-3&1&2&3\\0&1&1&0&0&0&0&0 \\ 0&0&1&-1&2&1&-1&0 \end{array}\end{math} & [1, 0, 2, 0, 6, 60, 20, 840, 3430, 7560, 94752] & $\FF_1$ + 0-st. \\
 
\hline

\hypertarget{eleven}{11} & Q38775 & P519409 & $\frac{1}{2}(1,1,1)$, $\frac{1}{3}(1,1,2)$, $\frac{1}{5}(1,1,4)$ & 3 & $xyz^2 + x + 2x/(yz) + x/(y^2z^3) + x/(y^3z^4) + y + 2z/y + 1/(y^2z) + 
2/(y^3z^2) + 1/(xy^3)$ & \begin{math} \begin{array}{cccccc|c} 1&3&-3&0&1&2&3\\0&-4&5&1&-1&-3&-3 \end{array}\end{math} & [1, 0, 0, 0, 0, 60, 0, 0, 3360, 0, 18900] & \hyperlink{Z1}{$Z_1$} \\ 

\hline

12 & Q38791 & P519951 & $2 \times \frac{1}{3}(1,1,2)$, $\frac{1}{4}(1,1,3)$ & 1 & $x + y + z + 1/(xz) + 1/(x^2y^2z^3) + 1/(x^3y^2z^2)$ & \begin{math}\begin{array}{cccccc|c} 1&0&0&2&1&1&2\\ 0&1&2&5&3&3& 6 \end{array} \end{math} &  [1, 0, 0, 6, 0, 0, 90, 0, 3360, 1680, 0] & $\FF_1$ (no~0-st.) \\

\hline

13& Q38946 & P427943 &  $\frac{1}{3}(1,1,2)$, $\frac{1}{5}(1,1,4)$ & 1 & $x + y^3z + y + z + y^3z^2/x + 1/x + 1/(xyz)$ & \begin{math}\begin{array}{ccccccc|c} 0&0&1&1&0&0&1&1\\ 1&0&-3&-1&1&1&-2&-1 \\ 0&1&-8&-3&3&2&-4&-2 \end{array}\end{math} & [1, 0, 2, 0, 30, 0, 380, 0, 5950, 7560, 101052] & $\FF_1$ + 0-st \\

\hline

14 & Q38948 & P428047 &  $\frac{1}{2}(1,1,1)$, $\frac{1}{3}(1,1,2)$, $\frac{1}{4}(1,1,3)$ & 1 & $x + xz^3/y^2 + y^2/z + y + y/(xz) + y/(xz^2) + y^3/(x^2z^3)$ & \begin{math}\begin{array}{ccccccc|c} 1&0&1&1&-1&1&1&2\\ 0&1&-1&-1&1&-1&0&-1\\ 0&0&-2&-1&2&-3&2&-1 \end{array}\end{math} & [1, 0, 0, 0, 24, 0, 0, 420, 2520, 7560, 0] &  $\FF_1$ + 0-st \\

\hline
15 & Q39259 & P505841 &  $\frac{1}{2}(1,1,1)$, $\frac{1}{3}(1,1,2)$ & 2 & $x^3y/z^2 + 2xy/z + x + y + z + 1/y + y/x$ & \begin{math}\begin{array}{ccccccc|c} 1&0&0&2&-2&1&3&4\\ 0&1&0&1&-1&0&2&2\\ 0&0&1&0&1&0&0&0 \end{array}\end{math} & [1, 0, 2, 6, 6, 180, 110, 2100, 8470, 19320, 258552] & $\FF_2$ + 0-st.   \\

16 & & P507415 & & 2 & $x^2yz^3 + xyz^2 + x + y + 2z + 1/x + 1/(x^2yz)$ & \begin{math}\begin{array}{cccccc|c} 1&0&2&-2&1&1&2\\0&1&0&1&0&1&1 \end{array}\end{math} & [1, 0, 2, 0, 6, 120, 20, 2100, 70, 22680, 88452] & $\FF_1$ + 0-st.   \\

\hline
17 & Q39277 & P520149 &  $3 \times \frac{1}{2}(1,1,1)$ & 1 & $xyz^2 + x + yz + y + z/x + 1/(xyz)
$ & \begin{math}\begin{array}{cccccc|c} 0&-2&0&1&1&-1&0\\ 1&3&1&-1&-1&3&2 \end{array}\end{math} & [1, 0, 0, 6, 0, 120, 90, 0, 6720, 1680, 88200] & $\FF_1$ + 0-st. \\

\hline

18& Q39321 & P255545 & $\frac{1}{3}(1,1,2)$ & 2 & $xy^3/z + xy^2/z + x + y + y/z + z + 1/(xy)$
 &  \begin{math}\begin{array}{ccccccc|c} 1&0&0&0&-2&1&1&0 \\ 0&1&0&0&1&0&0& 1 \\ 0&0&1&1&1&-1&1&1  \end{array}\end{math}  &  [1, 0, 0, 6, 24, 60, 90, 1680, 5880, 16800, 94500] & $\PP^2$ + 0-st.\\

\hline

19 & Q39328 & P254809 &  $2\times \frac{1}{2}(1,1,1)$ & 1 & $x^2z^3/y^2 + x^2z^4/y^3 + xz + x + xz^2/y + xz^3/y^2 + y + y/z + y/(xz^2)$ & \begin{math}\begin{array}{ccccccc|c} 1&0&0&-1&1&2&-1&1\\ 0&0&1&0&1&1&0&1\\ 0&1&1&1&-1&0&1&1 \end{array} \end {math} & [1, 0, 2, 6, 30, 60, 470, 2100, 7630, 42000, 195552] & $\PP^1 \times \PP^1$ + 0-st.
\\

\hline

\multicolumn{9}{| l |}{\textbf{Codim 21}}\\

\hline 


\hline

20 & Q38745 & P515798 &  $2\times \frac{1}{4}(1,1,3)$ & 2 & $x + y + z + 1/y + 1/(xy) + 1/(xy^2z) + 2/(x^2y^3z^2) + 1/(x^3y^5z^4)$ & \begin{math} \begin{array}{ccccccc|c} 1&0&0&2&-2&-1&3&2 \\ 0&1&1&0&0&0&0& 0 \\ 0&0&1&-1&2&1&-1&0 \end{array} \end{math} & [1, 0, 2, 6, 6, 120, 110, 1260, 8470, 11760, 189252] & $\FF_1$ + 0-st. \\
\hline

21 & Q38764 & P429611 &  $2\times \frac{1}{2}(1,1,1)$, $\frac{1}{5}(1,1,4)$ & 2 & $xy^4z^5 + 2xy^3z^4 + xy^2z^3 + x + 2yz^2 + y + 2z + 1/(xy^2z)$ & \begin{math}\begin{array}{cccccc|c} 1&0&0&1&-3&2&0\\0&1&1&1&2&-1&2 \end{array}\end{math} & [1, 0, 0, 0, 0, 120, 0, 0, 3360, 0, 88200] & $\PP^2$ + 0-st.\\

\hline
 
22 & Q38772 & P428621 &  $\frac{1}{2}(1,1,1)$, $\frac{1}{3}(1,1,2)$, $\frac{1}{4}(1,1,3)$ & 2 & $xy^3z^4 + x + 2yz^2 + y + z + 1/(xy) + 1/(xy^2z)$ & \begin{math}\begin{array}{cccccc|c} 1&0&2&-2&-1&3&2 \\ 0&1&-1&2&1&-1&0 \end{array}\end{math} & [1, 0, 0, 6, 0, 60, 90, 0, 6720, 1680, 18900] & $\FF_1 $ + 0-st.\\

\hline

23 & Q38917 & P512391 &  $\frac{1}{3}(1,1,2)$, $\frac{1}{4}(1,1,3)$ & 2 & $x^3y^4z + x^2y^4z + 2xy^2z + x + y + z + 2y/x + 2/(x^2y) + 1/(x^4y^2z)$ & \begin{math}\begin{array}{cccccc|c} 1&-3&-1&3&2&-2&2 \\ 0&2&1&-1&-1&2&0 \end{array}\end{math} & [1, 0, 0, 0, 24, 120, 0, 0, 2520, 37800, 88200] &  $\FF_1 $ (no~0-st.) \\

\hline 

24 & Q38945 & P399587 &  $\frac{1}{2}(1,1,1)$, $2 \times \frac{1}{3}(1,1,2)$ & 1 & $x^3z^4/y^3 + x^2z + x + xz^3/y^2 + xz^2/y^2 + y + y/z + y/(xz^2) $ & \begin{math}\begin{array}{cccccc|c} 1&0&1&-3&-1&-2& -2 \\ 0&1&1&1&1&1&2 \end{array}\end{math} &  [1, 0, 0, 12, 0, 0, 540, 420, 0, 33600, 75600] & $\FF_1$ (no~0-st.) \\
\hline

25 & Q39250 & P428844 &  $\frac{1}{3}(1,1,2)$ & 1 & $x^2z^2/y + 2xz + x + xz^3/y^2 + y + y/z + 1/z + y/(xz^2)$ & \begin{math}\begin{array}{cccccc|c} 1&0&1&-1&1&1&2\\ 0&1&-1&2&1&-1&0 \end{array}\end{math} & [1, 0, 0, 6, 48, 0, 90, 2520, 11760, 1680, 100800] & $\FF_1$ + 0-st.\\

26 & & P430410 & & 1 & $xyz^2 + xyz + x + y + z/x + 1/x + 1/(xyz)$ & \begin{math} \begin{array}{cccccc|c} 1&2&1&1&1&-1&2 \\ 0&1&1&1&0&1&2 \end{array} \end{math} & [1, 0, 4, 0, 36, 120, 400, 4200, 4900, 105840, 151704]
 & $\PP^1 \times \PP^1$ (no~0-st.) \\

\hline

27 & Q39258 & P401207 &   $2\times \frac{1}{2}(1,1,1)$ & 3 & $x + y + z + 1/y + 1/(xy) + 2/(xy^2z) + 1/(xy^3z^2)$ & \begin{math}\begin{array}{cccccc|c} 1&0&1&-1&-1&2&1 \\ 0&1&0&1&1&0&1 \end{array}\end{math} & [1, 0, 2, 6, 6, 180, 110, 2520, 8470, 26880, 283752] & $\FF_1$ + 0-st.\\

\hline
 
28 & Q39320 & P430098 &  $\frac{1}{2}(1,1,1)$ & 1 & $x^3yz + x^2yz + x + y + z + 1/(x^2z) + 1/(x^2y)$ & \begin{math}\begin{array}{cccccc|c} 1&0&1&2&3&1&4 \\ 0&1&0&0&1&0&1 \end{array}\end{math} & [1, 0, 0, 0, 48, 60, 0, 0, 11760, 30240, 18900] & $\PP^1 \times \PP^1$ (no~0-st.)\\

 \hline
 
 \multicolumn{9}{|l|}{\textbf{Codim 20}}\\
 
 \hline
 
 29 & Q37923 & P544136 &  $\frac{1}{2}(1,1,1)$, $\frac{1}{4}(1,1,3)$, $\frac{1}{5}(1,1,4)$ & 5 & $x + y + z + 1/(yz^2) + 1/(xz^2) + 1/(xy^2z^3) + 3/(x^2yz^3) + 3/(x^3z^3) + y/(x^4z^3)$ & \begin{math} \begin{array}{cccccc|c} 1&2&1&-1&3&1&4 \\ 0&-3&-2&3&-5&-1&-6 \end{array} \end{math} & [1, 0, 0, 0, 24, 0, 0, 2100, 2520, 0, 0] & $\FF_2$ (no 0-st) \\

 \hline
 30 & Q37992 & P519452 &   $\frac{1}{2}(1,1,1)$, $\frac{1}{4}(1,1,3)$, $\frac{1}{5}(1,1,4)$  & 2 & $x + x/(yz) + y + z + z/(xy^2) + 2/(xy^3) + 1/(xy^4z) + z^3/(x^2y) + 4z^2/(x^2y^2) + 3z/(x^2y^3) + 2z^4/(x^3y) + 3z^3/(x^3y^2) + z^5/(x^4y)$ & \begin{math} \begin{array}{cccccc|c} 1&0&2&-1&1&1&2 \\ 0&1&-1&3&1&-1&0 \end{array}\end{math} & [1, 0, 0, 0, 0, 60, 0, 1680, 0, 7560, 18900] & $\FF_1$ + 0-st. \\
  
 31 & & P544275 & & 4 & $x^3yz^5 + 2xyz^2 + x + y + 2z^2 + y/(xz) + 2/(x^2z) + 1/(x^3yz)$ & \begin{math} \begin{array}{cccccc|c} 0&3&4&1&3&2&6 \\1&0&1&0&1&-1&0 \end{array} \end{math} & [1, 0, 0, 0, 0, 0, 0, 2100, 0, 15120, 0] & $\PP^2$ (no~0-st.)\\
 
\hline 

32 & Q38005 & P544272 &   $2 \times \frac{1}{2}(1,1,1)$, $2\times \frac{1}{4}(1,1,3)$ & 2 & $xz^2 + 2xz + x + y + z^2/y + 1/x + z^2/(x^3y^3) + 2/(x^4y^2) + 1/(x^5yz^2)$ & \begin{math} \begin{array}{cccccc|c} 0&3&4&2&3&1&6 \\ 1&0&0&0&1&0&0 \end{array} \end{math} & [1, 0, 2, 0, 6, 0, 20, 2100, 70, 37800, 252] & $\FF_1$ (no~0-st.)\\

\hline

33 & Q38329 & P480334 &  $\frac{1}{2}(1,1,1)$, $\frac{1}{3}(1,1,2)$, $\frac{1}{5}(1,1,4)$ & 4 &  $xyz^2 + x + x/(yz^3) + y + 2/(yz^2) + 2/(y^3z^5) + 1/(xyz) + 1/(xyz^2) + 2/(xy^3z^4) + 1/(xy^5z^7)$  & \begin{math} \begin{array}{ccccccc|c} 0&0&3&5&1&-1&2&6 \\ 1&0&2&4&0&-1&1&4 \\ 0&1&-3&-5&-1&2&-2&-6 \end{array}\end{math} & 
[1, 0, 2, 0, 6, 60, 380, 1260, 10150, 15120, 170352] & $\FF_1$ (no~0-st.)\\

\hline 

34 & & P480338 & & 5 & $x^4y^7z^2 + 2x^3y^4z + 2x^2y^4z + x^2y^3z + x^2y + 2xy + x + y + z 
+ 1/z + 1/(xy^3z^2)$ & \begin{math}\begin{array}{ccccccc|c} 0&0&3&1&5&-1&2&6 \\ 1&0&-3&-1&-5&2&-2&-6 \\ 0&1&0&-1&1&0&0&0 \end{array} \end{math} & [1, 0, 2, 0, 6, 0, 380, 420, 10150, 15120, 151452] & $\FF_1$ (no~0-st.) \\

\hline

35 & Q38517 & P518018 &  $2 \times \frac{1}{4}(1,1,3)$ & 3 & $x + y + y/z + z + 1/z + 2/z^2 + 1/(yz^3) + 2y/x^2 + 2/(x^2z) + yz/x^4$ & \begin{math}\begin{array}{ccccccc|c} 0&0&-2&1&-1&0&1&0 \\ 1&0&2&-1&1&1&-1&0\\ 0&1&4&0&2&1&-1&2 \end{array} \end{math}  & [1, 0, 2, 6, 30, 60, 470, 2940, 7630, 64680, 271152] & $\FF_1$ + 0-st.\\

\hline

36 & Q38737 & P511389 &  $\frac{1}{3}(1,1,2)$, $\frac{1}{4}(1,1,3)$ & 31 & $x^3yz^4 + 3x^3z^3 + 3x^3z^2/y + x^3z/y^2 + x^2yz^3 + 4x^2z^2 + 5x^2z/y + 2x^2/y^2 + xz + x + 2x/y + x/(y^2z) + y + 1/x$ & \begin{math} \begin{array}{ccccccccc|ccccc}  1&0&0&-1&1&1&-1&-2&-1&0&-1&-2&-1&0\\0&1&0&0&2&3&1&-1&-1&3&2&0&1&2\\ 0&0&1&1&-1&-1&1&2&1&0&1&2&1&0 \end{array}\end{math} & [1, 0, 2, 12, 6, 180, 560, 1680, 20230, 46200, 359352] &  dP$_7$ + 0-st  \\

37 & & P429716 & & 3 & $xy^3z^4 + xy^2z^3 + x + 2yz^2 + y + 2z + 1/y + 1/(xy) + 1/(xy^2z)$ & \begin{math} \begin{array}{ccccccc|c} 1&0&-1&1&-1&0&2&2 \\ 0&1&0&-1&1&-2&1&-1 \\ 0&0&1&0&1&2&-2&0 \end{array} \end{math} & [1, 0, 2, 6, 6, 180, 110, 2100, 11830, 19320, 334152] & $\FF_1$ + 0-st.\\  
 
\hline

38 & Q38763 & P254789 &  $2\times\frac{1}{2}(1,1,1)$, $\frac{1}{4}(1,1,3)$ & 1 & $x + y + z + 1/(xz) + 1/(xyz^2) + 1/(x^2yz) + 1/(x^2y^2z^3) + 1/(x^3y^2z^2)$ & \begin{math} \begin{array}{cccccc|c}  
1&0&-1&1&-2&2&0 \\ 0&1&1&1&2&-1&2 \end{array}\end{math} & [1, 0, 0, 6, 0, 120, 90, 0, 10080, 1680, 88200]
 &  $\FF_1$ (no~0-st.)  \\

39 &  & P519989 &  & 3 & $x + xz^2/y + y + z + 2/y^2 + 1/(xyz) + 1/(xy^3z^2)$
 & \begin{math} \begin{array}{cccccc|c} 1&-3&2&-2&-1&-1&-2 \\ 0&2&0&1&1&1&2  \end{array}\end{math} & [1, 0, 0, 6, 24, 0, 90, 2520, 2520, 1680, 151200] &  $\FF_1$ + 0-st.  \\

\hline

40 & Q38769 & P429717 &   $\frac{1}{2}(1,1,1)$, \newline $2 \times \frac{1}{3}(1,1,2)$ & 24 & $x + y + z + y^4/(xz^3) + 3y^3/(xz^2) + 3y^2/(xz) + y/x + 2y^2/(x^2z^2) + 
5y/(x^2z) + 4/x^2 + z/(x^2y) + 1/(x^3z) + 2/(x^3y) + z/(x^3y^2)$ & \begin{math}\begin{array}{cccccccc|ccccc} 1&0&1&1&-1&-1&1&2& 0&0&1&2&1 \\ 0&1&2&1&2&1&1&1& 2&3&2&2&3 \end{array} \end{math} & [1, 0, 0, 12, 0, 60, 540, 0, 10080, 33600, 18900] &  dP$_7$ + 0-st. \\

\hline

41 & Q38916 & P246021 &  $2 \times \frac{1}{3}(1,1,2)$ & 2 & $xy^3z^2 + x + y^3z + yz + y + z + y/x + 1/x + 1/(xyz)$ &  \begin{math} \begin{array}{ccccccc|c} 1&0&0&-1&1 &1&-1&0 \\ 0&1&0&0&1&0&0&0 \\ 0&0&1&2&-1&-1&3&2 \end{array} \end{math}  & [1, 0, 2, 6, 30, 120, 470, 2940, 10990, 72240, 290052] & $ \FF_1$ + 0-st. \\

42 & & P387779 & & 2 & $x + x/(y^2z) + yz^2 + y + z + y^2z^2/x + 2yz/x + y^2z/x^2$ & \begin{math}\begin{array}{cccccc|c} 1&0&2&-1&3&1&4 \\ 0&1&-1&1&-2&1&-1 \end{array}\end{math} & [1, 0, 0, 18, 0, 0, 1350, 420, 0, 141120, 100800] & $\FF_1$ (no~0-st.) \\

43 & & P411025 & & 2  & $x^2y^3z + xy^3z + xy^2z + x + 2yz + y + z + 2y/x + 1/x + z/(xy) + 
2/(x^2y) + 1/(x^3yz)$ & \begin{math} \begin{array}{ccccccc|c} 1&0&0&2&-1&2&-1&0 \\ 0&1&0&1&0&0&1&1\\ 0&0&1&0&1&-1&2&2 \end{array} \end{math} & [1, 0, 2, 0, 30, 120, 380, 2100, 5950, 60480, 189252] & $\FF_1$ (no~0-st.) \\ 

\hline


\hline

44 & Q38909 & P412951 &   $\frac{1}{2}(1,1,1)$, $\frac{1}{4}(1,1,3)$ & 3 & $x^4y^3z^4 + 2x^3y^2z^3 + 2x^2y^2z^2 + x^2yz^2 + x^2yz + 2xyz + x + x/z + y + z + 1/(xy)$ & \begin{math}\begin{array}{ccccccc|c} 1&0&0&2&-1&2&-1&0 \\ 0&1&0&1&-1&1&0&0\\ 0&0&1&0&1&-1&2&2 \end{array}\end{math} & [1, 0, 0, 6, 24, 120, 90, 1680, 9240, 39480, 163800] & $\FF_1$ (no 0-st.)\\

\hline 

45 & Q38936 & P220424 &  $2 \times \frac{1}{2}(1,1,1)$, $\frac{1}{3}(1,1,2)$ & 2 & $xy^3z^2 + x + y^2z + yz + y + z + 1/x + 1/(xz) + 1/(xy)$ & \begin{math} \begin{array}{ccccccc|c} 0&0&1&1&0&0&1&1 \\ 1&0&1&2&-1&1&-1&1 \\ 0&1&0&0&1&1&-1&1 \end{array} \end{math} & [1, 0, 2, 12, 6, 180, 560, 2100, 16870, 53760, 359352] & $\PP^1 \times \PP^1$ + 0-st.\\

46 & & P231730 & & 2 & $x^2yz + x^2z + x + x/y + y + z + z/y + y/(xz)$ & \begin{math} \begin{array}{ccccccc|c} 1&0&4&-3&3&-2&-1&1 \\ 0&0&1&0&1&0&0&1 \\ 0&1&-1&1&-2&1&1&0 \end{array} \end{math} & [1, 0, 0, 12, 24, 0, 540, 2940, 2520, 33600, 327600] & $\PP^2 \times \PP^1$ \\
\hline

47 & Q39052 & P336592 &  $\frac{1}{2}(1,1,1)$, $\frac{1}{3}(1,1,2)$ & 6 & $x^3z/y^2 + x + 2x/y + y + z + 1/z + 1/(xz)$ & \begin{math}\begin{array}{ccccccc|c} 1&0&-1&2&-2&1&3&4 \\ 0&1&0&1&-1&0&2&2 \\ 0&0&1&0&1&0&0&0 \end{array} \end{math} & [1, 0, 2, 6, 54, 60, 830, 2940, 20230, 64680, 523152] & $\FF_2$ + 0-st. \\ 

\hline

48 & & P336593 & & 2 & $x + y + z + 2y^2z^2/x + y^2z/x + 2yz/x + y/x + 1/x + 1/(xyz) + 
y^4z^3/x^2 + 2y^3z^2/x^2 + y^2z/x^2$ & \begin{math} \begin{array}{cccccccc|c} 0&0&0&-1&1&-1&0&1&0 \\ 1&0&0&0&0&1&0&0&1 \\ 0&1&0&1&-1&1&1&-1&0 \\ 0&0&1&1&0&2&-1&1&2 \end{array} \end{math} & [1, 0, 2, 0, 54, 60, 740, 1260, 18550, 45360, 422352] & $\FF_1$ + 0-st. \\

\hline

49 & Q39159 & P255193 &  $\frac{1}{3}(1,1,2)$ & 2  & $x + y + z + y/(xz) + 2y/(xz^2) + y/(xz^3) + 2/(xz) + 2/(xz^2) + 1/(xyz)$ & \begin{math}\begin{array}{ccccccc | cc}  1&0&0&-1&1&1&1&0&2 \\ 0&1&1&2&0&1&-1&2&0 \end{array} \end{math} & [1, 0, 0, 12, 48, 0, 540, 5040, 11760, 33600, 504000] & $\FF_1$ + 0-st.\\ 

\hline
 
 50 & Q39168 & P515512 &   $2 \times \frac{1}{2}(1,1,1)$ & 2 & $x + x/(yz) + x/(y^2z) + y + z + 2/(yz) + 1/(xz)$ & \begin{math}\begin{array}{cccccc|c} 1&0&1&-1&1&2&3 \\ 0&1&0&1&1&-1&0 \end{array} \end{math} & [1, 0, 0, 18, 0, 60, 1350, 0, 10080, 141120, 18900] & $\FF_2$ + 0-st.\\
 
51 & & P516400 & & 2 &  $x + y + z + 1/z + 1/(yz) + 2/(xz) + y/(x^2z)$ & \begin{math}\begin{array}{cccccc|c} 1&0&1&1&1&2&3 \\ 0&1&-1&-1&-1&-2&-3 \end{array} \end{math} & [1, 0, 2, 18, 6, 180, 1370, 1260, 25270, 148680, 283752] & $\FF_1$ (no~0-st.)\\

\hline

52 & Q39249 & P254856 & $\frac{1}{2}(1,1,1)$ & 5 & $x + y + z + 1/(xy) + 1/(xyz) + z/(x^2y) + 2/(x^2y) + 1/(x^2yz)
$ &  \begin{math}\begin{array}{cccccc|c} 0&0&1&0&0&1&1 \\ 1&1&-1&2&1&-1&1  \end{array}\end{math}  &  [1, 0, 0, 6, 48, 60, 90, 2520, 15120, 31920, 119700]
  & $\FF_1$ + 0-st.\\

53 &  & P255732 &  & 2 & $x^2y/z + x + x/z + y + z + 1/y + 1/x$
 &  \begin{math}\begin{array}{cccccc|c} 1&0&-1&1&-1&1&0 \\ 0&1&1&0&1&1&2 \end{array}\end{math}  &  [1, 0, 4, 6, 36, 240, 490, 6300, 14980, 142800, 592704]
 & $\PP^1\times\PP^1$ (no~0-st.)\\

54 & & P255860 & & 1 & $xy/z + x + x/z + x/y + yz + y + z + 1/x$
 &  \begin{math}\begin{array}{ccccccc|c} 1&0&0&1&-1&-1&2&1 \\ 0&1&0&0&0&-1&1&0 \\ 0&0&1&0&1&1&0&1  \end{array}\end{math}  & [1, 0, 2, 12, 6, 240, 560, 2940, 23590, 61320, 554652]
  & $\PP^1\times\PP^1$ + 0-st.\\

\hline
\end{longtable}
\end{center}
\end{landscape}

\pagestyle{plain}
\restoregeometry

\section{Exploring the database to diversify examples}

Constructing varieties using Laurent inversion appears to work best for varieties in higher GRDB-codimension, at least given the set of shapes we have considered so far. In this section, we explore examples in GRDB-codimensions ranging between $19$ and $5$, as well as three examples in codimension $4$ which have been constructed previously in \cite{TomStephen}, albeit in a different format.

In some cases, several 3-dimensional shapes appear, with arguably the most difficult example we obtain being discussed at length in section \ref{3DPf}, a $\PP^1$-bundle over a Pfaffian. The remaining cases are are all towers of projective bundles. 

\subsection{Interesting 3D shapes}

The shape $Z$ can be any smooth toric variety, and so far surfaces have been sufficient in the construction of all but two of the examples of Table \ref{tab1}:  variety \hyperlink{eight}{8} was scaffolded with shape $\PP^3$ and variety \hyperlink{eleven}{11} was scaffolded with shape $Z_1$, which we introduce below. We also list the $3$D polarising polytopes of the toric varieties $Z_i$ which occur as shapes for the constructions in Table \ref{tab2}.

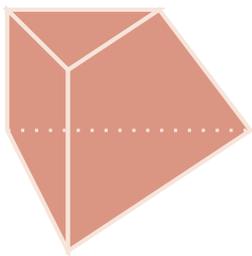
\begin{figure}[!htbp]

\begin{minipage}[t]{.45\textwidth} \vspace{0pt}%
\label{3DShapes}
\centering
\begin{tikzpicture}[scale=0.8,every node/.style={scale=0.9}]
\filldraw[color=BrickRed!10, ultra thick, fill=BrickRed!70, fill opacity =0.7, ultra thick]  (-1,1) -- (1.5,1) -- (0,0) -- cycle;
\filldraw[color=BrickRed!10, ultra thick, fill=BrickRed!70, fill opacity =0.7, ultra thick]  (-1,1) -- (0,0) -- (0,-3) -- (-1,-1) -- cycle;
\filldraw[color=BrickRed!10, ultra thick, fill=BrickRed!70, fill opacity =0.7, ultra thick]  (0,0) -- (1.5,1) -- (3,-1) -- (0,-3) -- cycle;
\draw[BrickRed!10, loosely dotted, very thick] (-1,-1) -- (3,-1);
\end{tikzpicture}
\caption*{\hypertarget{Z1}{$Z_1$} is a $\PP^1$-bundle over $\PP^2$ obtained as a toric blow up of a point on $\PP^3$. Shape for \hyperlink{eleven}{11} in Table \ref{tab1}}
\end{minipage}\hfill
\begin{minipage}[t]{.45\textwidth} \vspace{0pt}%
\centering
\begin{tikzpicture}[scale=0.5,every node/.style={scale=0.9}]
\filldraw[color=BrickRed!10, ultra thick, fill=BrickRed!70, fill opacity =0.7, ultra thick] (1.5,3.3) -- (5.5,1.3) -- (4,0.3) -- (0.5,0.3) -- (-1.5,1.3) -- cycle;
\filldraw[color=BrickRed!10, ultra thick, fill=BrickRed!70, fill opacity =0.7, ultra thick] (5.5,1.3) -- (4,0.3) -- (2.5,-3) -- (5.5,-1) -- cycle;
\filldraw[color=BrickRed!10, ultra thick, fill=BrickRed!70, fill opacity =0.7, ultra thick] (4,0.3) -- (2.5,-3) -- (-1,-3) -- (0.5,0.3) -- cycle;
\filldraw[color=BrickRed!10, ultra thick, fill=BrickRed!70, fill opacity =0.7, ultra thick] (-1,-3) -- (0.5,0.3) -- (-1.5,1.3) -- (-3,-2) -- cycle;
\draw[BrickRed!10, loosely dotted, very thick] (1.5,0.8) -- (5.5,-1);
\draw[BrickRed!10, loosely dotted, very thick] (1.5,0.8) -- (-3,-2);
\draw[BrickRed!10, loosely dotted, very thick] (1.5,0.8) -- (1.5,3.3);
\end{tikzpicture}
\caption*{\hypertarget{Z2}{$Z_2$} is a  $\PP^1$-bundle over dP$_7$. Can be obtained by blowing up a curve on $Z_3$. Shape for \hyperlink{3Dfour}{4} in Table \ref{tab2}, discussed in Section \ref{3DPf}}
\end{minipage}

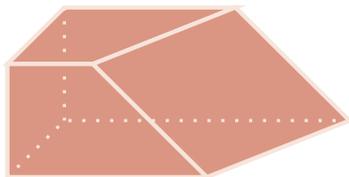
\begin{minipage}[t]{.45\textwidth} \vspace{10pt}%
\centering
\begin{tikzpicture}[scale=0.75,every node/.style={scale=0.9}]
\filldraw[color=BrickRed!10, ultra thick, fill=BrickRed!70, fill opacity =0.7, ultra thick]  (-1,0) -- (0,1) -- (3,1) -- (0.5,0) -- cycle;
\filldraw[color=BrickRed!10, ultra thick, fill=BrickRed!70, fill opacity =0.7, ultra thick]  (-1,0) -- (0.5,0) -- (2.5,-2) -- (-1,-2) -- cycle;
\filldraw[color=BrickRed!10, ultra thick, fill=BrickRed!70, fill opacity =0.7, ultra thick]  (0.5,0) -- (3,1) -- (5,-1) -- (2.5,-2) -- cycle;
\draw[BrickRed!10, loosely dotted, very thick] (0,1) -- (0,-1);
\draw[BrickRed!10, loosely dotted, very thick] (-1,-2) -- (0,-1);
\draw[BrickRed!10, loosely dotted, very thick] (5,-1) -- (0,-1);
\end{tikzpicture}
\vspace{0.5cm}
\caption*{\hypertarget{Z3}{$Z_3$} is a $\PP^1$-bundle over $\FF_1$ and a $\PP^1\times\PP^1$-bundle on $\PP^1$. Shape for \hyperlink{3Dfive}{5}, \hyperlink{3Dsix}{6} and \hyperlink{3Dseventeen}{16} in Table \ref{tab2}}
\end{minipage}\hfill
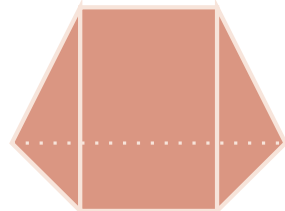
\begin{minipage}[t]{.45\textwidth} \vspace{10pt}%
\centering
\begin{tikzpicture}[scale=0.9,every node/.style={scale=0.9}]
\filldraw[color=BrickRed!10, ultra thick, fill=BrickRed!70, fill opacity =0.7, ultra thick] (-1,1) -- (-2,-1) -- (-1,-2) -- cycle;
\filldraw[color=BrickRed!10, ultra thick, fill=BrickRed!70, fill opacity =0.7, ultra thick] (-1,1) -- (1,1) -- (1,-2) -- (-1,-2) -- cycle;
\filldraw[color=BrickRed!10, ultra thick, fill=BrickRed!70, fill opacity =0.7, ultra thick] (1,1) -- (1,-2) -- (2,-1) -- cycle;

\draw[BrickRed!10, loosely dotted, very thick] (-2,-1) -- (2,-1);
\end{tikzpicture}
\caption*{\hypertarget{Z4}{$Z_4$} is a $\PP^2$-bundle over $\PP_1$ obtained by blowing up a toric curve on $\PP^3$. Shape for \hyperlink{3Deleven}{11} in Table \ref{tab2}.}
\end{minipage}
\end{figure}

Other 3D shapes that occur come from products: $\PP^2\times \PP^1$ and $\FF_1 \times \PP^1$. We also encounter many of the $2$D shapes in Table \ref{tab1}.  

\subsection{The 3D Pfaffian: Q33018}
\label{3DPf}
Consider the rigid MMLP \begin{align*}
    f &=xy^2z^3 + xy^2z^2 + xyz^2 + 2xyz + x + 3yz^2 + 2yz + y + 4z \\ &+ \frac{1}{y} + 3\frac{z}{x} + \frac{3}{x} + \frac{1}{xz} + \frac{4}{xy} + \frac{1}{xy^2z} + \frac{1}{x^2y} + \frac{3}{x^2yz} + \frac{1}{x^2y^2z} + \frac{1}{x^3y^2z^2} \end{align*} supported on the polytope P$31470$. This is mutation equivalent (via $4$ one-step mutations) to the Laurent polynomial:
\begin{align*}
g = \frac{x}{yz} + \frac{x}{yz^2} + yz + 2z + \frac{z}{y} + \frac{1}{y} + \frac{1}{yz} + \frac{2}{yz^2}+ \frac{1}{yz^3} + \frac{yz^2}{x} \\ + \frac{2yz}{x} + \frac{y}{x} + \frac{z^2}{x} + \frac{4z}{x} + \frac{4}{x} + \frac{1}{xz} + \frac{z}{xy} + \frac{2}{xy} + \frac{1}{xyz} 
\end{align*}
with Newton polytope:

\begin{figure}[ht]
\centering
\begin{tikzpicture}[every node/.style={scale=0.8}]

\draw[gray!80, thick] (-6,-1.5-0.2) -- (-3,-2);
\draw[gray!80, thick] (-3,-2) -- (-1,-3);
\draw[gray!80, thick] (-3,-2) -- (0,0);
\draw[gray!80, thick] (0,0) -- (0,2.5-0.2); 
\draw[gray!80, thick] (0,0) -- (3.5,0);
\draw[gray!80, thick] (3.5,0) -- (5.5,-1);
\draw[gray!80, thick] (3.5,0) -- (0,2.5-0.2);

\filldraw[gray] (-1.5,1.5-0.2) circle (1.5pt) node[anchor=south] {};

\draw[black, very thick] (-3,3-0.2) -- (-4.5,2-0.2);
\draw[black, very thick] (-3,3-0.2) -- (0,2.5-0.2);
\draw[black, very thick] (-3,3-0.2) -- (4,0.5-0.2);
\draw[black, very thick] (-4.5,2-0.2) -- (4,0.5-0.2);
\draw[black, very thick] (-4.5,2-0.2) -- (-6,-1.5-0.2);
\draw[black, very thick] (-6,-1.5-0.2) -- (-1,-3);
\draw[black, very thick] (-6,-1.5-0.2) -- (2.5,-3);
\draw[black, very thick] (4,0.5-0.2) -- (0,2.5-0.2);
\draw[black, very thick] (4,0.5-0.2) -- (5.5,-1);
\draw[black, very thick] (2.5,-3) -- (4,0.5-0.2);
\draw[black, very thick] (2.5,-3) -- (5.5,-1);
\draw[black, very thick] (2.5,-3) -- (-1,-3);


\filldraw[black] (-6,-1.5-0.2) circle (1.5pt) node[anchor=north east] {$(0,-1,-3)$};
\filldraw[gray] (-3,-2) circle (1.5pt) node[anchor=west] {\textcolor{black}{$(-1,-1,-1)$}};
\filldraw[black] (-1,-3) circle (1.5pt) node[anchor=north] {$(-1,0,-1)$};
\filldraw[gray] (0,0) circle (1.5pt) node[anchor=north west] {\textcolor{black}{$(-1,-1,1)$}};
\filldraw[black] (0,2.5-0.2) circle (1.5pt) node[anchor=south] {$(0,-1,1)$};
\filldraw[gray] (3.5,0) circle (1.5pt) node[anchor=north east] {\textcolor{black}{$(-1,0,2)$}};
\filldraw[black] (5.5,-1) circle (1.5pt) node[anchor=west] {$(-1,1,2)$};
\filldraw[black] (-3,3-0.2) circle (1.5pt) node[anchor=south] {$(1,-1,-1)$};
\filldraw[black] (-4.5,2-0.2) circle (1.5pt) node[anchor=east] {$(1,-1,-2)$};
\filldraw[black] (4,0.5-0.2) circle (1.5pt) node[anchor=south west] {$(0,1,1)$};
\filldraw[black] (2.5,-3) circle (1.5pt) node[anchor=north] {$(-1,1,0)$};

\filldraw[gray] (-3/2,-1) circle (1.5pt) node[anchor=south] {};
\filldraw[gray] (0.5,-2) circle (1.5pt) node[anchor=south] {};
\filldraw[gray] (2,-1) circle (1.5pt) node[anchor=south] {};
\filldraw[black] (4,-2) circle (1.5pt) node[anchor=south] {};
\filldraw[gray] (-3,0.5-0.2) circle (1.5pt) node[anchor=south] {};
\filldraw[gray] (-4.5,-0.5-0.2) circle (1.5pt) node[anchor=south] {};
\filldraw[black] (2,1.3) circle (1.5pt) node[anchor=south] {};
\filldraw[black] (-1,-0.7) circle (1.5pt) node[anchor=south] {};

\draw[gray, thick] (0.55,0.35) -- (0.5-0.05,0.3-0.05);
\draw[gray, thick] (0.55,0.3-0.05) -- (0.5-0.05,0.3+0.05);

\end{tikzpicture}
\caption{Newton polytope of $g$}
\label{P31470}
\end{figure}
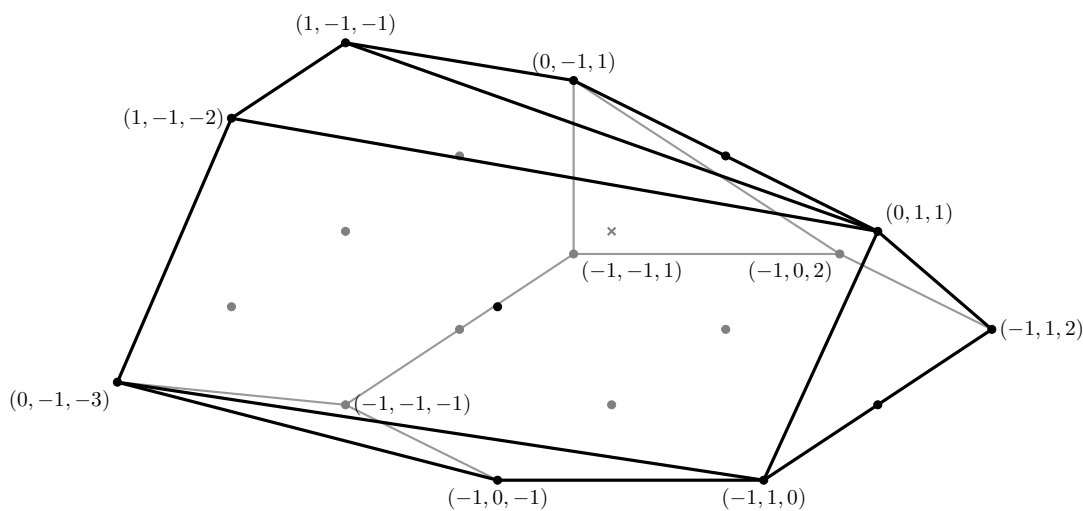

\subsubsection{The weight matrix and the Pfaffian format}
\label{Pfaffians}
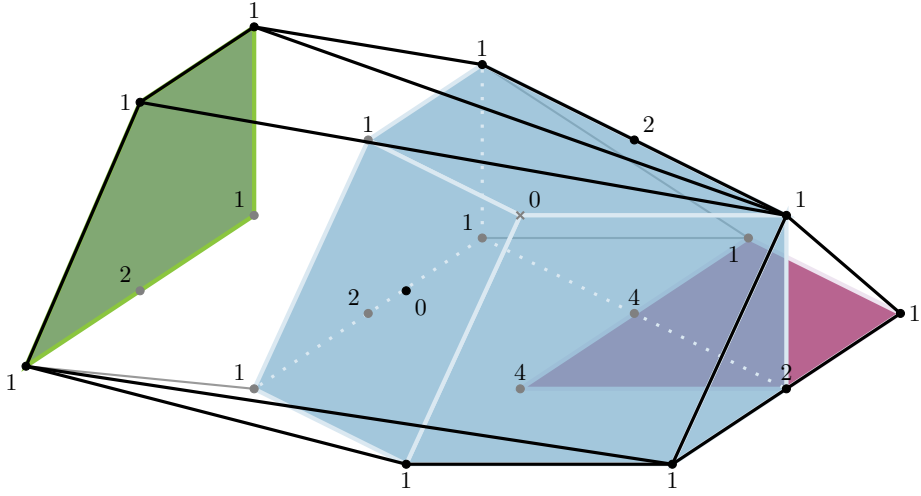
\begin{figure}[ht]
\centering
\begin{tikzpicture}[every node/.style={scale=0.9}]

\draw[gray!80, thick] (-6,-1.5-0.2) -- (-3,-2);
\draw[gray!80, thick] (-3,-2) -- (-1,-3);
\draw[gray!80, thick] (0,0) -- (3.5,0);
\draw[gray!80, thick] (3.5,0) -- (5.5,-1);
\draw[gray!80, thick] (3.5,0) -- (0,2.5-0.2);

\filldraw[color=Plum!10, ultra thick, fill=RedViolet!70, ultra thick] (5.5,-1) -- (4,-2) -- (0.5,-2) -- (3.5,0) -- cycle;

\filldraw[color=LimeGreen!100, ultra thick, fill=OliveGreen!60, ultra thick] (-3,3-0.2) --  (-3,0.3) -- (-6,-1.7) -- (-4.5,2-0.2) -- cycle;


\filldraw[color=MidnightBlue!10, ultra thick, fill=MidnightBlue!40, fill opacity =0.7, ultra thick] (0,2.3) -- (4,0.3) -- (0.5,0.3) -- (-1.5,1.3) -- cycle;
\filldraw[color=MidnightBlue!10, ultra thick, fill=MidnightBlue!40, fill opacity =0.7, ultra thick] (4,0.3) -- (2.5,-3) -- (4,-2) -- cycle;
\filldraw[color=MidnightBlue!10, ultra thick, fill=MidnightBlue!40, fill opacity =0.7, ultra thick] (4,0.3) -- (2.5,-3) -- (-1,-3) -- (0.5,0.3) -- cycle;
\filldraw[color=MidnightBlue!10, ultra thick, fill=MidnightBlue!40, fill opacity =0.7, ultra thick] (-1,-3) -- (0.5,0.3) -- (-1.5,1.3) -- (-3,-2) -- cycle;

\draw[MidnightBlue!10, loosely dotted, very thick] (0,0) -- (4,-2);
\draw[MidnightBlue!10, loosely dotted, very thick] (0,0) -- (-3,-2);
\draw[MidnightBlue!10, very thick, loosely dotted] (0,0) -- (0,2.3);

\filldraw[gray] (-1.5,1.5-0.2) circle (1.5pt) node[anchor=south] {\textcolor{black}{$1$}};

\draw[black, very thick] (-3,3-0.2) -- (-4.5,2-0.2);
\draw[black, very thick] (-3,3-0.2) -- (0,2.5-0.2);
\draw[black, very thick] (-3,3-0.2) -- (4,0.5-0.2);
\draw[black, very thick] (-4.5,2-0.2) -- (4,0.5-0.2);
\draw[black, very thick] (-4.5,2-0.2) -- (-6,-1.5-0.2);
\draw[black, very thick] (-6,-1.5-0.2) -- (-1,-3);
\draw[black, very thick] (-6,-1.5-0.2) -- (2.5,-3);
\draw[black, very thick] (4,0.5-0.2) -- (0,2.5-0.2);
\draw[black, very thick] (4,0.5-0.2) -- (5.5,-1);
\draw[black, very thick] (2.5,-3) -- (4,0.5-0.2);
\draw[black, very thick] (2.5,-3) -- (5.5,-1);
\draw[black, very thick] (2.5,-3) -- (-1,-3);

\filldraw[black] (-6,-1.5-0.2) circle (1.5pt) node[anchor=north east] {$1$};
\filldraw[gray] (-3,-2) circle (1.5pt) node[anchor=south east] {\textcolor{black}{$1$}};
\filldraw[black] (-1,-3) circle (1.5pt) node[anchor=north] {$1$};
\filldraw[gray] (0,0) circle (1.5pt) node[anchor=south east] {\textcolor{black}{$1$}};
\filldraw[black] (0,2.5-0.2) circle (1.5pt) node[anchor=south] {$1$};
\filldraw[gray] (3.5,0) circle (1.5pt) node[anchor=north east] {\textcolor{black}{$1$}};
\filldraw[black] (5.5,-1) circle (1.5pt) node[anchor=west] {$1$};
\filldraw[black] (-3,3-0.2) circle (1.5pt) node[anchor=south] {$1$};
\filldraw[black] (-4.5,2-0.2) circle (1.5pt) node[anchor=east] {$1$};
\filldraw[black] (4,0.5-0.2) circle (1.5pt) node[anchor=south west] {$1$};
\filldraw[black] (2.5,-3) circle (1.5pt) node[anchor=north] {$1$};

\filldraw[gray] (-3/2,-1) circle (1.5pt) node[anchor=south east] {\textcolor{black}{$2$}};
\filldraw[gray] (0.5,-2) circle (1.5pt) node[anchor=south] {\textcolor{black}{$4$}};
\filldraw[gray] (2,-1) circle (1.5pt) node[anchor=south] {\textcolor{black}{$4$}};
\filldraw[black] (4,-2) circle (1.5pt) node[anchor=south] {$2$};
\filldraw[gray] (-3,0.5-0.2) circle (1.5pt) node[anchor=south east] {\textcolor{black}{$1$}};
\filldraw[gray] (-4.5,-0.5-0.2) circle (1.5pt) node[anchor=south east] {\textcolor{black}{$2$}};
\filldraw[black] (0.5,0.5-0.2) node[anchor=south west] {$0$};
\filldraw[black] (2,1.3) circle (1.5pt) node[anchor=south west] {$2$}; 
\filldraw[black] (-1,-0.7) circle (1.5pt) node[anchor=north west] {$0$};

\draw[gray, thick] (0.55,0.35) -- (0.5-0.05,0.3-0.05);
\draw[gray, thick] (0.55,0.3-0.05) -- (0.5-0.05,0.3+0.05);

\end{tikzpicture}
\caption{a scaffolding of $f_4$}
\label{LP31470}
\end{figure}

Figure \ref{LP31470} contains a scaffolding compatible with the polynomial $g$: two $2$-dimensional quadrilateral struts and one full-dimensional strut. Note that we cannot choose the $3$-dimensional strut itself to be the polarising polytope of the shape $Z$, as it has a non-$\QQ$-factorial singularity (coming from the $4$-valent vertex of the strut). To obtain the correct object we make a small resolution and we recover $Z_2$ below, described by its ample polytope. Indeed, this shape is a $\PP^{1}$-bundle over $\textup{dP}_7$, whose torus-invariant divisors are labeled as in the picture:

\begin{figure}[ht]
\label{SP31470}
\centering
\begin{tikzpicture}[scale=0.7,every node/.style={scale=0.9}]

\draw[{Circle[MidnightBlue]}-{Stealth[BrickRed,scale=2]}] (3.5,1.05)--(5.4,2.6) node[right]{$\Delta_3$};
\draw[{Circle[MidnightBlue]}-{Stealth[BrickRed,scale=2]}] (0,1.15)-- (-2,2.35) node[left]{$\Delta_4$};
\draw[{Circle[MidnightBlue,]}-{Stealth[BrickRed,scale=2]}] (2,-1.3)--(2,-4) node[right]{$\Delta_2$};

\filldraw[color=BrickRed!10, ultra thick, fill=BrickRed!70, fill opacity =0.7, ultra thick] (1.5,3.3) -- (5.5,1.3) -- (4,0.3) -- (0.5,0.3) -- (-1.5,1.3) -- cycle;
\filldraw[color=BrickRed!10, ultra thick, fill=BrickRed!70, fill opacity =0.7, ultra thick] (5.5,1.3) -- (4,0.3) -- (2.5,-3) -- (5.5,-1) -- cycle;
\filldraw[color=BrickRed!10, ultra thick, fill=BrickRed!70, fill opacity =0.7, ultra thick] (4,0.3) -- (2.5,-3) -- (-1,-3) -- (0.5,0.3) -- cycle;
\filldraw[color=BrickRed!10, ultra thick, fill=BrickRed!70, fill opacity =0.7, ultra thick] (-1,-3) -- (0.5,0.3) -- (-1.5,1.3) -- (-3,-2) -- cycle;

\draw[BrickRed!10, loosely dotted, very thick] (1.5,0.8) -- (5.5,-1);
\draw[BrickRed!10, loosely dotted, very thick] (1.5,0.8) -- (-3,-2);
\draw[BrickRed!10, loosely dotted, very thick] (1.5,0.8) -- (1.5,3.3);

\draw[{Circle[MidnightBlue]}-{Stealth[BrickRed,scale=2]}] (2,1.3)--(2,3.8) node[right]{$\Delta_1$};
\draw[{Circle[MidnightBlue]}-{Stealth[BrickRed,scale=2]}] (3/2,-2.5/2)-- (3/2-0.8,-2.5/2-2.4)  node[left]{$\Delta_6$};
\draw[{Circle[MidnightBlue]}-{Stealth[BrickRed,scale=2]}] (-1.2,-0.8)--  (-1.2-2.1,-0.8+0.8)  node[left]{$\Delta_5$};
\draw[{Circle[MidnightBlue]}-{Stealth[BrickRed,scale=2]}] (9.5/2,-0.7/2)-- (9.5/2+2,-0.7/2-1)  node[right]{$\Delta_7$};
\end{tikzpicture}
\end{figure}
The variety we obtain is of codimension $4$, albeit in a special configuration because of the $\PP^{1}$-bundle structure, which corresponds to the vertical directions $\Delta_{1}$ and $\Delta_{2}$ on the shape. The weight matrix and line bundles are:
\[\begin{blockarray}{*{10}{c} | c | ccccc}
 \BAmulticolumn{3}{c}{} & \Delta_1 & \Delta_2 & \Delta_3 & \Delta_4 & \Delta_5 & \Delta_6  &\Delta_7& H & L_{1} & L_{2} & L_{3} & L_{4} &L_{5}\\
\BAhline
\begin{block}{*{10}{c} | c |   ccccc}
  1 & 0& 0 & 1 & 0 &  -1 & 1  &3  & 2 & -1  &  1  & 2   & 0   & 3 &  1 &  2\\
  0 & 1& 0 & 0 & 1 &   1 & 1  &0 &  0 &  1  &  1  &  1  & 2 	&1   &1&	1\\
  0 & 0& 1 & -1 & 1 &  2 & 0 &-1 & -1 & 1   &  0  &   0&  1 & -1 &	1& 1	\\
\end{block}
\BAhline
 x_1 & x_2 & x_3 & x_4 & x_5 & x_6 & y & x_7 & x_8 & x_9& &&&&& 
\end{blockarray}
\]
where $H$ corresponds to a hypersurface on $Y$ which intersects a Pfaffian variety given by the $L_i$. The chosen ample class is $\omega=(1,1,1)=-K_X=-K_Y+L$, as discussed in \ref{PfK} (here $U=\emptyset$). Despite the fact that the chosen $Z_2$ is not a tower of projective bundles, we obtain $X_{P}$ as a toric Pfaffian complete intersection (TPCI). More precisely, in these coordinates $X_P$ is given by the following six equations:
\begin{equation}
\label{superfluous}
x_4x_5= x_1x_2\textup{ and } \\
\end{equation}
\begin{equation}
\label{PF}
x_7x_9=x_2x_3x_8, \ \ yx_9=x_2^2x_3, \ \  yx_8=x_2x_7, \ \ x_6x_8=x_1^2x_9,  \textup{ and }  x_6x_7=x_1^2x_2x_3. \end{equation}
 Here, equations \eqref{PF} correspond to sections of $L_1, \ L_2, \ L_3,\ L_4 \textup{ and } L_5$, respectively, and fit into the format:
\[
\begin{pmatrix}
 0 & x_1^2 & x_{6} & y & x_{2} \\
& 0 & 0 & x_{7} & x_{8}\\
 &\textup{-sym} & 0 & x_{2}x_{3} & x_{9}\\
 &&& 0& 0\\
 &&&& 0\\
\end{pmatrix}, \] i.e. they are the $4\times4$ Pfaffians of this $5\times 5$ antisymmetric matrix. The format is a tool for deforming this toric variety, however before we proceed we simplify the situation: since the first deformed equation is linear in $y$, we can solve for this coordinate, eliminate equation \eqref{superfluous} and obtain the weight matrix and line bundles appearing in Table \ref{tab1}:
\[\begin{blockarray}{*{9}{c} | ccccc}
 x_{1} & x_2 & x_3 & x_4 & x_5 & x_6  &x_7 &x_{8}& x_{9} & L_{1} & L_{2} & L_{3} & L_{4} &L_{5}\\
\BAhline
\begin{block}{*{9}{c} |   ccccc}
  1 & 0& 0 & 1 & 0 &  -1 &3  & 2 & -1  & 2  & 0   & 3 &  1 &  2\\
  0 & 1& 0 & 0 & 1 &   1 &0 &  0 &  1   &  1  & 2 	&1   &1&	1\\
  0 & 0& 1 & -1 & 1 &  2 &-1 & -1 & 1  & 0&  1 & -1 &	1& 1\\
\end{block}
\end{blockarray}
\] where this is now a bona fide toric Pfaffian variety. In what follows, we slightly abuse notation by still referring to the ambient space given by the new matrix as $Y$ and to the toric variety inside it as $X_P$.

$X_{P}$ is given by 5 binomial equations, obtained by substituting $y$ in \eqref{PF}. To aid in determining the deformation, we record the weights above each monomial and obtain the following for $X_P$:
\[ A_P= 
\begin{pmatrix}
 &\fbox{200} & \fbox{-112} & \fbox{110} & \fbox{010} \\
 0 & x_1^2 & x_{6} & x_{4}x_{5} & x_{2} \\
&  & \fbox{101} & \fbox{30-1}& \fbox{20-1} \\ 
& 0 & 0 & x_{7} & x_{8 }\\
 &&  & \fbox{011}&\fbox{-111} \\
 &\textup{-sym} & 0 & x_{2}x_{3} & x_{9}\\
 &&&& \fbox{11-1}\\
 &&& 0& 0 \\
  &&&& \\
 &&&& 0\\
\end{pmatrix} \]

The weights of the matrix entries remain fixed when moving from $X_P$ to $X$. We deform these equations by adding all monomials of each homogeneity type in the splitting (maintaining the coefficient convention throughout). In other words, deforming a Pfaffian variety is done not by adding all monomials in $H^0(Y,L_i)$ until we obtain general sections, as is the case for complete intersections, but only by adding monomials which respect the splitting as below:
\[ A= \left( 
\begin{tabular}{m{0.5cm}m{2.2cm}m{2cm}m{3.5cm}m{2.5cm}}
  &\hfil \fbox{200} &\hfil \fbox{-112} &\hfil \fbox{110} & \hfil \fbox{010} \\
 \hfil 0 & $x_1^2+x_3x_8+x_1x_3x_4+x_3^2x_4^2$ & \hfil $x_6+x_3x_9$ & $x_1x_2+x_4x_5+x_8x_9+x_4^2x_6x_1x_4x_9+x_2x_3x_4+x_3x_4^2x_9$& \hfil $   x_2+x_4x_9$ \\
&  &\hfil \fbox{101} &\hfil \fbox{30-1}& \hfil\fbox{20-1} \\ 
&\hfil 0 & \hfil$x_1x_3+x_3^2x_4$ & $x_7 + x_1x_8 + x_1^2x_4 + x_3^2x_4^3+x_3x_4x_8+x_1x_3x_4^2$ &\hfil $x_8+x_1x_4+x_3x_4^2$\\
 &&  &\hfil \fbox{011}&\hfil\fbox{-111} \\
 &\hfil \textup{-sym} &\hfil 0 & $x_2x_3+x_5+x_4x_6+x_1x_9+x_3x_4x_9$ & \hfil$x_9$\\
 &&&& \hfil\fbox{11-1}\\
 &&& \hfil0& \hfil$x_2x_4+x_4^2x_9$\\
  &&&& \\
 &&&& \hfil0\\
\end{tabular} \right) \]

The equations of $X$ are the Pfaffians of this matrix, which we omit writing down as all the required information can be determined directly from the more compact form of $A$. 

\subsubsection{Singularities} GRDB suggests that we should find three singular points: $2\times \frac{1}{2}(1,1,1)$ and $1\times \frac{1}{3}(1,1,2)$. 

The ambient space $Y$ is not an orbifold, its non-$\QQ$-factorial locus consisting of both curves and points. Given the stability condition is $\omega=(1,1,1)$ following the discussion in Section \ref{PfK}, $Y$ is covered by three non-orbifold charts: $U_{15}$, $U_{19}$, $U_{68}$ and  $17$ orbifold charts:  $U_{123},  U_{126}, U_{146}, U_{234}, U_{237}, U_{238}, \ U_{267}, U_{345}, U_{346}, U_{349}, U_{357}, U_{358}$, $U_{379}, \ U_{389}, U_{467}, U_{567},$ and $U_{679}$. We recall the notation $U_{i_1\ldots i_k}:=\{x_{i_1}=\ldots=x_{i_k}=1\}$.

\renewcommand\labelitemi{{\boldmath$\cdot$}}

We first discuss the six singular curves: \vspace{-0.2cm}
\begin{itemize}
\item $C_1=\{ x_2=x_3=x_4=x_6=x_7=x_8=0\}$. Note that this is a non-$\QQ$-factorial curve, belonging to the charts $U_{15}$ and $U_{19}$ (notice in particular that $x_1$ never vanishes here). We restrict the equations of $X$ to $C_1$ and they become the Pfaffians of 
\[ A|_{C_1}= \left( 
\begin{tabular}{m{0.5cm}m{1cm}m{1cm}m{1.5cm}m{1cm}}
  &\hfil \fbox{200} &\hfil \fbox{-112} &\hfil \fbox{110} & \hfil \fbox{010} \\
 \hfil 0 & \hfil $x_1^2$ & \hfil $0$ & \hfil $0$& \hfil $ 0$ \\
&  &\hfil \fbox{101} &\hfil \fbox{30-1}& \hfil\fbox{20-1} \\ 
&\hfil 0 & \hfil$0$ & \hfil $0 $ &\hfil $0$\\
 &&  &\hfil \fbox{011}&\hfil\fbox{-111} \\
 &\hfil \textup{-sym} &\hfil 0 &\hfil $x_5+x_1x_9$ & \hfil$x_9$\\
 &&&& \hfil\fbox{11-1}\\
 &&& \hfil0& \hfil$0$\\
  &&&& \\
 &&&& \hfil0\\
\end{tabular} \right), \]
i.e. all but two equations disappear and we obtain the system \[ \left\{ \begin{array}{l} x_9x_1^2=0 \\
x_1^2(x_5+x_1x_9)=0 \end{array} \right. \] This system has no solutions, as $x_5=x_9=0$ does not belong to either chart containing the curve $C_1$. 
\item $C_2=\{ x_1=x_2=x_4=x_5=x_7=0\}$. This belongs to $U_{68}$ and $U_{389}$, therefore $x_8\neq 0$. The matrix $A$ becomes:
\[ A|_{C_2}= \left( 
\begin{tabular}{m{0.5cm}m{1cm}m{1.5cm}m{1cm}m{1cm}}
  &\hfil \fbox{200} &\hfil \fbox{-112} &\hfil \fbox{110} & \hfil \fbox{010} \\
 \hfil 0 & \hfil $x_3x_8$ & \hfil $x_6+x_3x_9$ & \hfil $x_8x_9$& \hfil $ 0$ \\
&  &\hfil \fbox{101} &\hfil \fbox{30-1}& \hfil\fbox{20-1} \\ 
&\hfil 0 & \hfil$0$ & \hfil $0 $ &\hfil $x_8$\\
 &&  &\hfil \fbox{011}&\hfil\fbox{-111} \\
 &\hfil \textup{-sym} &\hfil 0 &\hfil $0$ & \hfil$x_9$\\
 &&&& \hfil\fbox{11-1}\\
 &&& \hfil0& \hfil$0$\\
  &&&& \\
 &&&& \hfil0\\
\end{tabular} \right) \] and we have the three surviving equations 
\[ \left\{ \begin{array}{l} x_8(x_6+x_3x_9)-x_3x_8x_9=0 \\
x_8^2x_9=0\\
x_8x_9^2=0 \end{array} \right. \] This time the contradiction is $x_6=x_9=0$.

\item $C_3=\{ x_1=x_2=x_4=x_5=x_8=0\}$ belongs to $U_{379}$ and $U_{679}$, therefore both $x_7$ and $x_9$ do not vanish here. The matrix is  

\[ A|_{C_3}= \left( 
\begin{tabular}{m{0.5cm}m{1cm}m{1.5cm}m{1cm}m{1cm}}
  &\hfil \fbox{200} &\hfil \fbox{-112} &\hfil \fbox{110} & \hfil \fbox{010} \\
 \hfil 0 & \hfil $0$ & \hfil $x_6+x_3x_9$ & \hfil $0$& \hfil $ 0$ \\
&  &\hfil \fbox{101} &\hfil \fbox{30-1}& \hfil\fbox{20-1} \\ 
&\hfil 0 & \hfil$0$ & \hfil $x_7 $ &\hfil $0$\\
 &&  &\hfil \fbox{011}&\hfil\fbox{-111} \\
 &\hfil \textup{-sym} &\hfil 0 &\hfil $0$ & \hfil$x_9$\\
 &&&& \hfil\fbox{11-1}\\
 &&& \hfil0& \hfil$0$\\
  &&&& \\
 &&&& \hfil0\\
\end{tabular} \right) \] and among the equations notice $x_7x_9=0$, which is impossible in the above charts.

\item $C_4=\{ x_1=x_4=x_6=x_7=x_9=0\}$ belongs to $U_{238}$ and $U_{358}$ and
\[ A|_{C_4}= \left( 
\begin{tabular}{m{0.5cm}m{1cm}m{1cm}m{2cm}m{1cm}}
  &\hfil \fbox{200} &\hfil \fbox{-112} &\hfil \fbox{110} & \hfil \fbox{010} \\
 \hfil 0 & \hfil $x_3x_8$ & \hfil $0$ & \hfil $0$& \hfil $ x_2$ \\
&  &\hfil \fbox{101} &\hfil \fbox{30-1}& \hfil\fbox{20-1} \\ 
&\hfil 0 & \hfil$0$ & \hfil $0 $ &\hfil $x_8$\\
 &&  &\hfil \fbox{011}&\hfil\fbox{-111} \\
 &\hfil \textup{-sym} &\hfil 0 &\hfil $x_5+x_2x_3$ & \hfil$0$\\
 &&&& \hfil\fbox{11-1}\\
 &&& \hfil0& \hfil$0$\\
  &&&& \\
 &&&& \hfil0\\
\end{tabular} \right). \] There are three remaining equations:
 \[ \left\{ \begin{array}{l} x_8(x_5+x_2x_3)=0 \\
x_2(x_5+x_2x_3)=0\\
x_3x_8(x_5+x_2x_3)=0
\end{array} \right. \] whose only solution is $x_5+x_3x_2=0$\footnote{This binomial appears in all three equations, with arbitrary but fixed coefficients --- sometimes expressions look the same because of the coefficient convention, but may lead to a contradiction. That is not the case in this example.}. This is one point belonging to both charts but which is not the origin of either. Since $U_{238}=\frac{1}{2}(1,1,0,1,1,1)_{145679}$, and $C_4$ is given by the $x_5$-axis which $X$ intersects transversally, the point is of type $\frac{1}{2}(1,1,1)$.
\item $C_5=\{ x_1=x_4=x_6=x_8=x_9=0\}$ is almost identical to $C_4$, except that we find a single $\frac{1}{3}(1,1,2)$ point. 
\item $C_6=\{ x_1=x_2=x_3=x_8=x_9=0\}$ is similar to $C_3$, and $X$ does not intersect it.
\end{itemize}

The remaining singular points of $Y$ which are not contained in $C_1, \ldots, C_6$ are $P_1=\{x_1=x_3=x_4=x_5=x_8=x_9=0\}$, $P_2=\{x_1=x_2=x_3=x_4=x_5=x_7=x_9=0\}$ and $P_3=\{x_3=x_4=x_5=x_7=x_8=x_9=0\}$, i.e. the origins of $U_{267}$, $U_{68}$ and $U_{126}$ respectively. The latter is the only point belonging to $X$, as all $4\times4$ Pfaffians of the matrix 
\[ A|_{P_2}= \left( 
\begin{tabular}{m{0.5cm}m{1cm}m{1cm}m{1cm}m{1cm}}
 &\hfil \fbox{200} &\hfil \fbox{-112} &\hfil \fbox{110} & \hfil \fbox{010} \\
\hfil $0 $& \hfil $x_1^2$ & \hfil $x_6$ & \hfil $x_1x_2$& \hfil $ x_2$ \\
&  &\hfil \fbox{101} &\hfil \fbox{30-1}& \hfil\fbox{20-1} \\ 
&\hfil 0 & \hfil$0$ & \hfil $0 $ &\hfil $0$\\
 &&  &\hfil \fbox{011}&\hfil\fbox{-111} \\
 &\hfil \textup{-sym} &\hfil 0 &\hfil $0$ & \hfil$0$\\
 &&&& \hfil\fbox{11-1}\\
 &&& \hfil0& \hfil$0$\\
  &&&& \\
 &&&& \hfil0\\
\end{tabular} \right)  \] vanish. It is the second expected point of type $\frac{1}{2}(1,1,1)$.

\subsubsection{Anticanonical degree}
\label{degree}
Since the shape $Z$ is neither a tower of projective bundles \cite[Theorem 2.21]{doran_harder_2016} \cite[\S 8]{LaurentInversion} nor a product of projective spaces \cite[Proposition 5.1]{Prince_2019}, existing theoretical results\footnote{Note that there are already some examples of successful Pfaffian and TCPI constructions. However, the shapes involved are simpler: either 2D \cite{LaurentInversion} \cite{1/3} or a product \cite{prince2019cracked}.} do not determine the equations of $X_P$ in $Y$. As a cross-check that we have constructed the correct variety $X$ (i.e.~a deformation of $X_P$), we therefore compute the anticanonical degree of $X$ and verify that it coincides with the prediction from GRDB, that is, $47/3$.

We compute this degree inside an orbifold $\tilde{Y}$ quasi-isomorphic to $Y$, which is an obtained by slightly changing the stability condition from $\omega=(1,1,1)$ to (for example) $\tilde{\omega}=(1+\varepsilon,1+\frac{\varepsilon}{2},1+\frac{\varepsilon}{3})$ where $\varepsilon$ is small, which is now contained in a maximal-dimensional chamber of the secondary fan. Since $X$ avoids the non-orbifold locus of $Y$, it is isomorphic to $\tilde{X}\subset \tilde{Y}$ given by the same equations and in particular it has the same degree. 

Write the Chow ring of $\tilde{Y}$ in terms of its generators $M =(1,0,0), \ N=(0,1,0) \textup{ and } P=(0,0,1)$. The relations are given by writing the components of the irrelevant ideal in terms of $M, N$ and $P$. In our case \begin{align*}I_{\tilde{\omega}} = &(x_1, x_3, x_6)\cdot(x_2, x_4, x_5, x_8, x_9)\cdot(x_2, x_5, x_6, x_9)\cdot(x_3, x_5, x_6, x_9)\cdot \\ &
(x_1, x_3, x_7)\cdot(x_1, x_4, x_7, x_8)\cdot(x_2, x_4, x_7, x_8),\end{align*}
resulting in the relations: 
\begin{align*}
&M P (-M+N+2P)=0\\
&N(M-P ) (N+P) (2M-P)(-M+N+P)=0\\
&N (N+P) (-M+N+2P)(-M+N+P)=0\\
&P (N+P) (-M+N+2P)(-M+N+P)=0\\
&MP(3M-P)=0\\
&M(M-P )(3M-P)(2M-P)=0\\ 
&N(M-P )(3M-P)(2M-P)=0.
\end{align*}
A key ingredient in finding $(-K_{\widetilde{X}})^3$ is computing the fundamental class of the determinantal variety $\widetilde{X}$ by using the squaring principle \cite[Proposition 9]{HARRIS198471} and the corresponding Porteous formula \cite[Theorem 10]{HARRIS198471}, for which the set-up is as follows: 

\begin{thm}[Porteous formula, \cite{HARRIS198471}]
Let $E$ be a rank $5$ homogeneous vector bundle and $L$ a line bundle on the simplicial toric variety $\widetilde{Y}$. If $s: E \otimes L \to E^\vee $ is a general skew-symmetric map, and $\widetilde{X}$ the degeneracy locus inside $\widetilde{Y}$ where $s$ drops rank by $2$, then the cohomology class of $\widetilde{X}$ is:
\[ \mathbbm{1}_{\widetilde{X}}=\begin{array}{|cc|}
c_2 & c_3 \\
1 & c_1
\end{array}\]
where $c_i=c_i(E^\vee \otimes \sqrt L^\vee)$. 
\end{thm}

The squaring principle guarantees that we are allowed to perform computations involving $\sqrt L$, as there exists a variety $\sigma : \widehat{Y} \to \widetilde{Y}$ on which  $\sqrt L$ is a line bundle and the induced map $\sigma^*: H^*(\widetilde{Y},\ZZ) \to H^*(\widehat{Y}, \ZZ)$ is injective. 

Recall that the matrix $A$ represents a map  \[ s: E\otimes L=\bigoplus\limits_{i=1}^5 L_i \otimes L\to E^\vee=\bigoplus\limits_{i=1}^5 L_i^\vee \] 
and the variety $\widetilde{X}$ is the degeneracy locus of the antisymmetric homomorphism of vector bundles $s$ on $\widetilde{Y}$. As all vector bundles involved are split, the Chern classes $c_i$ are easy to determine in terms of $M, N$ and $P$. 

From \S \ref{PfK} we have $L=-\dfrac{\sum_{i=1}^5 L_i}{2}= \cO(-4,-3,-1)$, so that $E\otimes L=\cO(-2,-2,-1)\oplus\cO(-4,-1,0)\oplus\cO(-1,-2,-2)\oplus\cO(-3,-2,0)\oplus\cO(-2,-2,0)$. Using the adjunction formula $-K_{\widetilde{X}}=-K_{\widetilde{Y}}+L =\cO_{\widetilde{Y}}(1,1,1)=M+N+P$, we compute the degree as follows (using Macaulay2 \cite{M2} for the very last equality):
\begin{align*}
(-K_{\widetilde{X}})^3&= (M+P+N)^3\cdot \mathbbm{1}_{\widetilde{X}}=  (M+P+N)^3\cdot (c_1c_2-c_3)=\\
 &= (M+P+N)^3(\sum \limits_{i} \widehat{L}_i \cdot \sum\limits_{i<j}  \widehat{L}_i \widehat{L}_j- \sum\limits_{i<j<k} \widehat{L}_i \widehat{L}_j \widehat{L}_k)= \dfrac{564}{2213}N^6,
\end{align*} 
where  \begin{align*}
\widehat{L}_1&=L_1^\vee\otimes \sqrt{L^\vee}=(-2M-N)+\frac{1}{2}(4M+3N+P)=\frac{1}{2}(N+P), \\
\widehat{L}_2&=L_2^\vee\otimes \sqrt{L^\vee}=\frac{1}{2} (4M-N-P) \\
\widehat{L}_3&=L_3^\vee\otimes \sqrt{L^\vee}=\frac{1}{2} (-2M+N+3P) \\
\widehat{L}_4&=L_4^\vee\otimes \sqrt{L^\vee}=\frac{1}{2} (2M+N-P) \\
\widehat{L}_5&=L_5^\vee\otimes \sqrt{L^\vee}=\frac{1}{2} (N-P). 
\end{align*}

Finally, we look at any chart, for example $U_{678}$, and determine that $\prod\limits_{i\notin\{6,7,8\}} D_i=\dfrac{36}{2213}N^6$. Since $\prod\limits_{i\notin\{6,7,8\}} D_i=1$ (the chart is smooth), we have that $(-K_X)^3=\frac{564}{36}=\frac{47}{3}$.


\newpage
\pagestyle{empty}
\newgeometry{top=1cm,bottom=1cm,left=1cm, right=1cm}
\begin{landscape}
\title{\textbf{Table \ref{tab2}: Exploring the GRDB}}
\centering
\maketitle
\arrayrulecolor{Gray}
\setlength{\arrayrulewidth}{0.2mm}
\rowcolors{1}{Gray!10}{Periwinkle!20}
\begin{center}
\begin{longtable}[t]{| p{0.3cm} | c | c |  >{\centering\arraybackslash}m{2.2cm} |  >{\centering\arraybackslash}m{0.5cm} |  >{\centering\arraybackslash}m{3.9cm} | c |  >{\centering\arraybackslash}m{3cm}| >{\centering\arraybackslash}m{1.6cm} |}
\hline
\label{tab2}
  &  \textbf{$\QQ$-Fano} &\textbf{Polytope} & \textbf{Codimension} \newline \textbf{and Singularities} & \textbf{MG \#} & \textbf{Rigid MMLP} & \textbf{Weight Matrix and Bundles} & \textbf{Period Sequence} & \textbf{Shape}\\
 \hline
\hline
1 & Q29792 & P542987  & codim. 13; \newline $\frac{1}{2}(1,1,1)$, $2\times\frac{1}{3}(1,1,2)$,  $\frac{1}{5}(1,1,4)$ & 28 & $x^2y^3z^4 + xy^3z^4 + 2xy^2z^2 + x + y^2z^3 + 2y^2z^2 + y + z + 
2yz/x + y/x + z/x + 3/(xyz) + 1/(x^2z) + 1/(x^2y) + 2/(x^2yz) + 
1/(x^2y^2z^3) + 2/(x^3y^2z^2) + 1/(x^3y^2z^3) + 1/(x^4y^3z^4)$ & \begin{math}\begin{array}{cccccccc | cc} 1&0&0&1&1&-3&2&2&2&1 \\ 0&1&0&1&1&1&1&-1&2&1 \\ 0&0&1&0&0&1&0&0&0&1 \end{array}\end{math} & [1, 0, 0, 0, 576, 2820, 4320, 2100, 2540160, 26490240, 110508300]& $\PP^2 \times \PP^1$\\

\hline
2 & Q29911 & P542990 &  codim. 13; \newline $3 \times \frac{1}{3}(1,1,2)$, $\frac{1}{4}(1,1,3)$  & 5 & $xy/z^2+ xy/z^3+xz+x/z^2+x/z^3+xz/y+yz+y/z^2+y/z^3+3z+z/y+yz/x+z/x$ & \begin{math}\begin{array}{cccccccc|c} 0&0&0&3&-2&1&-1&1&1 \\ 1&0&0&-1&1&1&0&0&1 \\0&1&0&-1&1&0&1&0&1 \\ 0&0&1&-1&1&0&0&1&1\end{array}\end{math} & [1, 0, 0, 54, 576, 0, 14850, 408660, 2540160, 5821200, 268077600] & $\PP^2 \times \PP^1$\\

\hline
3  & Q31617 & P543084& codim. 11; \newline $2 \times \frac{1}{2}(1,1,1)$, $\frac{1}{5}(1,2,3)$  & 19 & $x^2z + x^2z/y + xyz + 2xz + xz/y + yz + y/z + y/z^2 + z + 1/z + 2/z^2 + 
1/(yz^2) + y/(xz) + y/(xz^2) + 1/(xz) + 2/(xz^2) + 1/(xyz^2)$ & \begin{math} \begin{array}{cccccccc|cc} 0&0&2&-1&1&0&-1&1&1&0 \\ 1&0&-1&1&0&1&1&1&1&2 \\ 0&1&-1&1&-1&2&2&0&1&2 \end{array}\end{math} & [1, 0, 8, 72, 264, 6900, 44360, 552720, 7081480, 63161280, 887727708] & $\FF_1 \times \PP^1$\\

%

\hline

\hypertarget{3Dfour}{4} & Q33018 & P31470 & codim. 10; \newline $2\times\frac{1}{2}(1,1,1)$, $\frac{1}{3}(1,1,2)$ & 4 & $xy^2 + xy^2/z + 2xy + 4xy/z + xy/z^2 + x + 4x/z + 2x/z^2 + x/(yz) + 
x/(yz^2) + y^2/z + 2y/z + z + 1/z + 2/y + 1/(yz) + 1/(y^2z) + 1/x + 1/(xy)$ & \begin{math} \begin{array}{ccccccccc|ccccc}1&0&0&1&0&-1&3&2&-1&2&3&2&1&0 \\ 0&1&0&0&1&1&0&0&1&1&1&1&1&2 \\ 0&0&1&-1&1&2&-1&-1&1&1&-1&0&1&1 \end{array}\end{math} & [1, 0, 8, 108, 720, 8760, 111500, 1218000, 15156400, 193520880, 2416466808]
& \hyperlink{Z2}{$Z_2$}\\

\hline

\hypertarget{3Dfive}{5} & Q33019 & P542859 & codim. 11; \newline $\frac{1}{2}(1,1,1)$,  $\frac{1}{5}(1,2,3)$  & 6 & $x^3y/z^3 + x^2y^2/z^3 + x^2y/z^2 + x + 2x/z + y + 2y/z + z + y^2/(xz) + 
y/x + z^2/(xy) + z/(xy) + z/x^2$ & \begin{math}  \begin{array}{cccccccc | cc }    1&0 &0 &1 &1 &-2 &2 &1 & 1& 2 \\ 
  			        0 &1& 0 & 0& 1& 1 &0& 0& 1& 1 \\
		                 0& 0& 1& -1& 0& 1& 1& 1& 1& 1
		                  \end{array}  \end{math} 
& [1,  0, 0, 72, 432, 1020, 33840, 406140, 2222640, 28113120, 412416900] & \hyperlink{Z3}{$Z_3$} \\

 \hline
\hypertarget{3Dsix}{6} & Q33189& P543701 & codim. 11;  \newline  $3\times\frac{1}{2}(1,1,1)$, $\frac{1}{3}(1,1,2)$ & 6 & $xy/z^2 + xz + x + x/z^2 + xz/y + x/y + y/z^2 + 2z + 2z/y + z/y^2+ z/x+ z/(xy)$  & \begin{math} \begin{array}{cccccccc|cc} 1&0&2&-2&1&1&0&1&2&1 \\ 0&0&0&1&1&0&1&0&1&1 \\0&1&-1&1&0&-1&2&1&1&1  \end{array}\end{math}  & [1, 0, 0, 72, 192, 120, 32760, 176400, 342720, 21354480, 171397800]& \hyperlink{Z3}{$Z_3$} \\

\hline 

7& Q33470  &  P518676 & codim. 13; \newline  $3\times\frac{1}{2}(1,1,1)$, $\frac{1}{4}(1,1,3)$& 3 & $xyz^2 + x/z + yz^2 + y/z + z + 3/z + 1/(yz) + y/(xz) + 2/(xz) + 1/(xyz)$ & \begin{math}\begin{array}{cccccccc | cc} 1& 0& 0& -2& 1 &0 &1 &-1 & 1 & 0\\ 0 & 1 &0& 1& 1& 0& 0& 0& 1 & 0 \\ 0& 0& 1 &1 &0 &1 &1 &1 & 1 & 2 
 \end{array}
 \end{math}  
 &[1, 0, 6, 54, 90, 2880, 16710, 119700, 1679370, 10402560, 115479756]  & $\PP^1\times\PP^1$ + 0-st.\\
 
 \hline
8 & Q34466 & P515086 & codim. 12;  \newline $3 \times \frac{1}{2}(1,1,1)$, $\frac{1}{3}(1,1,2)$& 5 & $xyz^3 + xyz^2 + xz^3 + xz + x + xz/y + yz^3 + yz + y + 2z + 1/z + 
1/(yz) + yz/x + 1/(xz)$ & 
\begin{math}\begin{array}{cccccccccc|ccccc} 
1&0&0&0&-2&1&1&0&-1&0&1&0&0&0&1 \\
0&1&0&0&-1&1&1&0&-1&0&1&0&0&0&1 \\ 
0&0&1&0&1&1&0&0&1&1&1&1&1&2&1 \\
 0&0&0&1&1&0&1&1&1&0&1&2&1&1&1 \end{array}\end{math} & [1, 0, 8, 42, 264, 3360, 21590, 245280, 2202760, 20993280, 218056608]
 &  dP$_7$~+~0-st
\\

9 &&  & codim. 12; \newline $3 \times \frac{1}{2}(1,1,1)$, $\frac{1}{3}(1,1,2)$ & 2 & $x + x/z + x/(yz) + y + y/z + y/z^3 + z + 3/z + z/y + y^2/(xz^3) + 2y/(xz) + 
z/x$ & \begin{math}\begin{array}{ccccccccc|cc} 1&0&0&0&-2&1&-1&0&1&0&1 \\ 0&1&0&0&-1&1&-1&0&1&0&1 \\ 0&0&1&0&1&0&0&0&1&0&1 \\ 0&0&0&1&1&1&1&1&0&2&1\end{array}\end{math} & [1, 0, 10, 54, 342, 4680, 31150, 371700, 3524710, 34957440, 383835060] &  $\PP^1 \times \PP^1$ + 0-st
\\
\hline

10 & Q35361 & P516572  & codim. 12; \newline $2 \times \frac{1}{2}(1,1,1)$, $\frac{1}{3}(1,1,2)$& 9 & $x + x/(y^2z) + y + z + 1/y + 3/(yz) + 1/(yz^2) + yz/x + y/x + 1/x + 2/(xz) 
+ 1/(xz^2)$ & \begin{math} \begin{array}{cccccccc|cc} 1&0&-1&1&-1&2&1&-2&1&-1 \\ 0&1&-1&1&0&1&2&-1&1&1\\ 0&0&1&0&1&0&0&1&1&1  \end{array}\end{math} & [1, 0, 4, 60, 204, 2040, 24340, 160440, 1547980, 16366560, 139003704] & $\FF_1 \times \PP^1$ \\
\hline

\hypertarget{3Deleven}{11} & Q35364 & P516506 & codim. 13; \newline $2 \times\frac{1}{2}(1,1,1)$, $\frac{1}{4}(1,1,3)$ & 7 & $xy/z^2 + xy/z^3 + xz + x + x/z^2 + xz/y + yz + y + y/z^2 + 3z + z/y +yz/x + z/x$ & \begin{math}\begin{array}{cccccccc|cc} 1&0&0&-2&1&1&2&1&1&2\\ 0&1&0&1&1&-1&1&1&1&2\\ 0&0&1&1&0&0&0&0&1&0\end{array}\end{math}& [1, 0, 0, 54, 168, 120, 14850, 100380, 239400, 6063120, 61021800] & \hyperlink{Z4}{$Z_4$} \\

\hline
12 & Q35470  & P519948  & codim. 19; \newline $2 \times \frac{1}{4}(1,1,3)$, $\frac{1}{6}(1,1,5)$ & 10 & $x + y^3z^4 + 2y^2z^2 + y + y^2z^3/x^3 + 3yz^2/x^3 + 2yz/x^3 + 3z/x^3 + 
6/x^3 + 1/(x^3z) + 1/(x^3y) + 6/(x^3yz) + 3/(x^3yz^2) + 2/(x^3y^2z^2) + 
3/(x^3y^2z^3) + 1/(x^3y^3z^4)
$ & \begin{math} \begin{array}{c c c c c c | c} 1&0&2&-2&1&1&2 \\ 0&1&1&2&1&1&2 \end{array} \end{math} & [1, 0, 0, 0, 24, 120, 0, 0, 2520, 75600, 88200]
& $\PP^1 \times \PP^1$ + 0-st\\

\hline

13 & Q35496 & P518648 & codim. 13; \newline $3 \times \frac{1}{2}(1,1,1)$, $\frac{1}{3}(1,1,2)$ & 6 & $xyz^3 + xyz^2 + xz^3 + xz + x + xz/y + yz^3 + yz + y + 3z + 1/(yz) + 
yz/x + 1/(xz)$ & \begin{math} \begin{array}{c c c c c c c c | c c} 0&0&2&-1&2&-1&1&1&1&2\\ 1&0&-1&1&1&0&1&1&1&2\\ 0&1&-1&1&0&1&0&0&1&0 \end{array}\end{math} & [1, 0, 4, 54, 36, 1800, 15250, 44100, 992740, 6773760, 40547304]  & $\FF_1 \times \PP^1$ \\

\hline
 
14 & Q36639 &  P254485 &  codim. 15; \newline $6\times\frac{1}{2}(1,1,1)$ & 1 & $xy^2z^3 + xy^2z^2 + xyz^2 + 2xyz + xy + x + x/z + y + 1/y + 1/(xyz)$ &
 \begin{math}\begin{array}{ccccccc|c} 1&0&0&0&1&-1&1&-1\\0&1&0&0&-1&1&-1&1\\0&0&1&0&1&0&0&1\\0&0&0&1&0&1&1&0\end{array}\end{math}  & [1, 0, 6, 0, 114, 120, 3300, 7560, 114450, 393120, 4461156]
 & $\PP^1\times \PP^1$ + 0-st. \\

15 & & & codim. 15; \newline $6\times\frac{1}{2}(1,1,1)$ & 1 & $xy^2z^3 + xy^2z^2 + xyz^2 + 3xyz + xy + x + x/z + y + 1/y + 1/(xyz)$ & 
\begin{math}
 \begin{array}{ccccccc|c} 
1&0&0&0&-1&1&1&-2 \\ 
0&1&0&0&1&-1&-1&2\\ 
0&0&1&0&1&0&1&1\\ 
0&0&0&1&0&1&0&0 
\end{array}\end{math} & [1, 0, 8, 0, 168, 120, 5120, 9240, 190120, 559440, 7976808]  & $\FF_1$ + 0-st. \\ 

 \hline 
\hypertarget{3Dseventeen}{16}  & Q37383 & P516493  & codim. 13; \newline $\frac{1}{3}(1,1,2)$  & 4 & $xy/z + xz + yz + y + y/z + 2z + 1/z+ z/y+yz/x+y/x+2z/x+1/x+z/(xy)$ &  \begin{math}\begin{array}{cccccccc|cc} 1&0&1&-1&0&1&1&1&1&2\\ 0&0&0&1&1&0&1&-1&1&0\\0&1&-1&1&0&1&0&0&1&0 \end{array}\end{math} & [1, 0, 8, 18, 264, 1200, 12950, 82320, 778120, 5856480, 52669008]& \hyperlink{Z3}{$Z_3$} \\

\hline

17  & 37832  & P220586 & codim. 13; \newline $\frac{1}{2}(1,1,1)$ & 2 & $x + x/(y^2z) + y + z + 1/y + 1/(yz) + 2/(y^2z) + y/x + 2/x + 1/(xy) + 
1/(xyz) + 1/(xy^2z) $ & \begin{math} \begin{array}{c c c c c c c | c} 1&0&0&0&1&0&1&2 \\ 0&1&0&0&0&1&-1&0 \\ 0&0&1&1&-1&2&-1&0 \end{array} \end{math} & [1, 0, 6, 18, 186, 1200, 9330, 69720, 553770, 4450320, 36502956] & $\FF_1~+~0-st$\\

\hline
18  & Q37847& P515113& codim. 15; \newline  $\frac{1}{2}(1,1,1)$, $\frac{1}{3}(1,1,2)$ & 3  &  $x + x/(yz) + y + z + z/y + 2yz/x + 2z/x + 2y^2z/x^2 + yz/x^2 + 
y^2z^2/x^3 + y^3z^2/x^4$ & \begin{math}\begin{array}{ccccccccc|ccccc} 1&0&0&0&0&1&2&1&-1&2&2&1&1&0\\ 0&1&0&1&0&0&1&1&1&1&1&2&1&1\\ 0&0&1&1&1&1&0&0&1&1&1&1&1&2 \end{array}\end{math} & [1, 0, 6, 18, 90, 960, 3210, 39900, 197610, 1668240, 12109356]  &  dP$_7$~+~0-st\\

19 & &  & codim. 15; \newline  $\frac{1}{2}(1,1,1)$, $\frac{1}{3}(1,1,2)$ & 3 & $x + x/(yz) + y + z + z/y + 3yz/x + 2z/x + 2y^2z/x^2 + yz/x^2 + 
y^2z^2/x^3 + y^3z^2/x^4$ & \begin{math}\begin{array}{cccccccc|ccc} 1&0&0&0&2&0&-1&1&2&0\\ 0&1&0&1&1&0&0&0&1&0 \\ 0&0&1&1&0&1&1&1&1&2 \end{array}\end{math} & [1, 0, 4, 18, 60, 780, 2470, 27720, 151900, 1063440, 8348004] &  $\PP^1 \times \PP^1$ + 0-st\\
\hline

20 &  & P516442& codim. 15 \newline $\frac{1}{2}(1,1,1)$, $\frac{1}{3}(1,1,2)$ & 6 & $ x^2z/y + x^2z/y^2 + xy + x + 3xz/y + xz/y^2 + y + z + 2z/y + y/(xz) + 
z/x$ & \begin{math}\begin{array}{cccccccc|cc} 1&0&0&0&2&-1&1&1&1&2 \\ 0&1&0&1&1&0&0&1&1&1\\ 0&0&1&1&0&1&1&0&1&1 \end{array}\end{math} & [1, 0, 4, 18, 60, 720, 2470, 25200, 141820, 972720, 7447104]  & $\FF_1$ + 0-st.
\\ 

21 &   & P516475 & codim. 15; \newline $\frac{1}{2}(1,1,1)$, $\frac{1}{3}(1,1,2)$ & 11 & $xy^3z^2 + 2xy^2z^2 + xyz^2 + x + 2y^2z + 4yz + y + 3z + z/y + y/x + 
2/x + 1/(xz) + 1/(xy)$ & \begin{math}\begin{array}{ccccccccc|ccccc} 1&0&0&0&0&1&1&0&-1& 1&1&0&0&0 \\ 0&1&0&1&1&0&0&1&2& 1&1&2&2&2 \\ 0&0&1&1&2&1&0&0&1 &2&1&1&2&2 \end{array}\end{math} & [1, 0, 4, 24, 60, 780, 4000, 24360, 209020, 1228920, 9356004] & dP$_7$+0-st.
\\ 
 
\hline
 
22 & Q37859 & P424999 & codim. 16; \newline $2\times \frac{1}{3}(1,1,2)$ & 15 & $x + x/(yz) + y + z + 2/z + 1/(yz) + y^2/(xz) + 3y/(xz) + 1/(xz) + 
2y^3/(x^2z) + 2y^2/(x^2z) + y^4/(x^3z) $ & \begin{math}\begin{array}{ccccccccc|ccccc} 1&0&0&0&-1&1&2&1&-1& 1&2&1&0&0 \\ 0&1&0&1&0&0&1&1&0& 1&1&1&1&0 \\ 0&0&1&1&2&1&0&0&1&2&1&1&2&2 \end{array}\end{math} & [1, 0, 4, 12, 60, 480, 1660, 17220, 74620, 610680, 3641904] & dP$_7$~+~0-st. \\
\hline 
23 & Q37894& P425006 & codim. 16; \newline $2\times \frac{1}{2}(1,1,1)$, $\frac{1}{3}(1,1,2)$ & 14 & $x + y + z + y/(xz) + 2/x + z/(xy) + 1/(xy) + z/(xy^2) + 1/(x^2z) + 
1/(x^2y) + 1/(x^2yz) + 2/(x^2y^2) + z/(x^2y^3)$  & \begin{math}\begin{array}{ccccccccc|ccccc} 1&0&0&0&2&1&-1&-1&1 &1&0&0&1&2\\ 0&1&0&1&1&0&0&1&1&1&1&1&2&1\\ 0&0&1&1&0&1&1&1&0&1&2&1&1&1  \end{array}\end{math} & [1, 0, 4, 6, 60, 360, 1210, 14280, 41020, 505680, 1966104] & dP$_7$~+~0-st. 
  \\

24 & & & codim. 16; \newline $2\times \frac{1}{2}(1,1,1)$, $\frac{1}{3}(1,1,2)$ & 1 & $x^3/(yz^2) + x^2/(yz) + x + 2x/(yz) + y + z + 3/y + z/(xy) + 1/(xy) + 
z/(x^2y)$ & \begin{math}\begin{array}{cccccccc|cc} 1&0&0&0&2&-1&-1&1&1&0\\ 0&1&0&1&1&0&0&0&1&0 \\ 0&0&1&1&0&1&1&1&1&2 \end{array}\end{math} & [1, 0, 6, 6, 90, 420, 1950, 19740, 63210, 810600, 2968056]
&   $\PP^1 \times \PP^1$  + 0-st. 
  \\

\hline
25 & Q37905 & P515878 & codim. 17; \newline  $\frac{1}{2}(1,1,1), \newline  2\times\frac{1}{3}(1,1,1)$& 2 & $x + xz^2/y + y + z + z/y + 2/x + 2/(x^2z) + y/(x^3z^2) + y/(x^4z^3)$ & \begin{math}\begin{array}{ccccccc|c} 1&0&0&2&-2&-1&3&2 \\ 0&1&0&1&0&0&1&1\\ 0&0&1&0&1&1&0&1  \end{array}\end{math} & [1, 0, 4, 0, 60, 120, 1120, 6300, 24220, 249480, 655704]
 &  $\FF_1$ + 0-st
 \\

26 & & & codim. 17; \newline  $\frac{1}{2}(1,1,1), \newline  2\times\frac{1}{3}(1,1,1)$ & 2 & $x + xz^2/y + y + z + z/y + 3/x + 2/(x^2z) + y/(x^3z^2) + y/(x^4z^3)$ & \begin{math}\begin{array}{ccccccc|c} 1&0&0&2&-2&-1&1&0 \\ 0&1&0&1&0&0&0&0\\ 0&0&1&0&1&1&1&2 \end{array}\end{math} & [1, 0, 6, 0, 90, 120, 1860, 7980, 44730, 378000, 1260756]  &  $\PP^1 \times \PP^1$  + 0-st.
\\

\hline
 
27& Q38250 &  P254482 &  codim. 16; \newline $3\times\frac{1}{2}(1,1,1)$ & 1 & $x^2z^3/y^2 + x^2z^4/y^3 + xz + x + 2xz^2/y + xz^3/y^2 + y + z + y/(xz) + y/(xz^2)$ &
\begin{math}
\begin{array}{cccccccc|c}
1&0&0&0&1&-1&0&0&0\\ 0&1&0&0&-1&1&1&-1&0\\0&0&1&0&1&0&1&0&1\\0&0&0&1&0&1&0&1&1 
\end{array}
\end{math}
& [1, 0, 4, 6, 84, 300, 2290, 11760, 80500, 483000, 3182004]
& $\PP^1\times \PP^1$ + 0-st.

\\
\hline
28 & & & codim. 16; \newline $3\times\frac{1}{2}(1,1,1)$ & 1 & $x^2z^3/y^2 + x^2z^4/y^3 + xz + x + 3xz^2/y + xz^3/y^2 + y + z + y/(xz) + y/(xz^2)$ 
 
 &
 \begin{math}\begin{array}{cccccccc|c} 
 1&0&0&0&0&0&1&-1&0\\
  0&1&0&0&-1&1&-2&1&-1\\
  0&0&1&0&1&0&1&1&2\\
  0&0&0&1&0&1&0&0&0  
  \end{array}\end{math} & [1, 0, 6, 6, 114, 360, 3390, 16800, 126210, 765240, 5368356] &
 $\FF_1$ + 0-st \\

\hline
29 & Q38484 & P515100 & codim. 15; \newline $\frac{1}{2}(1,1,1)$  & 14 &  $x + y + z + z/y + y/x + y/(xz) + 2/x + 1/(xz) + z/(xy) + 2/(xy) + z/(xy^2)$ & \begin{math}\begin{array}{ccccccccc|ccccc} 1&0&0&0&0&1&1&1&0&1&2&1&1&1 \\ 0&1&0&1&0&0&1&1&1&1&1&2&1&1 \\ 0&0&1&1&1&1&0&0&1&1&1&1&1&2 \end{array}\end{math} &
[1, 0, 4, 24, 84, 780, 4360, 29400, 214900, 1433040, 10540404] & dP$_7$~+~0-st.
 \\

30 &  & & codim. 15; \newline $\frac{1}{2}(1,1,1)$ & 4 & $x + y + z + z/y + y/x + y/(xz) + 3/x + 1/(xz) + z/(xy) + 2/(xy) + z/(xy^2)$
&  \begin{math}\begin{array}{cccccccc|cc}  1&0&0&0&1&0&0&1&1&1\\0&1&0&1&1&0&0&0&1&0\\ 0&0&1&1&0&1&1&1&1&2 \end{array}\end{math} & [1, 0, 6, 24, 114, 1020, 5460, 42000, 299250, 2100840, 16160256] &  $\PP^1 \times \PP^1$ + 0-st
 
\\
\hline

 31 & Q38496 & P427129 & codim. 16; \newline $2 \times  \frac{1}{2}(1,1,1)$ & 1 & $x + y + z + y/x + y/(xz) + 2/x + 1/(xz) + z/(xy) + 2/(xy) + z/(xy^2)$
 & \begin{math} \begin{array}{ccccccc|c} 1&0&0&1&-1&1&1&2 \\ 0&1&0&1&0&1&0&1\\ 0&0&1&0&1&0&1&1 \end{array}\end{math} & [1, 0, 4, 18, 60, 600, 2470, 18900, 118300, 723240, 5242104] & $\FF_1$ + 0-st \\
 
 32 & & & codim. 16; \newline $2 \times\frac{1}{2}(1,1,1)$   &  1  & $x + y + z + y/x + y/(xz) + 3/x + 1/(xz) + z/(xy) + 2/(xy) + z/(xy^2)$
  & \begin{math}\begin{array}{ccccccc|c} 1&0&0&1&-1&0&1&1 \\ 0&1&0&1&0&0&0&0\\ 0&0&1&0&1&1&1&2 \end{array} \end{math} &  [1, 0, 6, 18, 90, 780, 3210, 28560, 164010, 1146600, 8247456]  &  $\PP^1 \times \PP^1$ + 0-st
 \\
 
 \hline

33  & Q38895 & P425118 & codim. 17; \newline $\frac{1}{2}(1,1,1)$ & 1 &  $x^2z^3/y^2 + xz + x + 3xz^2/y + xz^3/y^2 + y + y/z + 2z + z^2/y + y/(xz^2)
$ & \begin{math}\begin{array}{cccccccc|cc} 
1&0&0&0&1&0&-1&1&1&0\\ 
0&1&0&1&1&0&0&0&1&0\\
 0&0&1&1&0&1&1&1&1&2 \end{array}\end{math} &  [1, 0, 6, 12, 90, 540, 2400, 20160, 95130, 751800, 4291056] & $\PP^1 \times \PP^1$ + 0-st \\
 
 34 &&  & codim. 17; \newline $\frac{1}{2}(1,1,1)$ & 4 & $x + yz + y + z + y/x + y/(xz) + z/x + 2/x + 1/(xz) + z/(xy) + 1/(xy)$ & \begin{math}\begin{array}{cccccccc|cc} 1&0&0&0&1&0&0&1&1&1\\ 0&1&0&1&1&0&1&0&1&1\\ 0&0&1&1&0&1&0&1&1&1\end{array}\end{math} &  [1, 0, 4, 12, 60, 420, 1660, 13440, 64540, 451920, 2665404] & $\PP^1 \times \PP^1$ + 0-st
 \\
 
\hline
35 &  Q38906 & P516432 & codim. 18;  \newline  $2\times \frac{1}{2}(1,1,1)$& 2 & $x + xz^2/y + y + z + 2z/y + 2/x + 1/(xy) + y/(x^2z) + 1/(x^2z)$ & \begin{math}\begin{array}{ccccccc|c} 1&0&0&1&-1&1&0 & 1\\0&1&0&1&0&0&1&1\\0&0&1&0&1&1&0 &1  \end{array}\end{math} & [1, 0, 4, 6, 60, 240, 1210, 7980, 34300, 256200, 1172304]  &  $\FF_1$ + 0-st \\

36 & & & codim. 18;  \newline  $2\times \frac{1}{2}(1,1,1)$& 2 &  $x + xz^2/y + y + z + 2z/y + 3/x + 1/(xy) + y/(x^2z) + 1/(x^2z)$ & \begin{math}\begin{array}{ccccccc|c} 
  1&0&0&1&-1&-1&1&0\\ 0&1&0&1&0&0&0&0\\ 0&0&1&0&1&1&1&2 \end{array}\end{math} & [1, 0, 6, 6, 90, 300, 1950, 11760, 56490, 432600, 2023056] &  $\PP^1 \times \PP^1$   + 0-st 
 \\
  
\hline

37 & Q38935  & P424987  & codim. 19; \newline $3\times \frac{1}{2}(1,1,1)$& 1 & $x^2z^3/y^2 + 2x^2z^4/y^3 + x^2z^5/y^4 + xz + x + 2xz^2/y + xz/y + xz^3/y^2 + y + y/(xz^2) $ & \begin{math}\begin{array}{ccccccc|c} 1&0&0&-1&1&2&-1&1\\ 0&1&0&1&0&0&0&0\\0&0&1&0&1&1&1&2 \end{array}\end{math} & [1, 0, 4, 0, 60, 60, 1120, 2520, 24220, 90720, 586404]
 & $\PP^1\times \PP^1$ + 0-st. 
\\

38 & & &   codim. 19; \newline $3\times \frac{1}{2}(1,1,1)$ & 1 & $x^2z^3/y^2 + 2x^2z^4/y^3 + x^2z^5/y^4 + xz + x + 3xz^2/y + xz/y + xz^3/y^2 + y + y/(xz^2)$ & \begin{math}\begin{array}{ccccccc|c} 1&0&0&-1&1&-1&1&0 \\ 0&1&0&1&0&0&1&1\\ 0&0&1&0&1&1&0&1 \end{array}\end{math} & [1, 0, 6, 0, 90, 60, 1860, 3360, 44730, 143640, 1191456] & $\FF_1$ + 0-st. \\

\hline

39 & \hypertarget{40948}{Q40948} & P534760 &  codim 9; \newline $2\times \frac{1}{5}(1,2,4)$ &  5 & $xy^3z^5 + x + 2y^2z^5 + 2yz^2 + y + yz^5/x + 2z^2/x + 1/x + 1/(xyz)$ &  \begin{math} \begin{array}{cccccc|c} 1&-3&-5&-1&-2&2&-6 \\ 0&3&5&1&2&-1&6 \end{array} \end{math} & [1, 0, 2, 0, 6, 0, 380, 0, 10150, 0, 151452] & $\FF_1$ (no~0-st.) \\

\hline

40 & \hypertarget{40971}{Q40971} & P543951 & codim 9; \newline $2\times \frac{1}{5}(2,2,3)$ & 2 & $x^6/(y^4z^5) + 2x^3/(y^2z^2) +2x^2/(y^3z^2) + x + y + z + 1/y + 2z/(xy) + z/(x^2y^2)$ & \begin{math}\begin{array}{cccccc|c} 1&0&0&0&0&1&0 \\ 0&3&1&5&2&3&6  \end{array}\end{math} & [1, 0, 2, 0, 6, 0, 20, 0, 3430, 0, 75852] & $\FF_1$ (no~0-st.) \\

\hline

41 & \hypertarget{40988}{Q40988}  & P413267 &  codim 6; \newline $\frac{1}{5}(1,2,4)$ & 1 & $ x + y + z + 1/(yz^2) + 1/(y^2z^3) + 1/(xz^2) + 1/(xyz) + 1/(xyz^3) + 1/(xy^2z^2)$ & \begin{math}\begin{array}{cccccc|c} 1&0&1&-1&1&2&2\\ 0&1&1&2&-1&1&2 \end{array} \end{math}  &  [1, 0, 0, 0, 48, 0, 360, 0, 11760, 0, 226800] & $\PP^2$ (no~0-st)\\

\hline

42 & \hypertarget{40993}{Q40993}  & P473887 & codim 6; \newline $3 \times \frac{1}{3}(1,2,2)$ & 2 & $x + xz/y^3 + y + z + 2z/y^2 + 1/x + z/(xy) + 2/(xy^2) + 2/(x^2y) + 1/(x^3yz)$ & \begin{math}\begin{array}{ccccccc|c} 1&0&0&1&-1&-1&0&0 \\ 0&1&0&0&0&0&1&0\\ 0&0&1&-1&3&2&1&2 \end{array}\end{math} & [1, 0, 2, 0, 54, 0, 740, 0, 18550, 0, 403452] & $\FF_1$ + 0-st.\\


 \hline 
43  & \hypertarget{41200}{Q41200} & P547328 & codim. 4; \newline $2\times \frac{1}{2}(1,1,1)$, $\frac{1}{5}(1,3,4)$ & 6 & $xyz^2 + x + xz/y + x/(y^2z) + y + z + 2/(yz) + 2/(xz) + 1/(xyz^3) + 
1/(x^2z^3)$ & \begin{math}\begin{array}{ccccccc|cc} 2&1&3&1&2&1&0& 3&4 \\ 1&-1&-1&1&-1&0&1 &0&0\end{array}\end{math} & [1, 0, 0, 72, 0, 0, 32760, 0, 0, 21067200, 0]
 & $\FF_1$ (no~0-st.)
\\ 

\hline 
44 & \hypertarget{41218}{Q41218} & P544064  & codim. 4; \newline $3 \times \frac{1}{2}(1,1,1)$, $\frac{1}{4}(1,1,3)$& 3 & $x + x/(y^2z) + y + z + 1/(yz) + 3yz/x + 2y^2z/x^2 + y^4z^3/x^3 + y^5z^3/x^4$ & \begin{math}\begin{array}{ccccccc|cc} 1&0&3&0&1&-1&2& 3&0 \\ 0&1&-1&1&1&1&-1&0&2 \end{array}\end{math} & [1, 0, 0, 54, 0, 0, 14850, 0, 0, 5821200, 0]  & $\PP^1 \times \PP^1$ (no~0-st.)\\

\hline
45 & \hypertarget{41251}{Q41251} & P402202 & codim 5; \newline $2\times\frac{1}{2}(1,1,1)$, $\frac{1}{4}(1,3,3)$ & 5 & $xyz^2 + x + y + z + z^2/y + 1/(xz) + 1/(xy)$ &  \begin{math}\begin{array}{cccccc|c} 1&0&2&-1&-1&4&2 \\ 0&1&-1&1&1&-2&0  \end{array} \end{math}  & [1, 0, 0, 12, 0, 0, 900, 0, 0, 86520, 0] & $\PP^2$ (no~0-st.)\\


\hline
46 & \hypertarget{41334}{Q41334}  & P543852 & codim. 4; \newline  $3\times \frac{1}{3}(1,1,2)$, $\frac{1}{5}(1,4,4)$ & 3 & $x^2/(y^3z^2) + x + x/(y^2z^2) + y + z + 3yz/x + 2y^2z/x^2 + y^5z^4/x^4 + 
y^6z^4/x^5$ & \begin{math}\begin{array}{ccccccc|cc} 1&0&-1&1&1&1&-1& 0&2 \\ 0&1&4&0&1&-1&3&4&0 \end{array} \end{math} & [1, 0, 0, 0, 576, 0, 0, 0, 2540160, 0, 0]
 & $\PP^1 \times \PP^1$ (no 0-st.)\\

\hline
\end{longtable}
\end{center}

\end{landscape}
\restoregeometry

\bibliographystyle{plain}
\bibliography{bibLI}

\end{document}